\documentclass[a4paper,11pt,final]{amsart}

\usepackage[utf8]{inputenc}
\usepackage[T1]{fontenc}
\usepackage[english]{babel}
\usepackage[margin=2cm]{geometry}
\usepackage{amsmath}
\usepackage{amsthm}
\usepackage{amssymb}
\usepackage{dsfont}
\usepackage{stmaryrd}
\usepackage{mathabx}
\usepackage{mathrsfs}
\usepackage[all]{xy}
\usepackage{wasysym}
\usepackage{mathtools}
\usepackage{bbold}
\usepackage{empheq}
\usepackage{tensor}
\usepackage{bbm}
\usepackage{enumitem}
\usepackage{rotating}
\usepackage{multirow}
\usepackage{setspace}
\usepackage{fancybox}
\usepackage{xargs}
\usepackage{indentfirst}
\usepackage{xcolor}
\definecolor{InternalLinks}{rgb}{0.33, 0.29, 0.31}
\usepackage{hyperref}
\usepackage{showkeys}
\usepackage[colorinlistoftodos,prependcaption,textsize=tiny,disable]{todonotes}
\usepackage{draftwatermark}
\SetWatermarkScale{5}
\SetWatermarkColor[gray]{0.9}
\usepackage{fixme}
\usepackage{comment} 
\newif\ifshow 
\showtrue 
\ifshow
  \includecomment{wrap}
\else
  \excludecomment{wrap} 
\fi

\theoremstyle{plain}
\newtheorem{theo}{Theorem}[section]
\newtheorem*{theo*}{Theorem}
\newtheorem{lemm}[theo]{Lemma}
\newtheorem{prop}[theo]{Proposition}
\newtheorem{coro}[theo]{Corollary}

\theoremstyle{definition}
\newtheorem{defi}[theo]{Definition}
\newtheorem{exam}[theo]{Example}

\theoremstyle{remark}
\newtheorem{rema}[theo]{Remark}
\newtheorem*{rema*}{Remark}
\newtheorem{nota}[theo]{Notation}
\newtheorem*{nota*}{Notation}

\numberwithin{equation}{section}

\newcommand{\pr}{\noindent{\it Proof}\quad}
\newcommand{\fin}{\qed\smallskip}

\newcommand{\ku}{\Bbbk}
\newcommand{\id}{\mathrm{id}}

\newcommand{\W}{\mathbb{W}}
\newcommand{\V}{\mathbb{V}}
\newcommand{\U}{\mathbb{U}}

\newcommand{\He}{\mathcal{H}}
\newcommand{\YD}{\mathcal{YD}}
\newcommand{\aYD}{\mathfrak{YD}}
\newcommand{\alg}{\mathfrak{Alg}}
\newcommand{\D}{\mathcal{D}}
\newcommand{\T}{\mathcal{T}}
\newcommand{\p}{\mathbf{p}}

\newcommand{\m}{\mathfrak{m}}
\newcommand{\du}[1]{\widehat{#1}}
\newcommand{\com}{\Delta}
\newcommand{\dcom}{\du{\Delta}}
\newcommand{\dS}{\du{S}}

\newcommand{\cou}{\varepsilon}
\newcommand{\dcou}{\du{\varepsilon}}

\newcommand{\op}{\textrm{op}}
\newcommand{\co}{\textrm{co}}
\newcommand{\ops}{\circ}

\newcommand{\alp}{\alpha}

\newcommand{\odo}{\otimes}

\newcommand{\oti}{\otimes}
\newcommand{\bicroi}{\bowtie}
\newcommand{\bt}{\boxtimes}

\newcommand{\blhd}{\blacktriangleleft}
\newcommand{\brhd}{\blacktriangleright}

\newcommand{\ome}{\omega}
\newcommand{\ad}{\mathrm{Ad}}

\newcommand{\sii}{\Leftrightarrow}

\newcommand{\ida}{\Rightarrow}
\newcommand{\vuelta}{\Leftarrow}

\newcommand{\fle}[1]{\xymatrix@C=3pc@R=3pc{#1}}

\newcommand{\iso}{\cong}

\newcommand\restr[2]{{
  \left.\kern-\nulldelimiterspace 
  #1 
  \vphantom{\big|} 
  \right|_{#2} 
}}

\definecolor{ToDo}{RGB}{30,144,255}
\definecolor{Question}{RGB}{220,20,60}
\definecolor{Attention}{RGB}{255,215,0}

\newcommandx{\td}[2][1=]{\hfill\todo[inline,size=\footnotesize,linecolor=ToDo,backgroundcolor=ToDo!40,bordercolor=black,#1]{{\bf ToDo:} #2}}
\newcommandx{\qs}[2][1=]{\hfill\todo[inline,size=\normalsize,linecolor=Question,backgroundcolor=Question!40,bordercolor=black,#1]{{\bf Question:} #2}}
\newcommandx{\rmq}[2][1=]{\hfill\todo[inline,size=\normalsize,linecolor=Attention,backgroundcolor=Attention!40,bordercolor=black,#1]{{\bf Remark:} #2}}

\makeatletter
\newcommand{\addresseshere}{\enddoc@text\let\enddoc@text\relax}
\makeatother
\makeatletter
\newcommand{\pushright}[1]{\ifmeasuring@#1\else\omit\hfill$\displaystyle#1$\fi\ignorespaces}
\newcommand{\pushleft}[1]{\ifmeasuring@#1\else\omit$\displaystyle#1$\hfill\fi\ignorespaces}
\makeatother

\title[Notes on Yetter--Drinfeld algebras over Hopf algebras]{Notes on Yetter--Drinfeld algebras \protect\\ over Hopf algebras}
\author{Frank Taipe}
\address{Universit\'e Paris-Saclay, CNRS, Laboratoire de Math\'ematiques d'Orsay, 91405 Orsay, France}
\email{frank.taipe@universite-paris-saclay.fr}


\hypersetup{
pdfauthor={Frank Taipe},
pdftitle={Notes on Yetter-Drinfeld algebras over Hopf algebras},
pdfsubject={Notes on Yetter-Drinfeld algebras over Hopf algebras},
pdfkeywords={Yetter-Drinfeld algebras, Hopf algebras}}

\begin{document}

\maketitle

\begin{abstract}
In this work, we study another characterization of Yetter--Drinfeld algebras over finite-dimensional Hopf algebras. We show the equivalence between this characterization, called the ``only coaction'' characterization, and the standard ``action-coaction'' characterization. This modern approach for Yetter--Drinfeld algebras is one of the key ingredient in a self-dual theory of quantum transformation groupoids arising from actions of quantum groups.
\end{abstract}

\tableofcontents

\section{Introduction}

Motivated for a modern approach of Yetter--Drinfeld structures over a pairing of Van Daele's multiplier Hopf algebras, in this work, we aim to explore this approach in the case of finite-dimensional Hopf algebras. Finite-dimensional Hopf algebras make more easy the computations because we don't have to worry about the {\em covering technique} as is done in the case of multiplier Hopf algebras.

\medskip

The plan of this work is the following: Section~\ref{sec:hopf_algebras} is intended to fix notations and collect all the basics facts about finite-dimensional Hopf algebras and their actions and coactions. In Section \ref{sec:standard}, we give a survey on the standard characterization of Yetter--Drinfeld algebras in the case of finite-dimensional Hopf algebras and we give proofs only for important results. Then, in Section~\ref{sec:only_coaction}, we give the ``only coaction'' characterization of Yetter--Drinfeld algebras and we prove the equivalence with the standard characterization. Here, we give results on {\em duality} of Yetter--Drinfeld algebras over finite-dimensional Hopf algebras. In Section~\ref{sec:yd_dd}, we give the equivalence between Yetter--Drinfeld algebras, module algebras over Drinfeld doubles, and comodule algebras over Drinfeld codoubles. Section~\ref{sec:equi_yd_cat} is intended to give the equivalences between Yetter--Drinfeld categories in the standard and in the only-coaction characterization form. Finally, in Section~\ref{sec:examples_yd}, we show some examples. In the appendix, we recollect some topics relate to Yetter--Drinfeld algebras. 

\medskip

In order to improve this work, comments and/or suggestions are more than welcome.

\subsection*{Conventions and Notations}

An algebra in this work means a unital $\ku$-algebra where $\ku$ is a field. Let $X, Y$ be two algebras. For any invertible element $W$ in $X \odo Y$, there is an invertible element $\ad(W) \in \mathrm{End}(X \odo Y)$ defined by
$\ad(W)(x \odo y) = W(x \odo y)W^{-1}$ for all $x \in X$, $y \in Y$. Its inverse map is given by $\ad(W^{-1})$. Let $X, Y, Z$ be algebras. Through this work, given an element $E \in X \odo Y$, we will made use of the {\em legs notation}
\[
E_{12} := E \odo 1_{Z}, \quad E_{23} := 1_{Z} \odo E, \quad E_{13} := (\id_{X} \odo \Sigma_{Y,Z})(E_{12}) = (\Sigma_{Z,X} \odo \id_{Y})(E_{23})
\]
see as elements in $X \odo Y \odo Z$, $Z \odo X \odo Y$ and $X \odo Z \odo Y$, respectively. For any algebra $X$, its multiplication map is denoted by $\m_{X}: X \odo X \to X$ and its opposite algebra is denoted by $X^{\op}$, i.e. the $\ku$-module $X$ with multiplication map $\m_{X^{\op}} = \bar{\m}_{X} := \m_{X}\Sigma_{X,X}$.

\section{Basics on finite-dimensional Hopf algebras}\label{sec:hopf_algebras}

Given a Hopf algebra $(H,\com,\cou,S)$, we will denoted by $(\du{H},\dcom,\dcou,\dS)$ its dual Hopf algebra.  The {\em canonical pairing of $H$} is the bilinear form $\p: H \times \du{H} \to \ku$, $(h,\theta) \mapsto \theta(h)$. The bilinear form $\bar{\p}: \du{H} \times H \to \ku$, $(\theta,h) \mapsto \p(h,\theta)$ is called the {\em flip canonical pairing} and it is the canonical pairing of $\du{H}$. It follows directly that the canonical pairing $\p$ satisfies the equalities
\[
\p \odo \p(h \odo h',\dcom(\theta)) = \p(hh',\theta), \quad \p \odo \p(\com(h), \theta \odo \theta') = \p(h,\theta\theta'),
\]
\[
\p(1_{{H}},\theta) = \dcou(\theta), \quad \p(h,1_{\du{H}}) = \cou(h), \quad \p(S(h),\theta) = \p(h,\dS(\theta)).
\]
for all $h,h' \in H$ and $\theta, \theta' \in \du{H}$.

\begin{prop}
Let $H$ and $H'$ be two Hopf algebras. Given an Hopf algebras isomorphism $T: H \to H'$, then there exist an unique  Hopf algebra isomorphism $\hat{T} : \du{H'} \to \du{H}$ such that $\p_{H'}(T(h),\ome') = \p_{H}(h,\hat{T}(\ome'))$ for all $h \in H$ and $\ome' \in \du{H'}$.
\end{prop}
\pr
It is obvious that the linear map $\hat{T}:\ome' \in \du{H'} \mapsto \ome' T \in \du{H}$ is the unique Hopf algebra isomorphism which satisfies the proposition.
\fin

\subsection{The canonical element associated to a pairing}

Let $H$ be a Hopf algebra and $\p$ be its canonical pairing. Consider the bilinear form $\p^{2} : H \odo \du{H} \times \du{H} \odo H \to \ku$ defined by $\p^{2}(h \odo \ome,\ome' \odo h')= \p(h,\ome')\bar{\p}(\ome,h')$. There exist an unique element $U \in H \odo \du{H}$ such that $\p^{2}(U, \ome \odo h) = \p(h,\ome)$ for all $h \in H$ and $\ome \in \du{H}$. Moreover, $U$ is invertible in $H \odo \du{H}$ with $U^{-1} = (S \odo \id_{\du{H}})(U) = (\id_{H} \odo \du{S})(U)$. This unique element $U$ will be called {\em the canonical element for the pairing $\p$}, and it will be denoted by $U(\p)$ or simply by $U$. This canonical element satisfies the following relations:
\[
(\com \odo \id)(U) = U_{13}U_{23}, \quad (\id \odo \dcom)(U) = U_{12}U_{13},
\]
\[
(\p(\,\cdot\,,\ome) \odo \id)(U) = \ome, \quad (\id \odo \p(h,\,\cdot\,))(U) = h.
\]

\begin{rema}
Consider a basis $\{h_{i}\}$ in $H$ with dual basis $\{\du{h}_{i}\}$ in $\du{H}$, i.e. $\p(h_{i},\du{h}_{j}) = \delta_{ij}$ for $1 \leq i,j \leq \text{dim}(H)$. It can be proof that $\sum_{i}h_{i} \odo \du{h}_{i}$ do not depend on the choice of the basis $\{h_{i}\}$. In this case, the canonical element $U$ for the canonical pairing $\p$ is given by $\sum_{i}h_{i} \odo \du{h}_{i}$.
\end{rema}

\begin{rema}
It is obvious that the canonical element for the pairing $\bar{\p} = \p\circ\chi_{\du{H},H}$ is given by $\Sigma_{H,\du{H}}(U) \in \du{H} \odo H$. Here, we are using the bilinear flip map $\chi_{\du{H},H}: \du{H} \times H \to H \times \du{H}$, $(\ome,h) \mapsto (h,\ome)$ and the linear flip map $\Sigma_{H,\du{H}}: H \odo \du{H} \to \du{H} \odo H$, $h \odo \ome \mapsto \ome \odo h$.
\end{rema}

\begin{rema}\label{rem:Ucirc_U}
Consider the bijective anti-homomorphism of Hopf algebras ${}^{\op}: H \to H^{\op}$, $h \mapsto h^{\op}$ and define ${}^{\ops}U:= ({}^{\op} \odo \id_{\du{H}})(U) \in H^{\op} \odo \du{H}$. We have that ${}^{\ops}U$ is the canonical element for the canonical pairing of $H^{\op}$. Moreover, because $({}^{\ops}U)^{-1} = ({}^{\op} \odo \dS^{-1})(U)$ and ${}^{\ops}U = ({}^{\op} \odo \dS^{-1})(U^{-1})$, then we have
\[
\ad({}^{\ops}U) = ({}^{\op} \odo \id_{\du{H}}) T_{U} ({}^{\op} \odo \id_{\du{H}}), \quad \ad({}^{\ops}U)^{-1} = ({}^{\op} \odo \id_{\du{H}}) T_{U^{-1}} ({}^{\op} \odo \id_{\du{H}})
\]
where $T_{U} := (\id_{H} \odo \du{S}^{-1})\ad(U)(\id_{H} \odo \du{S})$ and $T_{U^{-1}} := (\id_{H} \odo \du{S}^{-1})\ad(U^{-1})(\id_{H} \odo \du{S})$, are invertible elements in $\mathrm{End}(H \odo \du{H})$, such that $T^{-1}_{U} = T_{U^{-1}}$.
\end{rema}

\subsection{Actions and coactions of Hopf algebras} We give the basic notions about actions and coactions in the framework of Hopf algebras that will be use very often along this work.

\begin{defi}
An algebra $X$ is called {\em a left $H$-module algebra} or simply {\em a left $H$-algebra}, if $X$ is a left $H$-module with module map $\rhd: H \odo X \to X$ satisfying $h \rhd (xy) = (h_{(1)} \rhd x)(h_{(2)} \rhd y)$ and $h \rhd 1_{X} = \cou_{H}(h)1_{X}$ for all $h \in H$, $x,y \in X$. The module map $\rhd$ is called {\em a left action of the Hopf algebra $H$ on $X$}. Similarly, we can define a right $H$-module algebra.
\end{defi}

For a left $H$-algebra $X$ with left action $\rhd$, we can associate a new algebra called {\em the smash product algebra} denotes by $X \,\#_{\rhd} H$ and sometimes simply denotes by $X\,\#\,H$ if it is understand the action map. This algebra is the $\ku$-vector space $X \odo H$ endowed with the algebra structure  given by $(x \odo h)(y \odo h') = x(h_{(1)} \rhd y) \odo h_{(2)}h'$ for all $x,y \in X$, $h,h' \in H$. We will write $x \,\# h$ instead of $x \odo h$ for the basic elements of this algebra. The unit element of the smash product algebra is given by $1_{X} \,\# 1_{H}$.

\begin{defi}
Let $H$ be a Hopf algebra. An algebra $X$ is called {\em a left $H$-comodule algebra} if there is an homomorphism $\alp: X \to H \odo X$ satisfying $(\id_{H} \odo \alp)\alp = (\com_{H} \odo \id_{X})\alp$ and $(\cou_{H} \odo \id_{X})\alp = \id_{X}$. The map $\alp$ is called {\em a left coaction of the Hopf algebra $H$ on $X$}. Similarly, we can define a right $H$-comodule algebra and right coactions.
\end{defi}

\begin{rema}
The condition $(\cou_{H} \odo \id_{X})\alp = \id_{X}$ is equivalent to the injectivity of the map $\alp$. Indeed, if we have  $(\cou_{H} \odo \id_{X})\alp = \id_{X}$, it is obvious that $\alp$ must be injective. Reciprocally, if $\alp$ is injective, because $\alp((\cou \odo \id)(\alp(x))) = (\cou \odo \id \odo \id)(\id \odo \alp)(\alp(x)) = ((\cou \odo \id)\com \odo \id)(\alp(x)) = \alp(x)$ for all $x \in X$, we have $(\cou \odo \id)\alp = \id_{X}$.
\end{rema}

\begin{nota}
We will use the following {\em Sweedler type leg notations} for coactions: For a left coaction $\alp : X \to H \odo X$, we will use the subscript notation $x_{[-1]} \odo x_{[0]} := \alp(x)$, and for a right coaction $\alp : X \to X \odo H$, we will use the subscript $x_{[0]} \odo x_{[1]} := \alp(x)$. If we deal with coactions of the dual Hopf algebra $\du{H}$, sometimes we will use superscript notations. 
Note that we omit the sum symbol along the tensor product elements, i.e. instead to write $\alp(x) = \sum_{i}x^{i}_{[-1]} \odo x^{i}_{[0]}$, we write simply $\alp(x) = x_{[-1]} \odo x_{[0]}$.
\end{nota}

\begin{nota}
For each invertible element $W \in H \odo X$ such that $(\com_{H} \odo \id_{X})(W)=W_{23}W_{13}$, we have the left coaction $\alp_{W}: X \to H \odo X$ defined by $\alp_{W}(x) = \ad(W)(1_{H} \odo x) = W(1_{H} \odo x)W^{-1}$ for all $x \in X$. Similarly, for each invertible element $W \in X \odo H$ such that $(\id_{X} \odo \com_{H})(W)=W_{12}W_{13}$, we have the right coaction $\alp_{W}: X \to X \odo H$ define by $\alp_{W}(x) = \ad(W)(x \odo 1_{H}) = W(x \odo 1_{H})W^{-1}$ for all $x \in X$.
\end{nota}

\begin{rema}\label{rema:op_coactions}
If $\alp: X \to H \odo X$ is a left coaction of $H$ on $X$, its {\em opposite coaction} is given by the map $\alp^{\op}: X^{\op} \to H^{\op} \odo X^{\op}$, $\alp^{\op}(x^{\op}) = ({}^{\op} \odo {}^{\op})(\alp(x))$ for any $x \in X$. This is a left coaction of the opposite Hopf algebra $H^{\op}$ on $X^{\op}$. If $\alp: X \to H^{\op} \odo X$ is a left coaction of $H^{\op}$ on $X$, its {\em opposite coaction} is given by the map $\alp^{\op}: X^{\op} \to H \odo X^{\op}$, $\alp^{\op}(x^{\op}) = ({}^{\op} \odo {}^{\op})(\alp(x))$ for any $x \in X$. This is a left coaction of the Hopf algebra $H$ on $X^{\op}$. Similar definitions can be done for right coactions.
\end{rema}

\begin{rema}\label{rema:co_coactions}
If $\alp: X \to H \odo X$ is a left coaction of $H$ on $X$, its {\em co-opposite coaction} is given by the map $\alp^{\co}: X^{\op} \to H^{\co} \odo X^{\op}$, $\alp^{\co}(x^{\op}) = (S \odo {}^{\op})(\alp(x))$ for any $x \in X$. This is a left coaction of the co-opposite Hopf algebra $H^{\co}$ on $X^{\op}$. If $\alp: X \to H^{\op} \odo X$ is a left coaction of $H^{\op}$ on $X$, its {\em co-opposite coaction} is given by the map $\alp^{\co}: X^{\op} \to H^{\op,\co} \odo X^{\op}$, $\alp^{\co}(x^{\op}) = (S^{\ops} \odo {}^{\op})(\alp(x))$ for any $x \in X$. This is a left coaction of the opposite co-opposite Hopf algebra $H^{\op,\co}$ on $X^{\op}$. Similar definitions can be done for right coactions.
\end{rema}

The canonical element associated to the canonical pairing $\p$ of $H$ allows to give the following equivalence of categories.

\begin{prop}
Let ${}_{H}\alg$ be the category of left $H$-module algebras, i.e. the category with objects given by $(X,\rhd: H \odo X \to X)$ and morphisms given by $f: X \to X'$ such $f\circ\rhd = \rhd'\circ(\id \odo f)$, and $\alg^{H}$ be the category of right $H$-comodule algebras, i.e. the category with objects given by $(X,\alp: X \to X \odo H)$ and morphisms given by $f: X \to X'$ such that $(f \odo \id)\circ\alp = \alp'\circ f$. The functors
\[
\begin{array}{rccc}
& \alg^{H} & \leadsto & {}_{\du{H}}\alg \\
& (X,\alp) & \mapsto & (X,\rhd_{\alp}) \\
& f & \mapsto & f 
\end{array}
\quad\text{ and }\quad
\begin{array}{rccc}
& {}_{H}\alg & \leadsto & \alg^{\du{H}} \\
& (X,\rhd) & \mapsto & (X,\alp_{\rhd}) \\
& f & \mapsto & f 
\end{array}
\]
are equivalences of categories. The action $\rhd_{\alp}$ and the coaction $\alp_{\rhd}$ are given by
\[
\ome \rhd_{\alp} x := (\id_{X} \odo \p(\,\cdot\,,\ome))(\alp(x))
\quad \text{ and } \quad
\alp_{\rhd}(x) := (U_{1} \rhd x) \odo U_{2},
\]
respectively. Here $U_{1} \odo U_{2} = U \in H \odo \du{H}$ denotes the canonical element associated to the canonical pairing $\p$. Similar equivalences can be done for right actions and left coactions.
\end{prop}

\subsection{Canonical actions associated to a Hopf algebra}

We have the following canonical actions associated in a natural way with any Hopf algebra.

\subsubsection*{The regular actions} Let $H$ be a Hopf algebra. Using the canonical pairing $\p: H \times \du{H} \to \ku$, we can obtain the following actions called regular actions:
\[
\begin{array}{lccc}
\brhd : & H \odo \du{H} & \to & \du{H} \\
& h \odo \ome & \mapsto & \ome_{(1)}\p(h,\ome_{(2)})
\end{array}
\qquad
\begin{array}{lccc}
\blhd : & \du{H} \odo H & \to & \du{H} \\
& \ome \odo h & \mapsto & \p(h,\ome_{(1)})\ome_{(2)}
\end{array},  
\]
\[
\begin{array}{lccc}
\brhd : & \du{H} \odo H & \to & H \\
& \ome \odo h & \mapsto & h_{(1)}\p(h_{(2)},\ome)
\end{array}
\qquad
\begin{array}{lccc}
\blhd : & H \odo \du{H} & \to & H \\
& h \odo \ome & \mapsto & \p(h_{(1)},\ome)h_{(2)}
\end{array}.
\]

\subsubsection*{The adjoint actions} The left and right adjoint action of $(H,\com)$ on $H$ are respectively:
\[
\begin{array}{lccc}
\brhd_{ad}: & H \odo H & \to & H \\
& h \odo g & \mapsto & h_{(1)}gS(h_{(2)})
\end{array}
\qquad
\begin{array}{lccc}
\blhd_{ad}: & H \odo H & \to & H \\
& g \odo h & \mapsto & S(h_{(1)})gh_{(2)}
\end{array}. 
\]

\subsubsection*{The coadjoint actions} To obtain this actions, we dualize the adjoint actions. For example, we take the left adjoint action of $(H,\com)$, and we construct its canonical dual right action of $(H,\com)$ on $\du{H}$, i.e. $(\ome \lhd h): h' \mapsto \ome(h \brhd_{ad} h')$ for all $h \in H$ and $\ome \in \du{H}$, then we change to a left action using the inverse of the antipode. The left and right coadjoint action of $(H,\com)$ on $\du{H}$ are respectively:

\[
\begin{array}{lccc}
\brhd_{coad}: & H \odo \du{H} & \to & \du{H} \\
& h \odo \ome & \mapsto & h_{(1)} \brhd \ome \blhd S^{-1}(h_{(2)})
\end{array}
\qquad
\begin{array}{lccc}
\blhd_{coad}: & \du{H} \odo H & \to & \du{H} \\
& \ome \odo h & \mapsto & S^{-1}(h_{(1)}) \brhd \ome \blhd h_{(2)}
\end{array}.
\]

\subsection{Calculations using the canonical element of a pairing}

In this section, we give some calculations using this canonical element of $\p$. Associated to each $h \in H$ and each $\ome \in \du{H}$, we have the canonical linear maps ${}_{h}\p: \du{H} \to \ku$, $\theta \mapsto \p(h,\theta)$ and $\p_{\ome}: H \to \ku$, $h' \mapsto \p(h',\ome)$.

\begin{prop}
For any $h,g \in H$, $\ome \in \du{H}$, it hold
\begin{enumerate}[label=\textup{(\roman*)}]
\item $(\id \odo {}_{h}\p)(U(g \odo \ome)U^{-1}) = (\ome \brhd h_{(1)})gS(h_{(2)}) = h_{(1)}gS(h_{(2)} \blhd \ome)$,
\item $(\id \odo {}_{h}\p)(U^{-1}(g \odo \ome)U) = S(h_{(1)})g(h_{(2)} \blhd \ome) =  S(\ome \brhd h_{(1)})gh_{(2)}$,
\item $(\id \odo {}_{h}\p)(T_{U}(g \odo \ome)) = S^{-1}(h_{(2)} \blhd \ome)gh_{(1)} = S^{-1}(h_{(2)})g(\ome \brhd h_{(1)})$,
\item $(\id \odo {}_{h}\p)(T_{U^{-1}}(g \odo \ome)) = h_{(2)}gS^{-1}(\ome \brhd h_{(1)}) = (h_{(2)} \blhd \ome)gS^{-1}(h_{(1)})$.
\end{enumerate}
Where
\[
T_{U} = (\id_{H} \odo \du{S}^{-1})\ad(U)(\id_{H} \odo \du{S}) \quad \text{ and } \quad T_{U^{-1}} = (\id_{H} \odo \du{S}^{-1})\ad(U^{-1})(\id_{H} \odo \du{S})
\]
are elements in $\mathrm{End}(H \odo \du{H})$.
\end{prop}
\pr
\begin{enumerate}[label=\textup{(\roman*)}]
\item Because
\begin{align*}
\p((\id \odo {}_{h}\p)(U(g \odo \ome)U^{-1}),\theta) & = \p^{2}(U(g \odo \ome)U^{-1},\theta \odo h) \\
& = \p^{2}(U,\theta_{(1)} \odo h_{(1)})\p^{2}(g \odo \ome,\theta_{(2)} \odo h_{(2)})\p^{2}(U^{-1},\theta_{(3)} \odo h_{(3)}) \\
& = \p^{2}(U,\theta_{(1)} \odo h_{(1)})\p(g,\theta_{(2)})\p(h_{(2)},\ome)\p^{2}((\id \odo \du{S})(U),\theta_{(3)} \odo h_{(3)}) \\
& = \p^{2}(U,\theta_{(1)} \odo h_{(1)})\p(g,\theta_{(2)})\p(h_{(2)},\ome)\p^{2}(U,\theta_{(3)} \odo S(h_{(3)})) \\
& = \p(h_{(1)},\theta_{(1)})\p(g,\theta_{(2)})p(h_{(2)},\ome)\p(S(h_{(3)}),\theta_{(3)}) \\
& = \p(h_{(1)}gS(h_{(3)}),\theta)\p(h_{(2)},\ome) \\
& = \p(h_{(1)}\p(h_{(2)},\ome)gS(h_{(3)}),\theta) \\
& = \p((\ome \brhd h_{(1)})gS(h_{(2)}),\theta)
\end{align*}
for all $\theta \in \du{H}$, we have the result.

\item Because
\begin{align*}
\p((\id \odo {}_{h}\p)(U^{-1}(g \odo \ome)U),\theta) & = \p^{2}(U^{-1}(g \odo \ome)U,\theta \odo h) \\
& = \p^{2}(U^{-1},\theta_{(1)} \odo h_{(1)})\p^{2}(g \odo \ome,\theta_{(2)} \odo h_{(2)})\p^{2}( U,\theta_{(3)} \odo h_{(3)}) \\
& = \p^{2}((\id \odo \du{S})(U),\theta_{(1)} \odo h_{(1)})\p(g,\theta_{(2)})\p(h_{(2)},\ome)\p^{2}(U,\theta_{(3)} \odo h_{(3)}) \\
& = \p^{2}(U,\theta_{(1)} \odo S(h_{(1)}))\p(g,\theta_{(2)})\p(h_{(2)},\ome)\p^{2}(U,\theta_{(3)} \odo h_{(3)}) \\
& = \p(S(h_{(1)}),\theta_{(1)})\p(g,\theta_{(2)})\p(h_{(2)},\ome)\p(h_{(3)},\theta_{(3)}) \\
& = \p(S(h_{(1)})gh_{(3)},\theta)p(h_{(2)},\ome) \\
& = \p(S(h_{(1)})g\p(h_{(2)},\ome)h_{(3)},\theta) \\
& = \p(S(h_{(1)})g(h_{(2)} \blhd \ome),\theta)
\end{align*}
for all $\theta \in \du{H}$, we have the result.

\item Because
\begin{align*}
\p((\id \odo {}_{h}\p)(T_{U}(g \odo \ome)),\theta) & = \p^{2}((\id \odo \du{S}^{-1})\ad(U)(\id \odo \du{S})(g \odo \ome),\theta \odo h) \\
& = \p^{2}(U(g \odo \du{S}(\ome))U^{-1},\theta \odo S^{-1}(h)) \\
& = \p((\id \odo {}_{S^{-1}(h)}\p)(U(g \odo \du{S}(\ome))U^{-1}),\theta) \\
& = \p((\du{S}(\ome) \brhd S^{-1}(h)_{(1)})gS(S^{-1}(h)_{(2)}),\theta) \\
& = \p((\du{S}(\ome) \brhd S^{-1}(h_{(2)}))gh_{(1)},\theta) \\
& = \p(S^{-1}(h_{(2)} \blhd \ome)gh_{(1)},\theta)
\end{align*}
for all $\theta \in \du{H}$, we have the result.

\item Because
\begin{align*}
\p((\id \odo {}_{h}\p)(T_{U^{-1}}(g \odo \ome)),\theta) & = \p^{2}((\id \odo \du{S}^{-1})\ad(U^{-1})(\id \odo \du{S})(g \odo \ome),\theta \odo h) \\
& = \p^{2}(U^{-1}(g \odo \du{S}(\ome))U,\theta \odo S^{-1}(h)) \\
& = \p((\id \odo {}_{S^{-1}(h)}\p)(U^{-1}(g \odo \du{S}(\ome))U),\theta) \\
& = \p(S(S^{-1}(h)_{(1)})g(S^{-1}(h)_{(2)} \blhd \du{S}(\ome)),\theta) \\
& = \p(h_{(2)}g S^{-1}(\ome \brhd h_{(1)}),\theta)
\end{align*}
for all $\theta \in \du{H}$, we have the result.
\end{enumerate}
\fin

\subsection{Heisenberg and Drinfeld double of a Hopf algebra}

In this subsection, we recall two important constructions related to any finite-dimensional Hopf algebra.

\subsubsection*{Heisenberg algebra}\label{sec:heisenberg-double}

We recall the construction of the Heisenberg algebra for a Hopf algebra.

\begin{defi}[\cite{STS92}]
The {\em Heisenberg algebra of a Hopf algebra $H$}, denoted by $\He(H)$, is by definition the smash product algebra $H \# \du{H}$ with respect to the left regular action of $\du{H}$ on $H$, i.e. for every $h,h' \in H$, $\theta, \theta' \in \du{H}$ we have $(h \# \theta)(h' \# \theta') := h(\theta_{(1)} \brhd h') \# \theta_{(2)}\theta'$.
\end{defi}
\begin{rema}
Because
\[
(h \# \theta)(h' \# \theta') = hh'_{(1)}\p(h'_{(2)},\theta_{(1)}) \# \theta_{(2)}\theta' = hh'_{(1)} \# \p(h'_{(2)},\theta_{(1)})\theta_{(2)}\theta' = hh'_{(1)} \# (\theta \blhd h'_{(2)})\theta
\]
The Heisenberg algebra can also be regarded as the smash product algebra $H \# \du{H}$ with respect to the right regular action of $H$ on $\du{H}$.
\end{rema}

\begin{prop}[\cite{L94}]\label{prop:anti_heisen}
The maps
\[
\begin{array}{lccc}
\mathcal{L}_{1} : & \He(\du{H}) & \to & \He(H) \\
& \theta \# a & \mapsto & S^{-1}(a) \# \du{S}(\theta)
\end{array}
\qquad
\begin{array}{lccc}
\mathcal{L}_{2} : & \He(\du{H}) & \to & \He(H) \\
& \theta \# a & \mapsto & S(a) \# \du{S}^{-1}(\theta)
\end{array}
\]
are algebra anti-isomorphisms.
\end{prop}

Using the canonical inclusions maps $\iota_{1}: H \to \He(H)$, $\iota_{2}: \du{H} \to \He(H)$, $\du{\iota}_{1}: \du{H} \to \He(\du{H})$, $\du{\iota}_{2}: H \to \He(\du{H})$, the above proposition gives us the following relations
\begin{align}\label{eq:heis-dual-inclu-1}
\mathcal{L}_{1}\du{\iota}_{1} = \iota_{2}\du{S}, \;\; \mathcal{L}_{1}\du{\iota}_{2} = \iota_{1} S^{-1}, \;\; \mathcal{L}_{2}\du{\iota}_{1} = \iota_{2}\du{S}^{-1}, \;\; \mathcal{L}_{2}\du{\iota}_{2} = \iota_{1} S.
\end{align}

\subsubsection*{Drinfeld double}\label{sec:drinfeld-double}

Let $U = U_{1} \odo U_{2}$ be the canonical element for the Hopf algebra $H$. Using the invertible element ${}^{\ops}U = U^{\op}_{1} \odo U_{2} = ({}^{\op} \odo \id)(U) \in H^{\op} \odo \du{H}$, we consider the twisted Hopf algebra of $H^{\op}$ and $\du{H}$, denoted by $\T_{R}(H) := H^{\op} \odo \du{H}$, defined by the algebra structure on $H^{\op} \odo \du{H}$ given by $\m_{\T_{R}(H)} = (\m^{\op}_{H} \odo \m_{\du{H}})(\id_{H^{\op}} \odo \Sigma \odo \id_{\du{H}})$ and the coalgebra structure on $H^{\op} \odo \du{H}$ given by
\[
\com_{\T_{R}(H)} = (\id_{H^{\op}} \odo \Sigma\ad({}^{\ops}U) \odo \id_{\du{H}})(\com_{H^{\op}} \odo \com_{\du{H}}).
\]
This Hopf algebra $\T_{R}(H)$ will be called the {\em Drinfeld codouble} of $H$.

\begin{prop}[\cite{L94}]
The {\em Drinfeld double} of a Hopf algebra $H$ is the dual Hopf algebra of $\T_{R}(H)$ and it will be denoted by $\D_{R}(H)$.
\end{prop}

\begin{rema}
Using the canonical pairing $\p_{\D} := \tilde{\p} \odo \p: \du{H}^{\co} \odo H \times H^{\op} \odo \du{H} \to \ku$ associate to $\D_{R}(H)$ and $\T_{R}(H)$, we can give the explicit structure of $\D_{R}(H)$.
\begin{itemize}
\item As vector space $\D_{R}(H)$ is the tensor product vector space $\du{H} \odo H$.
\item 
It follows from
\small
\begin{align*}
\p_{\D}(m_{\D_{R}(H)}(\ome \bicroi a, \hat{\ome} \bicroi b),c^{\op} \odo \theta) & = (\p_{\D} \odo \p_{\D})(\ome \bicroi a \odo \hat{\ome} \bicroi b,\com_{\T_{R}(H)}(c^{\op} \odo \theta)) \\
& = \bar{\p}(\ome,c_{(1)})\bar{\p}^{2}(\hat{\ome} \odo a,({}^{\ops}U(c^{\op}_{(2)} \odo \theta_{(1)})({}^{\ops}U)^{-1})\p(b,\theta_{(2)}) \\
& = \bar{\p}(\ome,c_{(1)})\bar{\p}(\hat{\ome},S^{-1}(a_{(2)})c_{(2)}(\theta_{(1)} \brhd a_{(1)}))\p(b,\theta_{(2)}) \\
& = \bar{\p}(\ome,c_{(1)})\bar{\p}(\hat{\ome}_{(1)},S^{-1}(a_{(2)})c_{(2)})\bar{\p}(\hat{\ome}_{(2)},\theta_{(1)} \brhd a_{(1)})\p(b,\theta_{(2)}) \\
& = \bar{\p}(\ome,c_{(1)})\bar{\p}(\hat{\ome}_{(1)},S^{-1}(a_{(2)})c_{(2)})\bar{\p}(\hat{\ome}_{(2)}\theta_{(1)},a_{(1)})\p(b,\theta_{(2)}) \\
& = \bar{\p}(\ome,c_{(1)})\bar{\p}(\hat{\ome}_{(1)},S^{-1}(a_{(3)})c_{(2)})\bar{\p}(\hat{\ome}_{(2)},a_{(1)})\bar{\p}(\theta_{(1)},a_{(2)})\p(b,\theta_{(2)}) \\
& = \bar{\p}(\ome,c_{(1)})\bar{\p}(\hat{\ome},S^{-1}(a_{(3)})c_{(2)}a_{(1)})\p(a_{(2)}b,\theta) \\
& = \bar{\p}(\ome,c_{(1)})\bar{\p}(a_{(1)} \brhd \hat{\ome} \blhd S^{-1}(a_{(3)}),c_{(2)})\p(a_{(2)}b,\theta) \\
& = \p_{\D}(\ome(a_{(1)} \brhd \hat{\ome} \blhd S^{-1}(a_{(3)})) \bicroi a_{(2)}b,c^{\op} \odo \theta),
\end{align*}
\normalsize
that, the algebra structure of $\D_{R}(H)$ is given by
\[
\m_{\D_{R}(H)} = (\m_{\du{H}^{\co}} \odo \m_{H})(\id_{\du{H}^{\co}} \odo \mathfrak{T} \odo \id_{H}),
\]
where $\mathfrak{T} : H \odo \du{H}^{\co} \to \du{H}^{\co} \odo H$ is the linear map given by
\begin{align*}
\mathfrak{T}(a \odo \hat{\ome}) & = (a_{(1)} \brhd \hat{\ome} \blhd S^{-1}(a_{(3)})) \odo a_{(2)}\\
& = \hat{\ome}_{(2)} \odo (\du{S}^{-1}(\hat{\ome}_{(1)}) \brhd a \blhd \hat{\ome}_{(3)}) \\
& = (a_{(1)} \brhd \hat{\ome}_{(2)}) \odo (\du{S}^{-1}(\hat{\ome}_{(1)}) \brhd a_{(2)})\\
& = (\hat{\ome}_{(1)} \blhd S^{-1}(a_{(2)})) \odo (a_{(1)} \blhd \hat{\ome}_{(2)}).
\end{align*}
for all $a \in H$ and $\hat{\ome} \in \du{H}$. The linear map $\mathfrak{T}$ can be given also by
\[
\mathfrak{T}(a \odo \hat{\ome}) = a_{(1)} \brhd_{coad} \hat{\ome}_{(2)} \odo a_{(2)} \blhd_{coad} \hat{\ome}_{(1)}
\]
for all $a \in H$ and $\hat{\ome} \in \du{H}$.

\item  It follows from
\footnotesize
\begin{align*}
(\p_{\D}\odo \p_{\D})(\com_{\D_{R}(H)}(\ome \bicroi a),b^{\op} \odo \theta \odo c^{\op} \odo \overline{\theta}) & = \p_{\D}(\ome \bicroi a, m_{\T_{R}(H)}(b^{\op} \odo \theta,c^{\op} \odo \overline{\theta})) \\
& = \p_{\D}(\ome \bicroi a, (cb)^{\op} \odo \theta\overline{\theta}) \\
& = \bar{\p}(\ome_{(1)},c)\bar{\p}(\ome_{(2)},b)\p(a_{(1)},\theta)\p(a_{(2)},\overline{\theta}) \\
& = (\p_{\D} \odo \p_{\D})(\ome_{(2)} \bicroi a_{(1)} \odo \ome_{(1)} \bicroi a_{(2)},b^{\op} \odo \theta \odo c^{\op} \odo \overline{\theta}) \\
& = (\p_{\D} \odo \p_{\D})((\id \odo \Sigma \odo \id)(\dcom^{\co} \odo \com)(\ome \bicroi a),b^{\op} \odo \theta \odo c^{\op} \odo \overline{\theta}),
\end{align*}
\normalsize
that the coalgebra structure of $\D_{R}(H)$ is given by
\[
\com_{\D_{R}(H)} = (\id_{\du{H}^{\co}} \odo \Sigma \odo \id_{H})(\dcom^{\co} \odo \com).
\]
\end{itemize}
In other words, the Radford's characterization of the Drinfeld double, is given by the double crossproduct construction $\D_{R}(H) = \du{H}^{\co} \bicroi H$.
\end{rema}

\begin{rema}[Majid's characterization of the Drinfeld double]
Using the invertible element $U^{-1} \in H^{\co} \odo \du{H}$, consider the twisted Hopf algebra of $H^{\co}$ and $\du{H}$, denoted by $\T_{M}(H) := H^{\co} \odo \du{H}$, defined by the algebra structure given by
\[
m_{\T_{M}(H)} = (m_{H^{\co}} \odo m_{\du{H}})(\id_{H^{\co}} \odo \Sigma \odo \id_{\du{H}})
\]
and the coalgebra structure given by
\[
\com_{\T_{M}(H)} = (\id_{H^{\co}} \odo \Sigma\ad(U^{-1}) \odo \id_{\du{H}})(\com_{H^{\co}} \odo \com_{\du{H}}).
\]
Let $\D_{M}(H)$ be the dual Hopf algebra of $\T_{M}(H)$ and $\mathds{P} := \tilde{\bar{\p}} \odo \p : \du{H}^{\op} \odo  H \times H^{\co} \odo \du{H} \to \ku$ be the canonical pairing associate to $\D_{M}(H)$ and $\T_{M}(H)$. The explicit structure of $\D_{M}(H)$ is
\begin{itemize}
\item As vector space $\D_{M}(H)$ is the tensor product vector space $\du{H}^{\op} \odo H$.
\item 
It follows from
\begin{align*}
\mathds{P}(m_{\D_{M}(H)}(\ome^{\op} \bicroi h, \hat{\ome}^{\op} \bicroi \hat{h}),g \odo \theta) & = (\mathds{P} \odo \mathds{P})(\ome^{\op} \bicroi h \odo \hat{\ome}^{\op} \bicroi \hat{h},\com_{\T_{M}(H)}(g \odo \theta)) \\
& = \bar{\p}(\ome,g_{(2)})\bar{\p}^{2}(\hat{\ome}^{\op} \odo h,U^{-1}(g_{(1)} \odo \theta_{(1)})U)\p(\hat{h},\theta_{(2)}) \\
& = \bar{\p}(\ome,g_{(2)})\bar{\p}(\hat{\ome},S(\theta_{(1)} \brhd h_{(1)})g_{(1)}h_{(2)})\p(\hat{h},\theta_{(2)}) \\
& = \bar{\p}(\ome,g_{(2)})\bar{\p}(h_{(2)} \brhd \hat{\ome}_{(2)},g_{(1)})\p(h_{(1)} \blhd \dS(\hat{\ome}_{(1)}),\theta_{(1)})\p(\hat{h},\theta_{(2)}) \\
& = \bar{\p}((h_{(2)} \brhd \hat{\ome}_{(2)})\ome,g)\p((h_{(1)} \blhd \dS(\hat{\ome}_{(1)}))\hat{h},\theta) \\
& = \mathds{P}(\ome^{\op}(h_{(2)} \brhd \hat{\ome}_{(2)})^{\op} \bicroi (h_{(1)} \blhd \dS(\hat{\ome}_{(1)}))\hat{h},g \odo \theta),
\end{align*}
that, the algebra structure of $\D_{M}(H)$ is given by
\[
\m_{\D_{M}(H)} = (\m_{\du{H}^{\op}} \odo \m_{H})(\id_{\du{H}^{\op}} \odo \mathfrak{T} \odo \id_{H}),
\]
where $\mathfrak{T} : H \odo \du{H}^{\op}\to \du{H}^{\op} \odo H$ is the linear map defined by
\begin{align*}
\mathfrak{T}(h\odo \hat{\ome}^{\op}) & = (h_{(3)} \brhd \hat{\ome} \blhd S(h_{(1)}))^{\op} \odo h_{(2)}\\
& = \hat{\ome}^{\op}_{(2)} \odo (\hat{\ome}_{(3)} \brhd h \blhd \dS(\hat{\ome}_{(1)})) \\
& = (h_{(2)} \brhd \hat{\ome}_{(2)})^{\op} \odo (h_{(1)} \blhd \dS(\hat{\ome}_{(1)})) \\
& = (\hat{\ome}_{(1)} \blhd S(h_{(1)}))^{\op} \odo (\hat{\ome}_{(2)} \brhd h_{(2)}).
\end{align*}
for all $h \in H$ and $\hat{\ome} \in \du{H}$.
\item  It follows from
\small
\begin{align*}
(\mathds{P} \odo \mathds{P})(\com_{\D_{M}(H)}(\ome^{\op} \bicroi h),g \odo \theta \odo \hat{g} \odo \hat{\theta}) & = \mathds{P}(\ome^{\op} \bicroi h, m_{\T_{M}(H)}(g \odo \theta,\hat{g} \odo \hat{\theta})) \\
& = \mathds{P}(\ome^{\op} \bicroi h, g\hat{g} \odo \theta\hat{\theta}) \\
& = \bar{\p}(\ome_{(1)},g)\bar{\p}(\ome_{(2)},\hat{g})\p(h_{(1)},\theta)\p(h_{(2)},\hat{\theta}) \\
& = (\mathds{P} \odo \mathds{P})(\ome^{\op}_{(1)} \bicroi h_{(1)} \odo \ome^{\op}_{(2)} \bicroi h_{(2)},g \odo \theta \odo \hat{g} \odo \hat{\theta}) \\
& = (\mathds{P} \odo \mathds{P})((\id \odo \Sigma \odo \id)(\dcom \odo \com)(\ome^{\op} \bicroi h),g \odo \theta \odo \hat{g} \odo \hat{\theta}),
\end{align*}
\normalsize
that the coalgebra structure of $\D_{M}(H)$ is given by
\[
\com_{\D_{M}(H)} = (\id_{\du{H}^{\op}} \odo \Sigma \odo \id_{H})(\com_{\du{H}^{\op}} \odo \com_{H}).
\]
\end{itemize}
In other words, the Majid's characterization of the Drinfeld double, is given by the double crossproduct construction $\D_{M}(H) = \du{H}^{\op} \bicroi H$.
\end{rema}

\begin{rema}[Taipe's characterization of the Drinfeld double \cite{TPHD18}]
Using the invertible element $U \in H \odo \du{H}^{\co}$, consider the twisted Hopf algebra of $H$ and $\du{H}^{\co}$, denoted by $\T_{T}(H) := H \odo \du{H}^{\co}$, defined by the algebra structure given by
\[
m_{\T_{T}(H)} = (m_{H} \odo m_{\du{H}^{\co}})(\id_{H} \odo \Sigma \odo \id_{\du{H}^{\co}})
\]
and the coalgebra structure given by
\[
\com_{\T_{T}(H)} = (\id_{H} \odo \Sigma\ad(U) \odo \id_{\du{H}^{\co}})(\com_{H} \odo \com_{\du{H}^{\co}}).
\]
Let $\D_{T}(H)$ be the dual Hopf algebra of $\T_{T}(H)$ and $\mathds{P} := \bar{\p} \odo \tilde{\p} : \du{H} \odo  H^{\op} \times H \odo \du{H}^{\co} \to \ku$ be the canonical pairing associate to $\D_{T}(H)$ and $\T_{T}(H)$. The explicit structure of $\D_{T}(H)$ is
\begin{itemize}
\item As vector space $\D_{T}(H)$ is the tensor product vector space $\du{H} \odo H^{\op}$.
\item 
It follows from
\begin{align*}
\mathds{P}(m_{\D_{T}(H)}(\ome \bicroi a^{\op}, \hat{\ome} \bicroi b^{\op}),c \odo \theta) & = (\mathds{P} \odo \mathds{P})(\ome \bicroi a^{\op} \odo \hat{\ome} \bicroi b^{\op},\com_{\T_{T}(H)}(c \odo \theta)) \\
& = \bar{\p}(\ome,c_{(1)})\bar{\p}^{2}(\hat{\ome} \odo a,U(c_{(2)} \odo \theta_{(2)})U^{-1})\p(b,\theta_{(1)}) \\
& = \bar{\p}(\ome,c_{(1)})\bar{\p}(\hat{\ome},(\theta_{(2)} \brhd a_{(1)})c_{(2)}S(a_{(2)}))\p(b,\theta_{(1)}) \\
& = \bar{\p}(\ome,c_{(1)})\bar{\p}(\hat{\ome}_{(1)}\theta_{(2)},a_{(1)})\bar{\p}(\hat{\ome}_{(2)},c_{(2)}S(a_{(2)}))\p(b,\theta_{(1)}) \\
& = \bar{\p}(\ome(S(a_{(2)}) \brhd \hat{\ome}_{(2)}),c)\p(b(a_{(1)} \blhd \hat{\ome}_{(1)}),\theta) \\ 
& = \mathds{P}(\ome(S(a_{(2)}) \brhd \hat{\ome}_{(2)}) \bicroi (a_{(1)} \blhd \hat{\ome}_{(1)})^{\op}b^{\op},c \odo \theta),
\end{align*}
that, the algebra structure of $\D_{T}(H)$ is given by
\[
\m_{\D_{T}(H)} = (\m_{\du{H}} \odo \m_{H^{\op}})(\id_{\du{H}} \odo \mathfrak{T} \odo \id_{H^{\op}}),
\]
where $\mathfrak{T} : H^{\op} \odo \du{H} \to \du{H} \odo H^{\op}$ is the linear map defined by
\begin{align*}
\mathfrak{T}(a^{\op} \odo \hat{\ome}) & = (S(a_{(3)}) \brhd \hat{\ome} \blhd a_{(1)}) \odo a^{\op}_{(2)}\\
& = \hat{\ome}_{(2)} \odo (\dS(\hat{\ome}_{(3)}) \brhd a \blhd \hat{\ome}_{(1)})^{\op} \\
& = (S(a_{(2)}) \brhd \hat{\ome}_{(2)}) \odo (a_{(1)} \blhd \hat{\ome}_{(1)})^{\op} \\
& = (\hat{\ome}_{(1)} \blhd a_{(1)}) \odo (\dS(\hat{\ome}_{(2)}) \brhd a_{(2)})^{\op}.
\end{align*}
for all $a \in H$ and $\hat{\ome} \in \du{H}$.
\item  It follows from
\small
\begin{align*}
(\mathds{P} \odo \mathds{P})(\com_{\D_{T}(H)}(\ome \bicroi a^{\op}),b \odo \theta \odo c \odo \hat{\theta}) & = \mathds{P}(\ome \bicroi a^{\op}, m_{\T_{T}(H)}(b \odo \theta,c \odo \hat{\theta})) \\
& = \mathds{P}(\ome \bicroi a^{\op}, bc \odo \theta\hat{\theta}) \\
& = \bar{\p}(\ome_{(1)},b)\bar{\p}(\ome_{(2)},c)\p(a_{(1)},\theta)\p(a_{(2)},\hat{\theta}) \\
& = (\mathds{P} \odo \mathds{P})(\ome_{(1)} \bicroi a^{\op}_{(1)} \odo \ome_{(2)} \bicroi a^{\op}_{(2)},b \odo \theta \odo c \odo \hat{\theta}) \\
& = (\mathds{P} \odo \mathds{P})((\id \odo \Sigma \odo \id)(\dcom \odo \com)(\ome \bicroi a^{\op}),b \odo \theta \odo c \odo \hat{\theta}),
\end{align*}
\normalsize
that the coalgebra structure of $\D_{T}(H)$ is given by
\[
\com_{\D_{T}(H)} = (\id_{\du{H}^{\co}} \odo \Sigma \odo \id_{H})(\com_{\du{H}^{\co}} \odo \com_{H}).
\]
\end{itemize}
In other words, the Taipe's characterization of the Drinfeld double, $\D_{T}(H)$, is given by the double crossproduct construction $\D_{T}(H) = \du{H}^{\co} \bicroi H$.
\end{rema}

\subsection{Twisting of Hopf algebras by $2$-cocycles}\label{sec:twist-hopf-algebr}

\begin{defi}[\cite{BCM86}]
Let $H$ be a Hopf algebra. A bilinear map $\sigma: H \times H \to \ku$ is called a {\em left $2$-cocycle on $H$} if it satisfies the left cocycle condition
\[
\sigma(h_{(1)},h'_{(1)})\sigma(h_{(2)}h'_{(2)},h'') = \sigma(h_{(1)},h'''_{(1)})\sigma(h,h'_{(2)}h''_{(2)})
\]
for all $h,h',h'' \in H$. Similarly, $\sigma$ is said to be a {\em right $2$-cocycle} if it satisfies the right cocycle condition
\[
\sigma(h_{(1)}h'_{(1)},h'')\sigma(h_{(2)},h'_{(2)}) = \sigma(h,h'_{(1)}h''_{(1)})\sigma(h'_{(2)},h''_{(2)})
\]
for all $h,h',h'' \in H$. It is said to be {\em normal} if $\sigma(1,h) = \sigma(h,1) = \cou(h)$ for all $h \in H$.
\end{defi}

\begin{exam}
The trivial cocycle on $H$ is the one given by
\[
\sigma_{\text{tr}}(h,h') = \p\odo\p(h \odo h',\cou \odo \cou) = \p(h,\cou)\p(h',\cou).
\]
\end{exam}

Given a bilinear map $\sigma: H \times H \to \ku$, we can define two new products $\cdot_{\sigma}$ and ${}_{\sigma}\cdot$ on the $\ku$-vector space $H$ as follows
\[
h \cdot_{\sigma} h' := \sigma(h_{(1)},h'_{(1)})h_{(2)}h'_{(2)}, \qquad h \,{}_{\sigma}\cdot h' := h_{(1)}h'_{(1)}\sigma(h_{(2)},h'_{(2)}).
\]
for all $h, h' \in H$.
\begin{prop}[\cite{L94}]
We have
\begin{itemize}
\item The product ${}_{\sigma}\cdot$ is associative if and only if $\sigma$ is a left $2$-cocycle.
\item The product $\cdot_{\sigma}$ is associative if and only if $\sigma$ is a right $2$-cocycle.
\item $1_{H}$ is the unit element for ${}_{\sigma}\cdot$ or $\cdot_{\sigma}$ if and only if $\sigma$ is normal.
\end{itemize}
\end{prop}

\begin{nota}
We will use the notation $H_{\sigma}$ (resp.  ${}_{\sigma}H$) for the $\ku$-vector space $H$ with the product $\cdot_{\sigma}$ (resp. ${}_{\sigma}\cdot$).
\end{nota}

\subsection{Relations between the Heisenberg double and the Drinfeld double}

We recall here important theorems about the relationship between the Heisenberg double construction and the Radford's characterization of the Drinfeld double for Hopf algebras.

\begin{theo}[\cite{L94}]\label{theo:Double_Heisenberg}
The Heisenberg double $\He(H)$ of the Hopf algebra $H$ is the left twist of the Hopf algebra $\T_{R}(H) = \du{\D_{R}(H)}$ by the left $2$-cocycle on $\T_{R}(H)$ given by
\[
\begin{array}{rccc}
\sigma : & \T(H) \odo \T(H) & \to & \ku \\
& h \odo \theta \odo h' \odo \theta' & \mapsto & \p(h,1)\bar{\p}(\theta,h')\p(1,\theta')
\end{array}
\]
i.e.
\[
\He(H) \iso {}_{\sigma}\T_{R}(H).
\]

The Heisenberg double $\He(\du{H})$ of the Hopf algebra $\du{H}$ is the right twist of the Hopf algebra $\D_{R}(H)$ by the right $2$-cocycle on $\D_{R}(H)$ given by
\[
\begin{array}{rccc}
\eta : & \D_{R}(H) \odo \D_{R}(H) & \to & \ku \\
& \theta \odo h \odo \theta' \odo h' & \mapsto & \bar{\p}(\theta,1)\p(h,\theta')\bar{\p}(1,h')
\end{array}
\]
i.e.
\[
\He(\du{H}) \iso \D_{R}(H)_{\eta}.
\]
\end{theo}

Moreover, we also have

\begin{theo}[\cite{L94}]
The Heisenberg double $\He(H)$ of the Hopf algebra $H$ is the right twist of the Hopf algebra $\T_{R}(H) = \du{\D_{R}(H)}$ by the right $2$-cocycle on $\T_{R}(H)$ given by
\[
\begin{array}{rccc}
\sigma^{-1} : & \T_{R}(H) \odo \T_{R}(H) & \to & \ku \\
& h \odo \theta \odo h' \odo \theta' & \mapsto & \p(h,1)\bar{\p}(\dS(\theta),h')\p(1,\theta')
\end{array}
\]
i.e.
\[
\He(H) \iso \T_{R}(H)_{\sigma^{-1}}.
\]

The Heisenberg double $\He(\du{H})$ of the Hopf algebra $\du{H}$ is the left twist of the Hopf algebra $\D_{R}(H)$ by the left $2$-cocycle on $\D_{R}(H)$ given by
\[
\begin{array}{rccc}
\eta^{-1} : & \D_{R}(H) \odo \D_{R}(H) & \to & \ku \\
& \theta \odo h \odo \theta' \odo h' & \mapsto & \bar{\p}(\theta,1)\p(h,\dS(\theta'))\bar{\p}(1,h')
\end{array}
\]
i.e.
\[
\He(\du{H}) \iso {}_{\eta^{-1}}\D_{R}(H).
\]
\end{theo}

\section{The standard characterization of Yetter--Drinfeld algebras}\label{sec:standard}

Here, we recall the standard characterization of Yetter--Drinfeld algebras over finite-dimensional Hopf algebras. Let $(H,\com,\cou,S)$ be a Hopf algebra and $X$ be an algebra. 

\begin{rema}
In all this section, for each definition we use the first adjective for the side of the action and the second adjective for the side of the coaction, e.g. for a {\em left-right} Yetter--Drinfeld algebra over $H$, we use a {\em left} action of $H$ and a {\em right} coaction of $H$.
\end{rema}

\subsection{(left-left) Yetter--Drinfeld algebras}\label{sec:ll-YD}

Consider a left action $\rhd : H \odo X \to X$ and a left coaction $\delta: X \to H \odo X$, denoted by $x \mapsto x_{[-1]} \odo x_{[0]}$.

\begin{defi}
The algebra $X$ is called a {\em (left-left) Yetter--Drinfeld algebra over $H$} if 
\begin{equation}\label{eq:ll-YD-1}
h_{(1)}x_{[-1]} \odo (h_{(2)} \rhd x_{[0]}) = (h_{(1)} \rhd x)_{[-1]}h_{(2)} \odo (h_{(1)} \rhd x)_{[0]}
\end{equation}
for all $h \in H$, $x \in X$.
\end{defi}

The expression~\eqref{eq:ll-YD-1} is called the {\em Yetter--Drinfeld condition} for the tuple $(X,\rhd,\delta)$ and sometimes we will denoted $(X,\rhd,\delta) \in {}^{H}_{H}\aYD$ for say that $X$ is a (left-left) Yetter--Drinfeld algebra over $H$. We will give another equivalent condition in the next proposition.

\begin{prop}
The algebra $X$ is a (left-left) Yetter--Drinfeld algebra over $H$ if and only if
\begin{equation}\label{eq:ll-YD-2}
\delta(h \rhd x) = h_{(1)}x_{[-1]}S(h_{(3)}) \odo (h_{(2)} \rhd x_{[0]})
\end{equation}
for all $h \in H$, $x \in X$.
\end{prop}
\pr
\begin{itemize}
\item[($\ida$)] Fix $h \in H$ and $x \in X$, then
\begin{align*}
h_{(1)}x_{[-1]}S(h_{(3)}) \odo (h_{(2)} \rhd x_{[0]}) & = (\bar{\m} \odo \id)(S(h_{(3)}) \odo h_{(1)}x_{[-1]} \odo (h_{(2)} \rhd x_{[0]})) \\
& = (\bar{\m} \odo \id)(S(h_{(3)}) \odo (h_{(1)} \rhd x)_{[-1]}h_{(2)} \odo (h_{(1)} \rhd x)_{[0]}) \\
& = (h_{(1)} \rhd x)_{[-1]}h_{(2)}S(h_{(3)}) \odo (h_{(1)} \rhd x)_{[0]} \\
& = \cou(h_{(2)})\delta(h_{(1)} \rhd x) \\
& = \delta(h \rhd x).
\end{align*}

\item[($\vuelta$)] Fix $h \in H$ and $x \in X$, then
\begin{align*}
(h_{(1)} \rhd x)_{[-1]}h_{(2)} \odo (h_{(1)} \rhd x)_{[0]} & = (\bar{\m} \odo \id)(h_{(2)} \odo \delta(h_{(1)} \rhd x)) \\
& = (\bar{\m} \odo \id)(h_{(2)} \odo h_{(1)_{(1)}}x_{[-1]}S(h_{(1)_{(3)}}) \odo (h_{(1)_{(2)}} \rhd x_{[0]})) \\
& = (\bar{\m} \odo \id)(h_{(3)_{(2)}} \odo h_{(1)}x_{[-1]}S(h_{(3)_{(1)}}) \odo (h_{(2)} \rhd x_{[0]})) \\
& = h_{(1)}x_{[-1]}S(h_{(3)_{(1)}})h_{(3)_{(2)}} \odo (h_{(2)} \rhd x_{[0]}) \\
& = h_{(1)}x_{[-1]} \odo (\cou(h_{(3)})h_{(2)} \rhd x_{[0]}) \\
& = h_{(1)}x_{[-1]} \odo (h_{(2)} \rhd x_{[0]}),
\end{align*}
because
\[
h_{(1)_{(1)}} \odo h_{(1)_{(2)}} \odo h_{(1)_{(3)}} \odo h_{(2)} = h_{(1)} \odo h_{(2)} \odo h_{(3)_{(1)}} \odo h_{(3)_{(2)}}.
\]
\end{itemize}
\fin

\subsubsection*{Braided commutativity}\label{sec:ll-bc}

Let $(X,\rhd,\delta) \in {}^{H}_{H}\aYD$, and consider the map
\[
\begin{array}{lccc}
\tau_{X} : & X \odo X & \to & X \odo X \\ & x \odo y & \mapsto & (x_{[-1]} \rhd y) \odo x_{[0]}
\end{array}
\]

\begin{defi}
A (left-left) Yetter--Drinfeld algebra over $H$, $(X,\rhd,\delta)$, is called {\em braided commutative or $H$-commutative} if $\m_{X}\tau_{X} = \m_{X}$, i.e. if
\begin{equation}\label{eq:ll-bc-1}
xy = (x_{[-1]} \rhd y) x_{[0]}
\end{equation}
for all $x, y \in X$.
\end{defi}

\begin{prop}
A (left-left) Yetter--Drinfeld algebra over $H$, $(X,\rhd,\delta)$, is braided commutative if and only if
\begin{equation}\label{eq:ll-bc-2}
xy = y_{[0]}(S^{-1}(y_{[-1]}) \rhd x)
\end{equation}
for all $x,y \in X$
\end{prop}
\pr
Define the map
\[
\begin{array}{lccc}
\rho_{X} : & X \odo X & \to & X \odo X \\ & x \odo y & \mapsto & y_{[0]} \odo (S^{-1}(y_{[-1]}) \rhd x)
\end{array}
\]
then, $\m_{X}\rho_{X} = \m_{X}$ if and only if $xy = y_{[0]}(S^{-1}(y_{[-1]}) \rhd x)$ for all $x,y \in X$. On the other hand, we have
\begin{align*}
\rho_{X}\tau_{X}(x \odo y) & = \rho_{X}((x_{[-1]} \rhd y) \odo x_{[0]}) = x_{[0]_{[0]}} \odo (S^{-1}(x_{[0]_{[-1]}}) \rhd (x_{[-1]} \rhd y)) \\
& = x_{[0]} \odo (S^{-1}(x_{[-1]_{(2)}})x_{[-1]_{(1)}} \rhd y) = \cou(x_{[-1]})x_{[0]} \odo y = x \odo y
\end{align*}
and
\begin{align*}
\tau_{X}\rho_{X}(x \odo y) & = \tau_{X}(y_{[0]} \odo (S^{-1}(y_{[-1]}) \rhd x)) = (y_{[0]_{[-1]}} \rhd (S^{-1}(y_{[-1]}) \rhd x)) \odo y_{[0]_{[0]}} \\
& = (y_{[-1]_{(2)}}S^{-1}(y_{[-1]_{(1)}}) \rhd x) \odo y_{[0]}  = x \odo \cou(y_{[-1]})y_{[0]} = x \odo y
\end{align*}
for all $x,y \in X$. Then, $\m_{X}\tau_{X} = \m_{X}$ if and only if $\m_{X}\rho_{X} = \m_{X}$.
\fin

\subsection{(left-right) Yetter--Drinfeld algebras}\label{sec:lr-YD}

Consider a left action $\rhd : H \odo X \to X$ and a right coaction $\delta : X \to X \odo H^{\op}$, denoted by $x \mapsto x_{[0]} \odo x^{\op}_{[1]}$.

\begin{defi}
The algebra $X$ is called a {\em (left-right) Yetter--Drinfeld algebra over $H$} if 
\begin{equation}\label{eq:lr-YD-1}
(h_{(1)} \rhd x_{[0]}) \odo h_{(2)}x_{[1]} = (h_{(2)} \rhd x)_{[0]} \odo (h_{(2)} \rhd x)_{[1]}h_{(1)}
\end{equation}
for all $h \in H$, $x \in X$.
\end{defi}

The expression~\eqref{eq:lr-YD-1} is called the {\em Yetter--Drinfeld condition} for the tuple $(X,\rhd,\delta)$ and sometimes we will denoted $(X,\rhd,\delta) \in {}_{H}\aYD^{H}$ for say that $X$ is a (left-right) Yetter--Drinfeld algebra over $H$. We will give another equivalent condition in the next proposition.

\begin{prop}
The algebra $X$ is a (left-right) Yetter--Drinfeld algebra over $H$ if and only if
\begin{equation}
\delta(h \rhd x) = (h_{(2)} \rhd x_{[0]}) \odo (h_{(3)}x_{[1]}S^{-1}(h_{(1)}))^{\op}
\end{equation}
for all $h \in H$, $x \in X$.
\end{prop}
\pr
\begin{itemize}
\item[($\ida$)] Fix $h \in H$ and $x \in X$, then
\begin{align*}
h_{(2)} \rhd x_{[0]} \odo (h_{(3)}x_{[1]}S^{-1}(h_{(1)}))^{\op} & = (\id \odo {}^{\op}\m)(h_{(2)} \rhd x_{[0]} \odo h_{(3)}x_{[1]} \odo S^{-1}(h_{(1)})) \\
& = (\id \odo {}^{\op}\m)((h_{(3)} \rhd x)_{[0]} \odo (h_{(3)} \rhd x)_{[1]}h_{(2)} \odo S^{-1}(h_{(1)})) \\
& = (\id \odo {}^{\op})((h_{(3)} \rhd x)_{[0]} \odo (h_{(3)} \rhd x)_{[1]}h_{(2)}S^{-1}(h_{(1)})) \\
& = (\id \odo {}^{\op})((h_{(2)} \rhd x)_{[0]} \odo (h_{(2)} \rhd x)_{[1]}\cou(h_{(1)})) \\
& = \delta(h_{(2)} \rhd x)\cou(h_{(1)}) \\
& = \delta(h \rhd x).
\end{align*}

\item[($\vuelta$)] Fix $h \in H$ and $x \in X$, then
\begin{align*}
(h_{(2)} \rhd x)_{[0]} \odo (h_{(2)} \rhd x)_{[1]}h_{(1)} & = (\id \odo \m)((\id \odo {}^{\op})(\delta(h_{(2)} \rhd x)) \odo h_{(1)}) \\
& = (\id \odo \m)(h_{(2)_{(2)}} \rhd x_{[0]} \odo h_{(2)_{(3)}}x_{[1]}S^{-1}(h_{(2)_{(1)}}) \odo h_{(1)}) \\
& = h_{(2)} \rhd x_{[0]} \odo h_{(3)}x_{[1]}S^{-1}(h_{(1)_{(2)}})h_{(1)_{(1)}} \\
& = \cou(h_{(1)})h_{(2)} \rhd x_{[0]} \odo h_{(3)}x_{[1]} \\
& = h_{(1)} \rhd x_{[0]} \odo h_{(2)}x_{[1]}.
\end{align*}
\end{itemize}
\fin

\subsubsection*{Braided commutativity}\label{sec:lr-bc}

Let $(X,\rhd,\delta) \in {}_{H}\aYD^{H}$, and consider the map
\[
\begin{array}{lccc}
\tau_{X} : & X \odo X & \to & X \odo X \\ & x \odo y & \mapsto & (S(x_{[1]}) \rhd y) \odo x_{[0]}
\end{array}
\]

\begin{defi}
A (left-right) Yetter--Drinfeld algebra over $H$, $(X,\rhd,\delta)$, is called {\em braided commutative or $H$-commutative} if $\m_{X}\tau_{X} = \m_{X}$, i.e. if
\begin{equation}\label{eq:lr-bc-1}
xy = (S(x_{[1]}) \rhd y) x_{[0]}
\end{equation}
for all $x, y \in X$.
\end{defi}

\begin{prop}
A (left-right) Yetter--Drinfeld algebra over $H$, $(X,\rhd,\delta)$, is braided commutative if and only if
\begin{equation}\label{eq:lr-bc-2}
xy = y_{[0]}(y_{[1]} \rhd x)
\end{equation}
for all $x,y \in X$
\end{prop}
\pr
Define the map
\[
\begin{array}{lccc}
\rho_{X} : & X \odo X & \to & X \odo X \\ & x \odo y & \mapsto & y_{[0]} \odo (y_{[1]} \rhd x)
\end{array}
\]
then, $\m_{X}\rho_{X} = \m_{X}$ if and only if $xy = y_{[0]}(y_{[1]} \rhd x)$ for all $x,y \in X$. On the other hand, we have
\begin{align*}
\rho_{X}\tau_{X}(x \odo y) & = \rho_{X}((S(x_{[1]}) \rhd y) \odo x_{[0]}) = x_{[0]_{[0]}} \odo (x_{[0]_{[1]}} \rhd (S(x_{[1]}) \rhd y)) \\
& = x_{[0]} \odo (x_{[1]_{(1)}}S(x_{[1]_{(2)}}) \rhd y) = x_{[0]}\cou(x_{[1]}) \odo y = x \odo y
\end{align*}
and
\begin{align*}
\tau_{X}\rho_{X}(x \odo y) & = \tau_{X}(y_{[0]} \odo (y_{[1]} \rhd x)) = (S(y_{[0]_{[1]}}) \rhd (y_{[1]} \rhd x)) \odo y_{[0]_{[0]}} \\
& = (S(y_{[1]_{(1)}})y_{[1]_{(2)}} \rhd x) \odo y_{[0]} = x \odo y_{[0]}\cou(y_{[1]}) = x \odo y
\end{align*}
for all $x,y \in X$. Then, $\m_{X}\tau_{X} = \m_{X}$ if and only if $\m_{X}\rho_{X} = \m_{X}$.
\fin

\subsection{(right-right) Yetter--Drinfeld algebras}\label{sec:rr-YD}

Consider a right action $\lhd : X \odo H \to X$ and a right coaction $\delta : X \to X \odo H$, denoted by $x \mapsto x_{[0]} \odo x_{[1]}$.

\begin{defi}
The algebra $X$ is called a {\em (right-right) Yetter--Drinfeld algebra over $H$} if 
\begin{equation}\label{eq:rr-YD-1}
(x_{[0]} \lhd h_{(1)}) \odo x_{[1]}h_{(2)} = (x \lhd h_{(2)})_{[0]} \odo h_{(1)}(x \lhd h_{(2)})_{[1]}
\end{equation}
for all $h \in H$, $x \in X$.
\end{defi}

The expression~\eqref{eq:rr-YD-1} is called the {\em Yetter--Drinfeld condition} for the tuple $(X,\lhd,\delta)$ and sometimes we will denoted $(X,\lhd,\delta) \in \aYD^{H}_{H}$ for say that $X$ is a (right-right) Yetter--Drinfeld algebra over $H$. We will give another equivalent condition in the next proposition.

\begin{prop}
The algebra $X$ is a (right-right) Yetter--Drinfeld algebra over $H$ if and only if
\begin{equation}\label{eq:rr-YD-2}
\delta(x \lhd h) = (x_{[0]} \lhd h_{(2)}) \odo S(h_{(1)})x_{[1]}h_{(3)}
\end{equation}
for all $h \in H$, $x \in X$.
\end{prop}
\pr
\begin{itemize}
\item[($\ida$)] Fix $h \in H$ and $x \in X$, then
\begin{align*}
(x_{[0]} \lhd h_{(2)}) \odo S(h_{(1)})x_{[1]}h_{(3)} & = (\id \odo \bar{\m})((x_{[0]} \lhd h_{(2)}) \odo x_{[1]}h_{(3)} \odo S(h_{(1)})) \\
& = (\id \odo \bar{\m})((x \lhd h_{(3)})_{[0]} \odo h_{(2)}(x \lhd h_{(3)})_{[1]} \odo S(h_{(1)})) \\
& = (x \lhd h_{(3)})_{[0]} \odo S(h_{(1)})h_{(2)}(x \lhd h_{(3)})_{[1]} \\
& = \cou(h_{(1)})\delta(x \lhd h_{(2)}) \\
& = \delta(x \lhd h).
\end{align*}

\item[($\vuelta$)] Fix $h \in H$ and $x \in X$, then
\begin{align*}
(x \lhd h_{(2)})_{[0]} \odo h_{(1)}(x \lhd h_{(2)})_{[1]} & = (\id \odo \bar{\m})(\delta(x \lhd h_{(2)}) \odo h_{(1)}) \\
& = (\id \odo \bar{\m})((x_{[0]} \lhd h_{(2)_{(2)}}) \odo S(h_{(2)_{(1)}})x_{[1]}h_{(2)_{(3)}} \odo h_{(1)}) \\
& = (x_{[0]} \lhd h_{(2)}) \odo h_{(1)_{(1)}}S(h_{(1)_{(2)}})x_{[1]}h_{(3)} \\
& = (x_{[0]} \lhd \cou(h_{(1)})h_{(2)}) \odo x_{[1]}h_{(3)} \\
& = (x_{[0]} \lhd h_{(1)}) \odo x_{[1]}h_{(2)}.
\end{align*}
\end{itemize}
\fin

\subsubsection*{Braided commutativity}\label{sec:rr-bc}

Let $(X,\lhd,\alp) \in \aYD^{H}_{H}$, and consider the map
\[
\begin{array}{lccc}
\tau_{X} : & X \odo X & \to & X \odo X \\ & x \odo y & \mapsto & (y \lhd S^{-1}(x_{[1]})) \odo x_{[0]}
\end{array}
\]

\begin{defi}
A (right-right) Yetter--Drinfeld algebra $(X,\lhd,\delta)$ is called {\em braided commutative} if $\m_{X}\tau_{X} = \m_{X}$, i.e.
\begin{equation}\label{eq:rr-bc-1}
xy = (y \lhd S^{-1}(x_{[1]}))x_{[0]},
\end{equation}
for all $x, y \in X$.
\end{defi}

\begin{prop}
The (right-right) Yetter--Drinfeld algebra over $H$, $(X,\lhd,\delta)$, is braided commutative if and only if
\begin{equation}\label{eq:rr-bc-2}
xy = y_{[0]}(x \lhd y_{[1]}),
\end{equation}
for all $x,y \in X$.
\end{prop}
\pr
Define the map
\[
\begin{array}{lccc}
\rho_{X} : & X \odo X & \to & X \odo X \\ & x \odo y & \mapsto & y_{[0]} \odo (x \lhd y_{[1]})
\end{array}
\]
then, $\m_{X}\rho_{X} = \m_{X}$ if and only if $xy = y_{[0]}(x \lhd y_{[1]})$ for all $x,y \in X$. On the other hand, we have
\begin{align*}
\rho_{X}\tau_{X}(x \odo y) & = \rho_{X}((y \lhd S^{-1}(x_{[1]})) \odo x_{[0]}) = x_{[0]_{[0]}} \odo ((y \lhd S^{-1}(x_{[1]}) \lhd x_{[0]_{[1]}}) \\
& = x_{[0]} \odo (y \lhd S^{-1}(x_{[1]_{(2)}})x_{[1]_{(1)}}) = \cou(x_{[1]})x_{[0]} \odo y = x \odo y
\end{align*}
and
\begin{align*}
\tau_{X}\rho_{X}(x \odo y) & = \tau_{X}(y_{[0]} \odo (x \lhd y_{[1]}))  = ((x \lhd y_{[1]}) \lhd S^{-1}(y_{[0]_{[1]}})) \odo y_{[0]_{[0]}} \\
& = (x \lhd y_{[1]_{(2)}}S^{-1}(y_{[1]_{(1)}})) \odo y_{[0]} = x \odo \cou(y_{[1]})y_{[0]}  = x \odo y
\end{align*}
for all $x,y \in X$. Then, $\m_{X}\tau_{X} = \m_{X}$ if and only if $\m_{X}\rho_{X} = \m_{X}$.
\fin

\subsection{(right-left) Yetter--Drinfeld algebras}\label{sec:rl-YD}

Consider a right action $\lhd : X \odo H \to X$ and a left coaction $\delta : X \to H^{\op} \odo X$, denoted by $x \mapsto x^{\op}_{[-1]} \odo x_{[0]}$.

\begin{defi}
The algebra $X$ is called a {\em (right-left) Yetter--Drinfeld algebra over $H$} if 
\begin{equation}\label{eq:rl-YD-1}
x_{[-1]}h_{(1)} \odo (x_{[0]} \lhd h_{(2)}) = h_{(2)}(x \lhd h_{(1)})_{[-1]} \odo (x \lhd h_{(1)})_{[0]}
\end{equation}
for all $h \in H$, $x \in X$.
\end{defi}

The expression~\eqref{eq:rl-YD-1} is called the {\em Yetter--Drinfeld condition} for the tuple $(X,\lhd,\delta)$ and sometimes we will denoted $(X,\lhd,\delta) \in {}^{H}\aYD_{H}$ for say that $X$ is a (right-left) Yetter--Drinfeld algebra over $H$. We will give another equivalent condition in the next proposition.

\begin{prop}
The algebra $X$ is a (right-left) Yetter--Drinfeld algebra over $H$ if and only if
\begin{equation}\label{eq:rl-YD-2}
\delta(x \lhd h) = (S^{-1}(h_{(3)})x_{[-1]}h_{(1)})^{\op} \odo (x_{[0]} \lhd h_{(2)})
\end{equation}
for all $h \in H$, $x \in X$.
\end{prop}
\pr
\begin{itemize}
\item[($\ida$)] Fix $h \in H$ and $x \in X$, then
\begin{align*}
(S^{-1}(h_{(3)})x_{[-1]}h_{(1)})^{\op} \odo (x_{[0]} \lhd h_{(2)}) & = ({}^{\op}\m \odo \id)(S^{-1}(h_{(3)}) \odo x_{[-1]}h_{(1)} \odo (x_{[0]} \lhd h_{(2)})) \\
& = ({}^{\op}\m \odo \id)(S^{-1}(h_{(3)}) \odo h_{(2)}(x \lhd h_{(1)})_{[-1]} \odo (x \lhd h_{(1)})_{[0]}) \\
& = ({}^{\op} \odo \id)(S^{-1}(h_{(3)})h_{(2)}(x \lhd h_{(1)})_{[-1]} \odo (x \lhd h_{(1)})_{[0]}) \\
& = ({}^{\op} \odo \id)(\cou(h_{(2)})(x \lhd h_{(1)})_{[-1]} \odo (x \lhd h_{(1)})_{[0]}) \\
& = \cou(h_{(2)})\delta(x \lhd h_{(1)}) \\
& = \delta(x \lhd h).
\end{align*}

\item[($\vuelta$)] Fix $h \in H$ and $x \in X$, then
\begin{align*}
h_{(2)}(x \lhd h_{(1)})_{[-1]} \odo (x \lhd h_{(1)})_{[0]} & = (\m \odo \id)(h_{(2)} \odo ({}^{\op} \odo \id)(\delta(x \lhd h_{(1)}))) \\
& = (\m \odo \id)(h_{(2)} \odo S^{-1}(h_{(1)_{(3)}})x_{[-1]}h_{(1)_{(1)}} \odo (x_{[0]} \lhd h_{(1)_{(2)}})) \\
& = h_{(3)_{(2)}}S^{-1}(h_{(3)_{(1)}})x_{[-1]}h_{(1)} \odo (x_{[0]} \lhd h_{(2)}) \\
& = x_{[-1]}h_{(1)} \odo (x_{[0]} \lhd \cou(h_{(3)})h_{(2)}) \\
& = x_{[-1]}h_{(1)} \odo (x_{[0]} \lhd h_{(2)}).
\end{align*}
\end{itemize}
\fin

\subsubsection*{Braided commutativity}\label{sec:rl-bc}

Let $(X,\lhd,\delta) \in {}^{H}\aYD_{H}$, and consider the map
\[
\begin{array}{lccc}
\tau_{X} : & X \odo X & \to & X \odo X \\ & x \odo y & \mapsto & (y \lhd x_{[-1]}) \odo x_{[0]}
\end{array}
\]

\begin{defi}
A (right-left) Yetter--Drinfeld algebra over $H$, $(X,\lhd,\delta)$, is called {\em braided commutative or $H$-commutative} if $\m_{X}\tau_{X} = \m_{X}$, i.e. if
\begin{equation}\label{eq:rl-bc-1}
xy = (y \lhd x_{[-1]}) x_{[0]}
\end{equation}
for all $x, y \in X$.
\end{defi}

\begin{prop}
A (right-left) Yetter--Drinfeld algebra over $H$, $(X,\lhd,\delta)$, is braided commutative if and only if
\begin{equation}\label{eq:rl-bc-2}
xy = y_{[0]}(x \lhd S(y_{[-1]}))
\end{equation}
for all $x,y \in X$
\end{prop}
\pr
Define the map
\[
\begin{array}{lccc}
\rho_{X} : & X \odo X & \to & X \odo X \\ & x \odo y & \mapsto & y_{[0]} \odo (x \lhd S(y_{[-1]}))
\end{array}
\]
then, $\m_{X}\rho_{X} = \m_{X}$ if and only if $xy = y_{[0]}(x \lhd S(y_{[-1]}))$ for all $x,y \in X$. On the other hand, we have
\begin{align*}
\rho_{X}\tau_{X}(x \odo y) & = \rho_{X}((y \lhd x_{[-1]}) \odo x_{[0]}) = x_{[0]_{[0]}} \odo ((y \lhd x_{[-1]}) \lhd S(x_{[0]_{[-1]}})) \\
& = x_{[0]} \odo (y \lhd x_{[-1]_{(1)}}S(x_{[-1]_{(2)}})) = \cou(x_{[-1]})x_{[0]} \odo y = x \odo y
\end{align*}
and
\begin{align*}
\tau_{X}\rho_{X}(x \odo y) & = \tau_{X}(y_{[0]} \odo (x \lhd S(y_{[-1]}))) = ((x \lhd S(y_{[-1]})) \lhd y_{[0]_{[-1]}}) \odo y_{[0]_{[0]}} \\
& = (x \lhd S(y_{[-1]_{(1)}})y_{[-1]_{(2)}}) \odo y_{[0]} = x \odo \cou(y_{[-1]})y_{[0]} = x \odo y
\end{align*}
for all $x,y \in X$. Then, $\m_{X}\tau_{X} = \m_{X}$ if and only if $\m_{X}\rho_{X} = \m_{X}$.
\fin

\section{The ``only coaction'' characterization of Yetter--Drinfeld algebras}\label{sec:only_coaction}

We will give another characterization of Yetter--Drinfeld algebras making use of only coactions of Hopf algebras instead of an action and a coaction as in the standard characterization. The use of this characterization was motivated by the development of a theory of algebraic quantum transformation groupoids \cite{TPHD18} and it is an algebraic adaptation of the definition given in the framework of operator algebraic quantum groups, see \cite{NV10}. Given a Hopf algebra $H$, in each definition of a $(H,\du{H})$-Yetter--Drinfeld algebra, we use the first part in the adjective for the side of the coaction of $H$ and the second part for the side of the coaction of the dual $\du{H}$, for example in the case of a left-right $(H,\du{H})$-Yetter--Drinfeld algebra, we will use a left coaction of $H$ and a right coaction of $\du{H}$.

\subsection{(left-right) Yetter--Drinfeld algebras}

Let $\alp : X \to H \odo X$ a left coaction of $H$ on $X$ and $\beta : X \to X \odo \du{H}$ right coaction of $\du{H}$ on $X$.

\begin{defi}\label{def:new-lr-YD}
We say that the tuple $(X,\alp,\beta)$ is a {\em left-right $(H,\du{H})$-Yetter--Drinfeld algebra} if
\begin{equation}
(\alp \odo \id_{\du{H}})\beta = \ad(U_{13})(\id_{H} \odo \beta)\alp.
\end{equation}
In other words, if the diagram
\[
\xymatrix@C=6pc@R=2pc{X\ar[r]^-{\beta}\ar[d]_-{\alp} & X \odo \du{H}\ar[r]^-{\alp \odo \id_{\du{H}}} & H \odo X \odo \du{H} \\ H \odo X \ar[r]_-{id_{H} \odo \beta} & H \odo X \odo \du{H}\ar[r]_-{\ad(U_{13})} & H \odo X \odo \du{H}\ar@{<=>}[u]_-{id}}
\]
is commutative.
\end{defi}

\begin{nota}
Sometimes we will denoted $(X,\alp,\beta) \in \aYD^{\,lr}(H,\du{H})$ for say that $(X,\alp,\beta)$ is a left-right $(H,\du{H})$-Yetter--Drinfeld algebra.
\end{nota}

\begin{prop}\label{prop:lr-YD-equiv}
$(X,\alp,\beta)$ is a left-right $(H,\du{H})$-Yetter--Drinfeld algebra if and only if $(X,\rhd_{\beta},\alp)$ is a (left-left) Yetter--Drinfeld algebra over $H$, where
\[
\begin{array}{lccc}
\rhd_{\beta} : & H \odo X & \to & X \\ & h \odo x & \mapsto & (\id_{X} \odo {}_{h}\p)(\beta(x))
\end{array}
\]
is the left action of the Hopf algebra $(H,\com)$ on $X$ induced by the coaction $\beta$.
\end{prop}
\pr
Fix $h \in H$ and $x \in X$. It is obvious that $\alp(h \rhd_{\beta} x) =  (\id_{H} \odo \id_{X} \odo {}_{h}\p)(\alp \odo \id_{\du{H}})(\beta(x))$. Now, if we use the notation $y \odo \theta = \beta(x_{[0]})$ then
\begin{align*}
h_{(1)}x_{[-1]}S(h_{(3)}) \odo (h_{(2)} \rhd_{\beta} x_{[0]}) & = h_{(1)}x_{[-1]}S(h_{(3)}) \odo (\id \odo {}_{h_{(2)}}\p)(y \odo \theta) \\
& = h_{(1)}\p(h_{(2)},\theta)x_{[-1]}S(h_{(3)}) \odo y \\
& = (\theta \brhd h_{(1)})x_{[-1]}S(h_{(2)}) \odo y \\
& = (\id \odo {}_{h}\p)(U(x_{[-1]} \odo \theta)U^{-1}) \odo y \\
& = (\id \odo {}_{h}\p \odo \id)\ad(U_{12})(x_{[-1]} \odo \theta \odo y) \\
& = (\id \odo \id \odo {}_{h}\p)\Sigma_{23}\ad(U_{12})\Sigma_{23}(x_{[-1]}\odo y \odo \theta) \\
& = (\id \odo \id \odo {}_{h}\p)\ad(U_{13})(\id \odo \beta)(\alp(x)) 
\end{align*}
This implies that the equality $\alp(h \rhd_{\beta} x) = h_{(1)}x_{[-1]}S(h_{(3)}) \odo (h_{(2)} \rhd_{\beta} x_{[0]})$ holds for all $h \in H$ and $x \in X$ if and only if $(\alp \odo \id_{\du{H}})\beta = \ad(U_{13})(\id_{H} \odo \beta)\alp$.
\fin

\begin{rema}
Consider a left-right $(H,\du{H})$-Yetter--Drinfeld algebra $(X,\alp,\beta)$. The linear map $\delta := \Sigma\beta : X \to \du{H} \odo X$ is a left coaction of $\du{H}^{\co}$ on $X$. Because $\ad(U_{12})\Sigma_{23} = \Sigma_{23}\ad(U_{13})$, we have
\begin{align*}
(\id_{\du{H}} \odo \alp)\delta & = (\id_{\du{H}} \odo \alp)\Sigma\beta = \Sigma_{12}\Sigma_{23}(\alp \odo \id_{\du{H}})\beta \\
& = \Sigma_{12}\Sigma_{23}\ad(U_{13})(\id_{H} \odo \beta)\alp \\
& = (\Sigma \odo \id_{X})(\ad(U) \odo \id_{X})(\id_{H} \odo \delta)\alp.
\end{align*}
Using a pictorial language, we have
\[
\xymatrix@C=6pc@R=2pc{& X \odo \du{H}\ar[d]^-{\Sigma}\ar[r]^-{\alp \odo \id_{\du{H}}} & H \odo X \odo \du{H} \\ X\ar[r]^-{\delta}\ar[d]_-{\alp}\ar[ur]^-{\beta} & \du{H} \odo X\ar[r]^-{id_{\du{H}} \odo \alp} & \du{H} \odo H \odo X\ar[u]_-{\Sigma_{23}\Sigma_{12}} \\ H \odo X\ar[r]_-{id_{H} \odo \delta}\ar[dr]_-{id_{H} \odo \beta} & H \odo \du{H} \odo X\ar[r]_-{\ad(U) \odo \id_{X}} & H \odo \du{H} \odo X\ar[u]_-{\Sigma \odo \id_{X}} \\ & H \odo X \odo \du{H}\ar[u]_-{\Sigma_{23}}\ar[r]_-{\ad(U_{13})} & H \odo X \odo \du{H}\ar[u]_-{\Sigma_{23}}}
\]
which means that
\[
\xymatrix@C=6pc@R=2pc{X\ar[r]^-{\delta}\ar[d]_-{\alp} & \du{H} \odo X\ar[r]^-{id_{\du{H}} \odo \alp} & \du{H} \odo H \odo X \\ H \odo X\ar[r]_-{id_{H} \odo \delta} & H \odo \du{H} \odo X\ar[r]_-{\ad(U) \odo \id_{X}} & H \odo \du{H} \odo X\ar[u]_-{\Sigma \odo \id_{X}}}
\]
is a commutative diagram.
\end{rema}

\subsubsection*{Braided commutativity}\label{sec:lr-bc-mod}

Given a left-right $(H,\du{H})$-Yetter--Drinfeld algebra $(X,\alp,\beta)$, we known that $(X,\rhd_{\beta},\alp)$ is a (left-left) Yetter--Drinfeld algebra over $H$. Observe that the braided commutativity condition~\eqref{eq:ll-bc-1} for $(X,\rhd_{\beta},\alp)$ can be written as the condition $\m_{X}\tau_{X} = \m_{X}$, where $\tau_{X} = (\p \odo \id_{X \odo X})\Sigma_{24}(\alp \odo \beta) : X \odo X \to X \odo X$. Similarly, the condition~\eqref{eq:ll-bc-2} can be written as the condition $\m_{X}\rho_{X} = \m_{X}$, where $\rho_{X} = (\p_{-} \odo \id_{X \odo X})\Sigma_{34}\Sigma_{23}\Sigma_{12}(\beta \odo \alp) : X \odo X \to X \odo X$. One can verified that $\tau_{X}\rho_{X} = \id_{X \odo X} = \rho_{X}\tau_{X}$.

\begin{rema}
We use the notation $\p: H \odo \du{H} \to \ku$ for the canonical pairing $\p(h \odo \theta)= \theta(h)$, and $\p_{-}: H \odo \du{H} \to \ku$ for the linear map $\p_{-}(h \odo \theta)= \theta(S^{-1}(h)) = \du{S}^{-1}(\theta)(h)$.
\end{rema}

Consider the algebra $\He(H) \odo X$ with multiplication given by
\[
\m_{\He(H) \odo X} = (\m_{\He(H)} \odo \m_{X})(\id_{\He(H)} \odo \Sigma_{X,\He(H)} \odo \id_{X}),
\]
here $\He(H)$ is the Heisenberg algebra of $H$ and $\m_{\He(H)}((h \# \theta) \odo (h' \# \theta')) = h(\theta_{(1)} \brhd h') \# \theta_{(2)}\theta'$. Denote by $\iota_{1,X} : H \odo X \to \He(H) \odo X$, $h \odo x \mapsto h \# 1 \odo x$ and by $\iota_{2,X}: \du{H} \odo X \to \He(H) \odo X$, $\theta \odo x \mapsto 1 \# \theta \odo x$, the canonical inclusions induced by the inclusions of $H$ and $\du{H}$ in the Heisenberg algebra $\He(H)$. Note that the canonical pairing can be see as a linear map on $\He(H)$, $\p : H \# \du{H} \to \ku$, $h \# \theta \mapsto \p(h,\theta)$.

\begin{defi}\label{df:lr-bc-mod}
We will say that a left-right $(H,\du{H})$-Yetter--Drinfeld algebra $(X,\alp,\beta)$ is {\em braided commutative} if
\begin{equation}\label{eq:lr-bc-mod}
\m_{\He(H) \odo X}(\iota_{1,X} \odo \iota_{2,X})(\alp(x) \odo \Sigma\beta(y)) = \m_{\He(H) \odo X}(\iota_{2,X} \odo \iota_{1,X})(\Sigma\beta(y) \odo \alp(x))
\end{equation}
for every $x,y \in X$.
\end{defi}

\begin{nota}
We will denote
\[
\aYD^{\,lr}_{bc}(H,\du{H}) = \bigl\{ (X,\alp,\beta) \in \aYD^{\,lr}(H,\du{H}) \;:\; (X,\alp,\beta) \text{ is braided commutative} \bigr\}
\]
\end{nota}

The braided commutativity condition~\eqref{eq:lr-bc-mod} for Yetter--Drinfeld algebras is equivalent to the braided commutativity condition in the standard characterization of Yetter--Drinfeld algebras as we prove in the next proposition. 

\begin{prop}\label{prop:lr-bc-equiv}
The left-right $(H,\du{H})$-Yetter--Drinfeld algebra $(X,\alp,\beta)$ is braided commutative if and only if the (left-left) Yetter--Drinfeld algebra over $H$, $(X,\rhd_{\beta},\alp)$, is braided commutative.
\end{prop}
\pr
Given $x,y \in X$, using the notations $x_{[-1]} \odo x_{[0]} := \alp(x)$ and $y^{[0]} \odo y^{[1]} := \beta(y)$, we have
\begin{align*}
\m_{\He(H) \odo X}(\iota_{1,X} \odo \iota_{2,X})(\alp(x) \odo \Sigma\beta(y)) & = \m_{\He(H) \odo X}(x_{[-1]} \# 1 \odo x_{[0]} \odo 1 \# y^{[1]} \odo y^{[0]}) \\
& = (\id_{\He(H)} \odo \m_{X})(x_{[-1]} \# y^{[1]} \odo x_{[0]} \odo y^{[0]})
\intertext{and}
\m_{\He(H) \odo X}(\iota_{2,X} \odo \iota_{1,X})(\Sigma\beta(y) \odo \alp(x)) & = \m_{\He(H) \odo X}(1 \# y^{[1]} \odo y^{[0]} \odo x_{[-1]} \# 1 \odo x_{[0]}) \\
& = ({y^{[1]}_{(1)}} \brhd x_{[-1]}) \# {y^{[1]}_{(2)}} \odo y^{[0]}x_{[0]} \\
& = x_{[-1]_{(1)}} \# {y^{[1]}_{(2)}} \odo \p(x_{[-1]_{(2)}} \odo {y^{[1]}_{(1)}}) y^{[0]}x_{[0]} \\
& = x_{[-1]} \# y^{[1]} \odo \p(x_{[0]_{[-1]}} \odo y^{[0]^{[1]}}) y^{[0]^{[0]}}x_{[0]_{[0]}} \\
& = x_{[-1]} \# y^{[1]} \odo (\p \odo \m_{X})\Sigma_{24}(\alp \odo \beta)(x_{[0]} \odo y^{[0]}) \\
& = (\id_{\He(H)} \odo \m_{X}\tau_{X})(x_{[-1]} \# y^{[1]} \odo x_{[0]} \odo y^{[0]}).
\end{align*}

\begin{itemize}
\item[($\vuelta$)] If $(X,\lhd_{\beta},\alp)$ is braided commutative, we have $\m_{X} = \m_{X}\tau_{X}$, then follows immediately from the above equalities that the condition~\eqref{eq:lr-bc-mod} is true for every $x,y \in X$, i.e. the tuple $(X,\alp,\beta)$ is braided commutative.

\item[($\ida$)] Let $x',y' \in X$. Take $x \odo y = \rho_{X}\Sigma(x' \odo y')$, thus we have $x' \odo y' = \Sigma\tau_{X}(x \odo y)$. Because $(X,\alp,\beta)$ satisfies the condition~\eqref{eq:lr-bc-mod} for $x,y \in X$, the equalities
\begin{align*}
(\p \odo \id_{X})\m_{\He(H) \odo X}(i_{1,X} \odo i_{2,X})(\alp(x) \odo \Sigma\beta(y)) & = (\p \odo \m_{X})(x_{[-1]} \# y^{[1]} \odo x_{[0]} \odo y^{[0]}) \\
& = \m_{X}(x_{[0]} \odo \p(x_{[-1]} \# y^{[1]})y^{[0]}) \\
& = \m_{X}\Sigma\tau_{X}(x \odo y) \\
& = \m_{X}(x' \odo y'), \\
(\p \odo \id_{X})\m_{\He(H) \odo X}(i_{2,X} \odo i_{1,X})(\Sigma\beta(y) \odo \alp(x)) & = (\p \odo \m_{X}\tau_{X})(x_{[-1]} \# y^{[1]} \odo x_{[0]} \odo y^{[0]}) \\
& = \m_{X}\tau_{X}\Sigma\tau_{X}(x \odo y) \\
& = \m_{X}\tau_{X}(x' \odo y'),
\end{align*}
imply that $\m_{X}(x' \odo y') = \m_{X}\tau_{X}(x' \odo y')$. Thus $(X,\lhd_{\beta},\alp)$ is braided commutative.
\end{itemize}
\fin

\begin{rema}
It holds
\begin{enumerate}[label=\textup{(\roman*)}]
\item $\alp: X \to H \odo X$ is a left coaction of $H$ on $X$ if and only if $\Sigma\alp^{\co}: X^{\op} \to X^{\op} \odo H$ is a right coaction of $H$ on $X^{\op}$.
\item $\beta: X \to X \odo \du{H}$ is a right coaction of $\du{H}$ on $X$ if and only if  $(\Sigma\beta)^{\co}: X^{\op} \to \du{H} \odo X^{\op}$ is a left coaction of $\du{H}$ on $X^{\op}$.
\end{enumerate}
\end{rema}

With the last remark, we have the following propositions

\begin{prop}\label{prop:dual-lr-bc-mod}
A left-right $(H,\du{H})$-Yetter--Drinfeld algebra $(X,\alp,\beta)$ is braided commutative if and only if
\small
\begin{align*}
\m_{\He(\du{H}) \odo X^{\op}}(\du{\iota}_{1,X^{\op}} \odo \du{\iota}_{2,X^{\op}})((\Sigma\beta)^{\co}(y^{\op}) \odo \alp^{\co}(x^{\op})) = \m_{\He(\du{H}) \odo X^{\op}}(\du{\iota}_{2,X^{\op}} \odo \du{\iota}_{1,X^{\op}})(\alp^{\co}(x^{\op}) \odo (\Sigma\beta)^{\co}(y^{\op}))
\end{align*}
\normalsize
for every $x,y \in X$.
\end{prop}
\pr
Consider the anti-isomorphism $\mathcal{L}_{1}: \He(\du{H}) \to \He(H)$, see Proposition~\ref{prop:anti_heisen}. Observe that
\[
(\mathcal{L}_{1} \odo {}^{\op})\hat\iota_{1,X^{\op}}(\Sigma\beta)^{\co}(x^{\op}) =  \iota_{2,X}\Sigma\beta(x) \quad \text{ and } \quad (\mathcal{L}_{1} \odo {}^{\op})\hat\iota_{2,X^{\op}}\alp^{\co}(x^{\op}) = \iota_{1,X}\alp(x)
\]
for all $x \in X$. Because $(\mathcal{L}_{1} \odo {}^{\op})$ is an anti-isomorphism between $\He(\du{H}) \odo X^{\op}$ and $\He(H) \odo X$, it hold
\small
\begin{align*}
(\mathcal{L}_{1} \odo {}^{\op})\m_{\He(\du{H}) \odo X^{\op}}(\hat\iota_{1,X^{\op}} \odo \hat\iota_{2,X^{\op}})((\Sigma\beta)^{\co}(y^{\op}) \odo \alp^{\co}(x^{\op})) = \m_{\He(H) \odo X}(\iota_{1,X} \odo \iota_{2,X})(\alp(x) \odo \Sigma\beta(y))
\end{align*}
\normalsize
and
\small
\[
(\mathcal{L}_{1} \odo {}^{\op})\m_{\He(\du{H}) \odo X^{\op}}(\du{\iota}_{2,X^{\op}} \odo \du{\iota}_{1,X^{\op}})(\alp^{\co}(x^{\op}) \odo (\Sigma\beta)^{\co}(y^{\op})) = \m_{\He(H) \odo X}(\iota_{2,X} \odo \iota_{1,X})(\Sigma\beta(y) \odo \alp(x))
\]
\normalsize
for any $x, y \in X$. The result follows from those equalities.
\fin

\begin{prop}
We have
\[
(X,\alp,\beta) \in \aYD^{\,lr}(H,\du{H}) \qquad\text{ if and only if }\qquad (X^{\op},(\Sigma\beta)^{\co},\Sigma\alp^{\co})  \in \aYD^{\,lr}(\du{H},H).
\]
Moreover, $(X,\alp,\beta)$ is braided commutative if and only if $(X^{\op},(\Sigma\beta)^{\co},\Sigma\alp^{\co})$ is braided commutative.
\end{prop}
\pr
First, note that
\[
((\Sigma\beta)^{\co} \odo \id_{H})\circ\Sigma\alp^{\co} = \Sigma_{23}\circ\Sigma_{12}\circ\Sigma_{23}\circ(S \odo {}^{\op} \odo \dS^{-1})\circ(\id \odo \beta)\circ\alp\circ{}^{\op}
\]
and
\begin{align*}
(\id_{\du{H}} \odo \Sigma\alp^{\co})\circ(\Sigma\beta)^{\co} & = \Sigma_{23}\circ(\dS^{-1} \odo S \odo {}^{\op})\circ\Sigma_{12}\circ\Sigma_{23}\circ(\alp \odo \id)\circ\beta\circ{}^{\op} \\
& = \Sigma_{23}\circ\Sigma_{12}\circ\Sigma_{23}\circ(S \odo {}^{\op} \odo \dS^{-1})\circ(\alp \odo \id)\circ\beta\circ{}^{\op}.
\end{align*}
Because, $(S \odo {}^{\op} \odo \dS^{-1})\ad(U_{13}) = \ad(U^{-1}_{13})(S \odo {}^{\op} \odo \dS^{-1})$, then the equality
\[
(\alp \odo \id_{\du{H}})\beta = \ad(U_{13})(\id_{H} \odo \beta)\alp
\]
is equivalent to the equality $\ad(\Sigma(U)_{13})(\id_{\du{H}} \odo \Sigma\alp^{\co})(\Sigma\beta)^{\co} = ((\Sigma\beta)^{\co} \odo \id_{H})\Sigma\alp^{\co}$. This implies, $(X,\alp,\beta) \in \aYD^{\,lr}(H,\du{H})$ if and only if $(X^{\op},(\Sigma\beta)^{\co},\Sigma\alp^{\co})  \in \aYD^{\,lr}(\du{H},H)$. Follows from Definition~\ref{df:lr-bc-mod} and Proposition~\ref{prop:dual-lr-bc-mod} the equivalence about the braided commutativity conditions.
\fin

\subsection{(right-right) Yetter--Drinfeld algebras}

Let $\alp : X \to X \odo H^{\op}$ be a right coaction of $H^{\op}$ on $X$ and $\beta : X \to X \odo \du{H}$ be right coaction of $\du{H}$ on $X$.

\begin{defi}\label{def:new-rr-YD}
We say that the tuple $(X,\alp,\beta)$ is a {\em right-right $(H,\du{H})$-Yetter--Drinfeld algebra} if
\begin{equation}
(\alp \odo \id_{\du{H}})\beta = (\id_{X} \odo \ad({}^{\ops}U)^{-1})(\id_{X} \odo \Sigma)(\beta \odo \id_{H^{\op}})\alp.
\end{equation}
In other words, if the diagram
\[
\xymatrix@C=6pc@R=2pc{X\ar[r]^-{\beta}\ar[d]_-{\alp} & X \odo \du{H}\ar[r]^-{\alp \odo \id_{\du{H}}} & X \odo H^{\op} \odo \du{H} \\ X \odo H^{\op}\ar[r]_-{\beta \odo \id_{H}} & X \odo \du{H} \odo H^{\op}\ar[r]_-{id_{X} \odo \Sigma} & X \odo H^{\op} \odo \du{H}\ar[u]_-{id_{X} \odo \ad({}^{\ops}U)^{-1}}}
\]
is commutative.
\end{defi}

\begin{nota}
Sometimes we will denoted $(X,\alp,\beta) \in \aYD^{\,rr}(H,\du{H})$ for say that $(X,\alp,\beta)$ is a right-right $(H,\du{H})$-Yetter--Drinfeld algebra.
\end{nota}

\begin{prop}\label{prop:rr-YD-equiv}
$(X,\alp,\beta)$ is a right-right $(H,\du{H})$-Yetter--Drinfeld algebra if and only if $(X,\rhd_{\beta},\alp)$ is a (left-right) Yetter--Drinfeld algebra over $H$. Here
\[
\begin{array}{lccc}
\rhd_{\beta} : & H \odo X & \to & X \\ & h \odo x & \mapsto & (\id \odo {}_{h}\p)(\beta(x))
\end{array}
\]
is the left action of the Hopf algebra $(H,\com)$ on $X$ induced by the right coaction $\beta$.
\end{prop}
\pr
Fix $h \in H$ and $x \in X$. It is obvious that $\alp(h \rhd_{\beta} x) = (\id_{X} \odo \id_{H^{\op}} \odo{}_{h}\p)(\alp \odo \id_{\du{H}})(\beta(x))$. Now, if we use the notation $y \odo \theta = \beta(x_{[0]})$ then
\begin{align*}
(h_{(2)} \rhd_{\beta} x_{[0]}) \odo h_{(3)}x_{[1]}S^{-1}(h_{(1)}) & = (\id_{X} \odo {}_{h_{(2)}}\p)(y \odo \theta) \odo h_{(3)}x_{[1]}S^{-1}(h_{(1)}) \\
& = y \odo \p(h_{(2)},\theta)h_{(3)}x_{[1]}S^{-1}(h_{(1)}) \\
& = y \odo (h_{(2)} \blhd \theta)x_{[1]}S^{-1}(h_{(1)}) \\
& = y \odo (\id_{H} \odo {}_{h}\p)(T_{U^{-1}}(x_{[1]} \odo \theta)) \\
& = (\id_{X} \odo \id_{H} \odo {}_{h}\p)(\id_{X} \odo T_{U^{-1}}\Sigma)(y \odo \theta \odo x_{[1]})\\
& = (\id_{X} \odo \id_{H} \odo {}_{h}\p)(\id_{X} \odo T_{U^{-1}}\Sigma(\id_{\du{H}} \odo {}^{\op}))(\beta \odo \id_{H^{\op}})(\alp(x)).
\end{align*}
This implies that the equality $\alp(h \rhd_{\beta} x) = (h_{(2)} \rhd_{\beta} x_{[0]}) \odo (h_{(3)}x_{[1]}S^{-1}(h_{(1)}))^{\op}$ holds for all $h \in H$ and $x \in X$ if and only if
\[
(\alp \odo \id_{\du{H}})\beta = (\id_{X} \odo ({}^{\op} \odo \id_{\du{H}}) T_{U^{-1}}({}^{\op} \odo \id_{\du{H}})\Sigma)(\beta \odo \id_{H^{\op}})\alp.
\]
Finally, the proposition holds, because $({}^{\op} \odo \id_{\du{H}}) T_{U^{-1}}({}^{\op} \odo \id_{\du{H}}) = \ad({}^{\ops}U)^{-1}$.
\fin

\subsubsection*{Braided commutativity}\label{sec:rr-bc-mod}

Given a right-right $(H,\du{H})$-Yetter--Drinfeld algebra $(X,\alp,\beta)$, we known that $(X,\rhd_{\beta},\alp)$ is a (left-right) Yetter--Drinfeld algebra over $H$. Observe that the braided commutativity condition~\eqref{eq:lr-bc-1} for $(X,\rhd_{\beta},\alp)$ can be written as the condition $\m_{X}\tau_{X} = \m_{X}$, where $\tau_{X} = (\p_{+} \odo \id_{X \odo X})\Sigma_{24}\Sigma_{12}(\dot\alp \odo \beta) : X \odo X \to X \odo X$. Similarly, the condition~\eqref{eq:lr-bc-2} can be written as the condition $\m_{X}\rho_{X} = \m_{X}$, where $\rho_{X} = (\p \odo \id_{X \odo X})\Sigma_{14}(\beta \odo \dot\alp) : X \odo X \to X \odo X$. One can verified that $\tau_{X}\rho_{X} = \id_{X \odo X} = \rho_{X}\tau_{X}$.

\begin{rema}
We use the notation $\p: H \odo \du{H} \to \ku$ for the canonical pairing $\p(h \odo \theta)= \theta(h)$, and $\p_{+}: H \odo \du{H} \to \ku$ for the linear map $\p_{+}(h \odo \theta)= \theta(S(h)) = \du{S}(\theta)(h)$.
\end{rema}


Consider the algebra $X \odo \He(H)$ with multiplication given by
\[
\m_{X \odo \He(H)} = (\m_{X} \odo \m_{\He(H)})(\id_{X} \odo \Sigma_{\He(H),X} \odo \id_{\He(H)}),
\]
here $\He(H) := H\#\du{H}$ is the Heisenberg algebra of $H$. Denote by $\iota_{X,1} : H \odo X \to \He(H) \odo X$, $h \odo x \mapsto h \# 1 \odo x$ and by $\iota_{X,2}: \du{H} \odo X \to \He(H) \odo X$, $\theta \odo x \mapsto 1 \# \theta \odo x$, the canonical inclusions induced by the inclusions of $H$ and $\du{H}$ in the Heisenberg algebra $\He(H)$. Note that the canonical pairing $\p$ can be see as a linear map on $\He(H)$, $\p : H \# \du{H} \to \ku$, $h \# \theta \mapsto \p(h,\theta)$.

\begin{defi}\label{df:rr-bc-mod}
We will say that a right-right $(H,\du{H})$-Yetter--Drinfeld algebra $(X,\alp,\beta)$ is {\em braided commutative} if
\begin{equation}\label{eq:rr-bc-mod}
\m_{X \odo \He(H)}(\iota_{X,1} \odo \iota_{X,2})((\id \odo S)\dot\alp(x) \odo \beta(y)) = \m_{X \odo \He(H)}(\iota_{X,2} \odo \iota_{X,1})(\beta(y) \odo (\id \odo S)\dot\alp(x))
\end{equation}
for every $x,y \in X$.
\end{defi}

\begin{nota}
We will denote
\[
\aYD^{\,rr}_{bc}(H,\du{H}) = \bigl\{ (X,\alp,\beta) \in \aYD^{\,rr}(H,\du{H}) \;:\; (X,\alp,\beta) \text{ is braided commutative} \bigr\}
\]
\end{nota}

The braided commutativity condition~\eqref{eq:rr-bc-mod} for Yetter--Drinfeld algebras is equivalent to the braided commutativity condition in the standard characterization of Yetter--Drinfeld algebras as we prove in the next proposition. 

\begin{prop}\label{prop:bc_rr_equivalence}
The right-right $(H,\du{H})$-Yetter--Drinfeld algebra $(X,\alp,\beta)$ is braided commutative if and only if the (left-right) Yetter--Drinfeld algebra over $H$, $(X,\rhd_{\beta},\alp)$, is braided commutative.
\end{prop}
\pr
Given $x,y \in X$, using the notations $x_{[0]} \odo x^{\op}_{[1]} := \alp(x)$ and $y^{[0]} \odo y^{[1]} := \beta(y)$, we have
\small
\begin{align*}
\m_{X \odo \He(H)}(\iota_{X,1} \odo \iota_{X,2})((\id \odo S)\dot\alp(x) \odo \beta(y)) & = \m_{X \odo \He(H)}(x_{[0]} \odo S(x_{[1]}) \# 1 \odo y^{[0]} \odo 1 \# y^{[1]}) \\
& = (\m_{X} \odo \id_{\He(H)})(x_{[0]} \odo y^{[0]} \odo S(x_{[1]}) \# y^{[1]})
\intertext{\normalsize and}
\m_{X \odo \He(H)}(\iota_{X,2} \odo \iota_{X,1})(\beta(y) \odo (\id \odo S)\dot\alp(x)) & = \m_{X \odo \He(H)}(y^{[0]} \odo 1 \# y^{[1]} \odo x_{[0]} \odo S(x_{[-1]}) \# 1) \\
& = y^{[0]}x_{[0]} \odo ({y^{[1]}_{(1)}} \brhd S(x_{[1]})) \# {y^{[1]}_{(2)}} \\
& = \p(S(x_{[1]_{(1)}}) \odo {y^{[1]}_{(1)}}) y^{[0]}x_{[0]} \odo S(x_{[1]_{(2)}}) \# {y^{[1]}_{(2)}} \\
& = \p_{+}(x_{[0]_{[1]}} \odo y^{[0]^{[1]}}) y^{[0]^{[0]}}x_{[0]_{[0]}} \odo S(x_{[1]}) \# y^{[1]}\\
& = (\p_{+} \odo \m_{X})(x_{[0]_{[1]}} \odo y^{[0]^{[1]}} \odo y^{[0]^{[0]}} \odo x_{[0]_{[0]}}) \odo S(x_{[1]}) \# y^{[1]}\\
& = (\p_{+} \odo \m_{X})\Sigma_{24}(x_{[0]_{[1]}} \odo x_{[0]_{[0]}} \odo y^{[0]^{[0]}} \odo y^{[0]^{[1]}}) \odo S(x_{[1]}) \# y^{[1]}\\
& = (\p_{+} \odo \m_{X})\Sigma_{24}\Sigma_{12}(\dot\alp \odo \beta)(x_{[0]} \odo y^{[0]}) \odo S(x_{[1]}) \# y^{[1]} \\
& = \m_{X}(\p_{+} \odo \id_{X \odo X})\Sigma_{24}\Sigma_{12}(\dot\alp \odo \beta)(x_{[0]} \odo y^{[0]}) \odo S(x_{[1]}) \# y^{[1]}\\
& = \m_{X}\tau_{X}(x_{[0]} \odo y^{[0]}) \odo S(x_{[1]}) \# y^{[1]} \\
& = (\m_{X}\tau_{X} \odo \id_{\He(H)})(x_{[0]} \odo y^{[0]} \odo S(x_{[1]}) \# y^{[1]}).
\end{align*}
\normalsize

\begin{itemize}
\item[($\vuelta$)] If $(X,\rhd_{\beta},\alp)$ is braided commutative, we have $\m_{X} = \m_{X}\tau_{X}$, then follows immediately from the above equalities that the condition~\eqref{eq:rr-bc-mod} is true for every $x,y \in X$, i.e. the tuple $(X,\alp,\beta)$ is braided commutative.

\item[($\ida$)] Let $x',y' \in X$. Take $x \odo y = \rho_{X}\Sigma(x' \odo y')$, thus we have $x' \odo y' = \Sigma\tau_{X}(x \odo y)$. Because $(X,\alp,\beta)$ satisfies the condition~\eqref{eq:rr-bc-mod} for $x,y \in X$, the equalities
\begin{align*}
(\id_{X} \odo \p)\m_{X \odo \He(H)}(\iota_{X,1} \odo \iota_{X,2})((\id \odo S)\dot\alp(x) \odo \beta(y)) & = (\m_{X} \odo \p)(x_{[0]} \odo y^{[0]} \odo S(x_{[1]}) \# y^{[1]}) \\
& = \m_{X}(x_{[0]} \odo \p(S(x_{[1]}) \# y^{[1]})y^{[0]}) \\
& = \m_{X}(x_{[0]} \odo (S(x_{[1]}) \rhd_{\beta} y)) \\
& = \m_{X}\Sigma\tau_{X}(x \odo y) \\
& = \m_{X}(x' \odo y'), \\
(\id_{X} \odo \p)\m_{X \odo \He(H)}(\iota_{X,1} \odo \iota_{X,2})((\id \odo S)\dot\alp(x) \odo \beta(y)) & = (\m_{X}\tau_{X} \odo \p)(x_{[0]} \odo y^{[0]} \odo S(x_{[1]}) \# y^{[1]}) \\
& = \m_{X}\tau_{X}\Sigma\tau_{X}(x \odo y) \\
& = \m_{X}\tau_{X}(x' \odo y'),
\end{align*}
imply that $\m_{X}(x' \odo y') = \m_{X}\tau_{X}(x' \odo y')$. Thus $(X,\rhd_{\beta},\alp)$ is braided commutative.
\end{itemize}
\fin

\begin{rema}
It is obvious that
\begin{enumerate}[label=\textup{(\roman*)}]
\item $\alp: X \to X \odo H^{\op}$ is a right coaction of $H^{\op}$ on $X$ if and only if $\alp^{\op}=(\id \odo S)\alp^{\co}: X^{\op} \to X^{\op} \odo H$ is a right coaction of $H$ on $X^{\op}$.
\item $\beta: X \to X \odo \du{H}$ is a right coaction of $\du{H}$ on $X$ if and only if $\beta^{\op}=(\id \odo \dS^{-1})\beta^{\co}: X^{\op} \to X^{\op} \odo \du{H}^{\op}$ is a right coaction of $\du{H}^{\op}$ on $X^{\op}$.
\end{enumerate}
\end{rema}

With the last remark, we have the following propositions

\begin{prop}\label{prop:dual-rr-bc-mod}
A right-right $(H,\du{H})$-Yetter--Drinfeld algebra $(X,\alp,\beta)$ is braided commutative if and only if
\begin{multline*}
\m_{X^{\op} \odo \He(\du{H})}(\hat\iota_{X^{\op},1} \odo \hat\iota_{X^{\op},2})((\id \odo \dS)\dot\beta^{\op}(y^{\op}) \odo \alp^{\op}(x^{\op})) \\ = \m_{X^{\op} \odo \He(\du{H})}(\hat\iota_{X^{\op},2} \odo \hat\iota_{X^{\op},1})(\alp^{\op}(x^{\op}) \odo (\id \odo \dS)\dot\beta^{\op}(y^{\op}))
\end{multline*}
for every $x,y \in X$.
\end{prop}
\pr
Consider the anti-isomorphism $\mathcal{L}_{2}: \He(\du{H}) \to \He(H)$, see Proposition~\ref{prop:anti_heisen}. Obverse that
\[
({}^{\op} \odo \mathcal{L}_{2})\hat\iota_{X^{\op},1}(\id \odo \dS)\dot\beta^{\op}(x^{\op}) = \iota_{X,2}\beta(x) \quad \text { and } \quad ({}^{\op} \odo \mathcal{L}_{2})\hat\iota_{X^{\op},2}\alp^{\op}(x^{\op}) = \iota_{X,1}(\id \odo S)\dot\alp(x)
\]
for all $x \in X$. Because $({}^{\op} \odo \mathcal{L}_{2})$ is an anti-isomorphism between $X^{\op} \odo \He(\du{H})$ and $X \odo \He(H)$, it hold
\footnotesize
\begin{align*}
({}^{\op} \odo \mathcal{L}_{2})\m_{X^{\op} \odo \He(\du{H})}(\hat\iota_{X^{\op},1} \odo \hat\iota_{X^{\op},2})((\id \odo \dS)\dot\beta^{\op}(y^{\op}) \odo \alp^{\op}(x^{\op})) = \m_{X \odo \He(H)}(\iota_{X,1} \odo \iota_{X,2})((\id \odo S)\dot\alp(x) \odo \beta(y))
\end{align*}
\normalsize
and similarly
\footnotesize
\[
({}^{\op} \odo \mathcal{L}_{2})\m_{X^{\op} \odo \He(\du{H})}(\du{\iota}_{X^{\op},2} \odo \du{\iota}_{X^{\op},1})(\alp^{\op}(x^{\op}) \odo (\id \odo \dS)\beta^{\op}(y^{\op})) = \m_{X \odo \He(H)}(\iota_{X,2} \odo \iota_{X,1})(\beta(y) \odo (\id \odo S)\dot\alp(x))
\]
\normalsize
for all $x, y \in X$. The result follows those equalities.
\fin

\begin{prop}
We have
\[
(X,\alp,\beta) \in \aYD^{\,rr}(H,\du{H}) \qquad\text{ if and only if }\qquad (X^{\op},\beta^{\op},\alp^{\op})  \in \aYD^{\,rr}(\du{H},H).
\]
Moreover, $(X,\alp,\beta)$ is braided commutative if and only if $(X^{\op},\beta^{\op},\alp^{\op})$ is braided commutative.
\end{prop}
\pr
First, note that $(\beta^{\op} \odo \id_{H})\circ\alp^{\op} = ({}^{\op} \odo {}^{\op} \odo {}^{\op})\circ(\beta \odo \id_{H^{\op}})\circ\alp\circ{}^{\op}$ and similarly
\begin{align*}
(\alp^{\op} \odo \id_{\du{H}^{\op}})\circ\beta^{\op} & = ({}^{\op} \odo {}^{\op} \odo {}^{\op})\circ(\alp \odo \id_{\du{H}})\circ\beta\circ{}^{\op}
\end{align*}
Because $({}^{\op} \odo {}^{\op})\Sigma\ad({}^{\ops}U) = \ad({}^{\ops}\Sigma(U)^{-1})\Sigma({}^{\op} \odo {}^{\op})$, then the equality
\[
(\beta^{\op} \odo \id_{H})\alp^{\op} = (\id_{X} \odo \ad({}^{\ops}\Sigma(U))^{-1})(\id_{X} \odo \Sigma)(\alp^{\op} \odo \id_{\du{H}^{\op}})\alp^{\ops}
\]
is equivalent to the equality $(\alp \odo \id_{\du{H}})\beta = (\id_{X} \odo \ad({}^{\ops}U)^{-1})(\id_{X} \odo \Sigma)(\beta \odo \id_{H^{\op}})\alp$. This implies, $(X,\alp,\beta) \in \aYD^{\,rr}(H,\du{H})$ if and only if $(X^{\op},\beta^{\op},\alp^{\op})  \in \aYD^{\,rr}(\du{H},H)$. Follows from Definition~\ref{df:lr-bc-mod} and Proposition~\ref{prop:dual-rr-bc-mod} the equivalence about the braided commutativity conditions.
\fin

\subsection{(right-left) Yetter--Drinfeld algebras}

Let $\alp : X \to X \odo H$ be a right coaction of $H$ on $X$ and $\beta : X \to \du{H} \odo X$ a left coaction of $\du{H}$ on $X$.

\begin{defi}\label{def:new-rl-YD}
We say that the tuple $(X,\alp,\beta)$ is a {\em right-left $(H,\du{H})$-Yetter--Drinfeld algebra} if
\begin{equation}
(\id_{\du{H}} \odo \alp)\beta = \ad(\Sigma(U)^{-1}_{13})(\beta \odo \id_{H})\alp.
\end{equation}
In other words, if the diagram
\[
\xymatrix@C=6pc@R=2pc{X\ar[r]^-{\beta}\ar[d]_-{\alp} & \du{H}\ar[r]^-{id_{\du{H}} \odo \alp} \odo X & \du{H} \odo X \odo H \\ X \odo H \ar[r]_-{\beta \odo \id_{H}} & \du{H} \odo X \odo H \ar[r]_-{\ad(\Sigma(U)^{-1}_{13})} & \du{H} \odo X \odo H \ar@{<=>}[u]_-{id}}
\]
is commutative.
\end{defi}

\begin{nota}
Sometimes we will denoted $(X,\alp,\beta) \in \aYD^{\,rl}(H,\du{H})$ for say that $(X,\alp,\beta)$ is a right-left $(H,\du{H})$-Yetter--Drinfeld algebra.
\end{nota}

\begin{prop}\label{prop:rl-YD-equiv}
$(X,\alp,\beta)$ is a right-left $(H,\du{H})$-Yetter--Drinfeld algebra if and only if $(X,\lhd_{\beta},\alp)$ is a (right-right) Yetter--Drinfeld algebra over $H$, where
\[
\begin{array}{lccc}
\lhd_{\beta} : & X \odo H & \to & X \\ & x \odo h & \mapsto & ({}_{h}\p \odo \id_{X})(\beta(x))
\end{array}
\]
is the right action of the Hopf algebra $(H,\com)$ on $X$ induced by the left coaction $\beta$.
\end{prop}
\pr
Fix $h \in H$ and $x \in X$. It is obvious that $\alp(x \lhd_{\beta} h) =  ({}_{h}\p \odo \id_{X} \odo \id_{H})(\id_{\du{H}} \odo \alp)(\beta(x))$. Now, if we use the notation $\theta \odo y = \beta(x_{[0]})$ then
\begin{align*}
(x_{[0]} \lhd_{\beta} h_{(2)}) \odo S(h_{(1)})x_{[1]}h_{(3)} & = ({}_{h_{(2)}}\p \odo \id_{X})(\theta \odo y) \odo S(h_{(1)})x_{[1]}h_{(3)} \\
& = y \odo S(h_{(1)})x_{[1]}\p(h_{(2)},\theta)h_{(3)} \\
& = y \odo S(h_{(1)})x_{[1]}(h_{(2)} \blhd \theta) \\
& = y \odo (\id \odo {}_{h}\p)(U^{-1}(x_{[1]} \odo \theta)U) \\
& = y \odo ({}_{h}\p \odo \id)(\Sigma(U)^{-1}(\theta \odo x_{[1]})\Sigma(U)) \\
& = (\id \odo {}_{h}\p \odo \id)\ad(\Sigma(U)^{-1}_{23})(y \odo \theta \odo x_{[1]}) \\
& = ({}_{h}\p \odo \id \odo \id)\Sigma_{12}\ad(\Sigma(U)^{-1}_{23})\Sigma_{12}(\theta \odo y \odo x_{[1]}) \\
& = ({}_{h}\p \odo \id \odo \id)\ad(\Sigma(U)^{-1}_{13})(\beta \odo \id)(\alp(x)) 
\end{align*}
This implies that the equality $\alp(h \lhd_{\beta} x) = (x_{[0]} \lhd_{\beta} h_{(2)}) \odo S(h_{(1)})x_{[1]}h_{(3)}$ holds for all $h \in H$ and $x \in X$ if and only if $(\id_{\du{H}} \odo \alp)\beta = \ad(\Sigma(U)^{-1}_{13})(\beta \odo \id_{H})\alp$.
\fin

\subsubsection*{Braided commutativity}\label{sec:rl-bc-mod}

Given a right-left $(H,\du{H})$-Yetter--Drinfeld algebra $(X,\alp,\beta)$, we known that $(X,\lhd_{\beta},\alp)$ is a (right-right) Yetter--Drinfeld algebra over $H$. Observe that the braided commutativity condition~\eqref{eq:rr-bc-1} for  $(X,\lhd_{\beta},\alp)$ can be written as the condition $m_{X}\tau_{X} = m_{X}$, where $\tau_{X} = (\p_{-} \odo \id_{X \odo X})\Sigma_{34}\Sigma_{23}\Sigma_{12}(\alp \odo \beta) : X \odo X \to X \odo X$. Similarly, the condition~\eqref{eq:rr-bc-2} can be written as the condition $m_{X}\rho_{X} = m_{X}$, where $\rho_{X} = (\p \odo \id_{X \odo X})\Sigma_{12}\Sigma_{24}(\beta \odo \alp) : X \odo X \to X \odo X$. One can verified that $\tau_{X}\rho_{X} = \id_{X \odo X} = \rho_{X}\tau_{X}$.

\begin{rema}
We use the notation $\p: H \odo \du{H} \to \ku$ for the canonical pairing $\p(h \odo \theta)= \theta(h)$, and $\p_{-}: H \odo \du{H} \to \ku$ for the linear map $\p_{-}(h \odo \theta)= \theta(S^{-1}(h)) = \du{S}^{-1}(\theta)(h)$.
\end{rema}

Consider the algebra $X \odo \He(\du{H})$ with multiplication given by
\[
\m_{X \odo \He(\du{H})} = (\m_{X} \odo \m_{\He(\du{H})})(\id \odo \Sigma_{\He(\du{H}),X} \odo \id),
\]
here $\He(\du{H})$ is the Heisenberg algebra of $\du{H}$. Recall that in $\He(\du{H})$, we have
\[
\m_{\He(\du{H})}((\theta \# h) \odo (\theta' \# h')) = \theta(h_{(1)} \brhd \theta') \# h_{(2)}h' = \theta\theta'_{(1)} \# (h \blhd \theta'_{(2)})h'.
\]
Denote by $\hat\iota_{X,1} : X \odo \du{H} \to X \odo \He(\du{H})$, $x \odo \theta \mapsto x \odo \theta \# 1$ and by $\hat\iota_{X,2}: X \odo H \to X \odo \He(\du{H})$, $x \odo h \mapsto x \odo 1 \# h$, the canonical inclusions induced by the inclusions of $\du{H}$ and $H$ in the Heisenberg algebra $\He(\du{H})$. Note that the canonical flip pairing can be see as a linear map on $\He(\du{H})$, $\bar{\p} : \du{H} \# H \to \ku$, $\theta \odo h \mapsto \p(h,\theta)$.

\begin{defi}\label{df:rl-bc-mod}
We will say that a right-left $(H,\du{H})$-Yetter--Drinfeld algebra $(X,\alp,\beta)$ is {\em braided commutative} if
\begin{equation}\label{eq:rl-bc-mod}
\m_{X \odo \He(\du{H})}(\hat\iota_{X,2} \odo \hat\iota_{X,1})(\alp(x) \odo \Sigma\beta(y)) = \m_{X \odo \He(\du{H})}(\hat\iota_{X,2} \odo \hat\iota_{X,1})(\Sigma\beta(y) \odo \alp(x))
\end{equation}
for every $x,y \in X$.
\end{defi}

\begin{nota}
We will denote
\[
\aYD^{\,rl}_{bc}(H,\du{H}) = \bigl\{ (X,\alp,\beta) \in \aYD^{\,rl}(H,\du{H}) \;:\; (X,\alp,\beta) \text{ is braided commutative} \bigr\}
\]
\end{nota}

The braided commutativity condition~\eqref{eq:rl-bc-mod} for Yetter--Drinfeld algebras is equivalent to the braided commutativity condition in the standard characterization of Yetter--Drinfeld algebras as we prove in the next proposition. 

\begin{prop}\label{rl-bc-equiv}
The right-left $(H,\du{H})$-Yetter--Drinfeld algebra $(X,\alp,\beta)$ is braided commutative if and only if the (right-right) Yetter--Drinfeld algebra over $H$, $(X,\lhd_{\beta},\alp)$, is braided commutative.
\end{prop}
\pr
Given $x,y \in X$, if we use the notations $x_{[0]} \odo x_{[1]} := \alp(x)$ and $y^{[-1]} \odo y^{[0]} := \beta(y)$, we have
\small
\begin{align*}
\m_{X \odo \He(\du{H})}(\hat\iota_{X,2} \odo \hat\iota_{X,1})(\alp(x) \odo \Sigma\beta(y)) & = \m_{X \odo \He(\du{H})}(x_{[0]} \odo 1 \# x_{[-1]} \odo y_{[0]} \odo y^{[-1]} \# 1) \\
& = x_{[0]}y^{[0]} \odo (1 \# x_{[1]})(y^{[1]} \# 1) \\
& = x_{[0]}y^{[0]} \odo (x_{[1]_{(1)}} \brhd y^{[-1]}) \# x_{[1]_{(2)}} \\
& = \p(x_{[1]_{(1)}} \odo y^{[-1]}_{(2)})x_{[0]}y^{[0]} \odo y^{[-1]}_{(1)} \# x_{[1]_{(2)}} \\
& = (\p \odo \m_{X})(x_{[0]_{[1]}} \odo y^{[0]^{[-1]}} \odo x_{[0]_{[0]}} \odo y^{[0]^{[0]}}) \odo y^{[-1]} \# x_{[1]} \\
& = (\p \odo \m_{X})\Sigma_{12}\Sigma_{24}(\beta \odo \alp)(y^{[0]} \odo x_{[0]}) \odo y^{[-1]} \# x_{[1]} \\
& = (\m_{X}\rho_{X} \odo \id_{\He(H)})(y^{[0]} \odo x_{[0]} \odo y^{[-1]} \# x_{[1]})
\intertext{\normalsize and}
\m_{X \odo \He(\du{H})}(\hat\iota_{X,1} \odo \hat\iota_{X,2})(\Sigma\beta(y) \odo \alp(x)) & = \m_{X \odo \He(\du{H})}(y^{[0]} \odo y^{[-1]} \# 1 \odo x_{[0]} \odo 1 \# x_{[1]}) \\
& = (\m_{X} \odo \id_{\He(\du{H})})(y^{[0]} \odo x^{[0]} \odo y^{[-1]} \# x_{[1]}).
\end{align*}
\normalsize

\begin{itemize}
\item[($\vuelta$)] If $(X,\lhd_{\beta},\alp)$ is braided commutative, we have $\m_{X} = \m_{X}\rho_{X}$, then follows immediately from the above equalities that the condition~\eqref{eq:rl-bc-mod} is true for every $x,y \in X$, i.e. the tuple $(X,\alp,\beta)$ is braided commutative.

\item[($\ida$)] Let $x',y' \in X$. Take $y \odo x = \tau_{X}\Sigma(y' \odo x')$, thus we have $y' \odo x' = \Sigma\rho_{X}(y \odo x)$. Because $(X,\alp,\beta)$ satisfies the condition~\eqref{eq:rl-bc-mod} for $x,y \in X$, then the equalities
\begin{align*}
(\id_{X} \odo \p)\m_{X \odo \He(\du{H})}(\hat\iota_{X,2} \odo \hat\iota_{X,1})(\alp(x) \odo \Sigma\beta(y)) & = (\m_{X}\rho_{X} \odo \p)(y^{[0]} \odo x^{[0]} \odo y^{[-1]} \# x_{[1]}) \\
& = \m_{X}\rho_{X}\Sigma\rho_{X}(y \odo x) \\
& = \m_{X}\rho_{X}(y' \odo x'), \\
\\
(\id_{X} \odo \p)\m_{X \odo \He(\du{H})}(\hat\iota_{X,1} \odo \hat\iota_{X,2})(\Sigma\beta(y) \odo \alp(x)) & = (\m_{X} \odo \p)(y^{[0]} \odo x^{[0]} \odo y^{[-1]} \# x_{[1]}) \\
& = \m_{X}(\p(x_{[1]} \# y^{[-1]})y^{[0]} \odo x_{[0]}) \\
& = \m_{X}((y \lhd_{\beta} x_{[1]}) \odo x_{[0]}) \\
& = \m_{X}\Sigma\rho_{X}(y \odo x) \\
& = \m_{X}(y' \odo x')
\end{align*}
imply that $\m_{X}(y' \odo x') = \m_{X}\rho_{X}(y' \odo x')$. Thus $(X,\lhd_{\beta},\alp)$ is braided commutative.
\end{itemize}
\fin

We have equivalent conditions for the braided commutativity. For that consider the algebra $X \odo \He(H^{\op})$ with multiplication given by
\[
\m_{X \odo \He(H^{\op})} = (\m_{X} \odo \m_{\He(H^{\op})})(\id \odo \Sigma_{\He(H^{\op}),X} \odo \id),
\]
here $\He(H^{\op})$ is the Heisenberg algebra of $H^{\op}$. Recall that in $\He(H^{\op})$, we have
\[
\m_{\He(H^{\op})}((h^{\op} \# \theta) \odo (h'^{\op} \# \theta')) = ((\theta_{(2)} \brhd h')h)^{\op} \# \theta_{(1)}\theta' = (h'_{(1)}h)^{\op} \# (h'_{(2)} \brhd \theta)\theta'.
\]

Denote by $\iota^{\ops}_{X,1} : X \odo H^{\op} \to X \odo \He(H^{\op})$, $x \odo h^{\op} \mapsto x \odo h^{\op} \# 1$ and by $\iota^{\ops}_{X,2}: X \odo \du{H}^{\co} \to X \odo \He(H^{\op})$, $x \odo \theta \mapsto x \odo 1^{\op} \# \theta$, the canonical inclusions induced by the inclusions of $H^{\op}$ and $\du{H}^{\co}$ in the Heisenberg algebra $\He(H^{\op})$. Note that the canonical pairing can be see as a linear map on $\He(H^{\op})$, $\p : H^{\op} \# \du{H}^{\co} \to \ku$, $h^{\op} \odo \theta \mapsto \p(h,\theta)$.

\begin{prop}
Let $(X,\alp,\beta)$ be a right-left $(H,\du{H})$-Yetter--Drinfeld algebra. We have
\begin{equation}\label{eq:rl-bc-mod-2}
\m_{X \odo \He(H^{\op})}(\iota^{\ops}_{X,1} \odo \iota^{\ops}_{X,2})((\id \odo S^{\ops})\dot\alp(x) \odo \Sigma\beta(y)) = \m_{X \odo \He(H^{\op})}(\iota^{\ops}_{X,2} \odo \iota^{\ops}_{X,1})(\Sigma\beta(y) \odo (\id \odo S^{\ops})\dot\alp(x))
\end{equation}
for all $x,y \in X$, if and only if $(X,\lhd_{\beta},\alp)$ is braided commutative.
\end{prop}
\pr
Given $x,y \in X$, using the notations $x_{[0]} \odo x_{[1]} := \alp(x)$ and $y^{[-1]} \odo y^{[0]} := \beta(y)$, we have
\footnotesize
\begin{align*}
\m_{X \odo \He(H^{\op})}(\iota^{\ops}_{X,1} \odo \iota^{\ops}_{X,2})((\id \odo S^{\ops})\dot\alp(x) \odo \Sigma\beta(y)) & = \m_{X \odo \He(H^{\op})}(x_{[0]} \odo S^{-1}(x_{[1]})^{\op} \# 1 \odo y^{[0]} \odo 1^{\op} \# y^{[-1]}) \\
& = (\m_{X} \odo \id_{\He(H^{\op})})(x_{[0]} \odo y^{[0]} \odo S^{-1}(x_{[1]})^{\op} \# y^{[-1]}) \\
\intertext{\normalsize and}
\m_{X \odo \He(H^{\op})}(\iota^{\ops}_{X,2} \odo \iota^{\ops}_{X,1})(\Sigma\beta(y) \odo (\id \odo S^{\ops})\dot\alp(x)) & = \m_{X \odo \He(H^{\op})}(y^{[0]} \odo 1^{\op} \# y^{[-1]} \odo x_{[0]} \odo  S^{-1}(x_{[1]})^{\op} \# 1) \\
& = y^{[0]}x_{[0]} \odo (y^{[-1]}_{(2)} \brhd S^{-1}(x_{[1]}))^{\op} \# y^{[-1]}_{(1)} \\
& = \p(S^{-1}(x_{[1]_{(1)}}),y^{[-1]}_{(2)})y^{[0]}x_{[0]} \odo S^{-1}(x_{[1]_{(2)}})^{\op} \# y^{[-1]}_{(1)} \\
& = (\p_{-} \odo \m_{X})(x_{[0]_{[1]}} \odo y^{[0]^{[-1]}} \odo y^{[0]^{[0]}} \odo x_{[0]_{[0]}}) \odo S^{-1}(x_{[1]})^{\op} \# y^{[-1]} \\
& = (\p_{-} \odo \m_{X})\Sigma_{34}\Sigma_{23}\Sigma_{12}(\alp \odo \beta)(x_{[0]} \odo y^{[0]}) \odo S^{-1}(x_{[1]})^{\op} \# y^{[-1]} \\
& = (\m_{X}\tau_{X} \odo \id_{\He(H)})(x_{[0]} \odo y^{[0]} \odo S^{-1}(x_{[1]})^{\op} \# y^{[-1]}).
\end{align*}
\normalsize

\begin{itemize}
\item[($\vuelta$)] If $(X,\lhd_{\beta},\alp)$ is braided commutative, we have $\m_{X} = \m_{X}\tau_{X}$, then follows immediately from the above equalities that the condition~\eqref{eq:rl-bc-mod-2} is true for every $x,y \in X$, i.e. the tuple $(X,\alp,\beta)$ is braided commutative.

\item[($\ida$)] Let $x',y' \in X$. Take $x \odo y = \rho_{X}\Sigma(x' \odo y')$, thus we have $x' \odo y' = \Sigma\tau_{X}(x \odo y)$. Because $(X,\alp,\beta)$ satisfies the condition~\eqref{eq:rl-bc-mod-2} for $x,y \in X$, then the equalities
\small
\begin{align*}
(\id_{X} \odo \p)\m_{X \odo \He(H^{\op})}(\iota^{\ops}_{X,1} \odo \iota^{\ops}_{X,2})((\id \odo S^{\ops})\dot\alp(x) \odo \Sigma\beta(y)) & = (\m_{X} \odo \p)(x_{[0]} \odo y^{[0]} \odo S^{-1}(x_{[1]})^{\op} \# y^{[-1]}) \\
& = \m_{X}(x_{[0]} \odo \p(S^{-1}(x_{[1]}) \odo y^{[-1]})y^{[0]}) \\
& = \m_{X}\Sigma\tau_{X}(x \odo y) \\
& = \m_{X}(x' \odo y'),
& \\
(\id_{X} \odo \p)\m_{X \odo \He(H^{\op})}(\iota^{\ops}_{X,2} \odo \iota^{\ops}_{X,1})(\Sigma\beta(y) \odo (\id \odo S^{\ops})\dot\alp(x)) & = (\m_{X}\tau_{X} \odo \p)(x_{[0]} \odo y^{[0]} \odo S^{-1}(x_{[1]})^{\op} \# y^{[-1]}) \\
& = \m_{X}\tau_{X}\Sigma\tau_{X}(x \odo y) \\
& = \m_{X}\tau_{X}(x' \odo y'),
\end{align*}
\normalsize
imply that $\m_{X}(x' \odo y') = \m_{X}\tau_{X}(x' \odo y')$. Thus $(X,\lhd_{\beta},\alp)$ is braided commutative.
\end{itemize}
\fin

\begin{prop}\label{prop:dual-rl-bc-mod}
A right-left $(H,\du{H})$-Yetter--Drinfeld algebra $(X,\alp,\beta)$ is braided commutative if and only if
\begin{multline*}
\m_{X^{\op} \odo \He(H)}(\iota_{X^{\op},1} \odo \iota_{X^{\op},2})(({}^{\op} \odo S^{-1})\alp(x) \odo \Sigma(\dS \odo {}^{\op})(\beta(y)) \\
= \m_{X^{\op} \odo \He(H)}(\iota_{X^{\op},2} \odo \iota_{X^{\op},1})(\Sigma(\dS \odo {}^{\op})(\beta(y) \odo ({}^{\op} \odo S)\alp(x))
\end{multline*}
for every $x,y \in X$.
\end{prop}
\pr
Consider the anti-isomorphism $\mathcal{L}_{1}: \He(\du{H}) \to \He(H)$, see Proposition~\ref{prop:anti_heisen}. Observe that
\[
({}^{\op} \odo \mathcal{L}_{1})\hat\iota_{X,1}\Sigma\beta = \iota_{X^{\op},2}({}^{\op} \odo \du{S})\Sigma\beta \quad \text{ and } \quad ({}^{\op} \odo \mathcal{L}_{1})\hat\iota_{X,2}\alp = \iota_{X^{\op},1}({}^{\op} \odo S^{-1})\alp.
\]
Because $({}^{\op} \odo \mathcal{L}_{1})$ is an anti-isomorphism between $X \odo \He(\du{H})$ and $X^{\op} \odo \He(H)$, it hold
\footnotesize
\begin{align*}
({}^{\op} \odo \mathcal{L}_{1})\m_{X \odo \He(\du{H})}(\hat\iota_{X,2} \odo \hat\iota_{X,1})(\alp(x) \odo \Sigma\beta(y)) = \m_{X^{\op} \odo \He(H)}(\iota_{X^{\op},2} \odo \iota_{X^{\op},1})(({}^{\op} \odo \du{S})\Sigma\beta(y) \odo ({}^{\op} \odo S^{-1})\alp(x))
\end{align*}
\normalsize
and
\footnotesize
\begin{align*}
({}^{\op} \odo \mathcal{L}_{1})\m_{X \odo \He(\du{H})}(\hat\iota_{X,1} \odo \hat\iota_{X,2})(\Sigma\beta(y) \odo \alp(x))  = \m_{X^{\op} \odo \He(H)}(\iota_{X^{\op},1} \odo \iota_{X^{\op},2})(({}^{\op} \odo S^{-1})\alp(x) \odo ({}^{\op} \odo \du{S})\Sigma\beta(y)).
\end{align*}
\normalsize
for all $x,y \in X$. The result follows from those equalities.
\fin

\begin{rema}
It is obvious that
\begin{enumerate}[label=\textup{(\roman*)}]
\item $\alp: X \to X \odo H$ is a right coaction of $H$ on $X$ if and only if $(\Sigma\alp)^{\co}: X^{\op} \to H \odo X^{\op}$ is a left coaction of $H$ on $X^{\op}$.
\item $\beta: X \to \du{H} \odo X$ is a left coaction of $\du{H}$ on $X$ if and only if  $\Sigma\beta^{\co}: X^{\op} \to  X^{\op} \odo \du{H}$ is a right coaction of $\du{H}$ on $X^{\op}$.
\end{enumerate}
\end{rema}

With the last remark, we have the following proposition.

\begin{prop}
We have
\[
(X,\alp,\beta) \in \aYD^{\,rl}(H,\du{H}) \qquad\text{ if and only if }\qquad (X^{\op},\Sigma\beta^{\co},(\Sigma\alp)^{\co})  \in \aYD^{\,rl}(\du{H},H).
\]
Moreover, $(X,\alp,\beta)$ is braided commutative if and only if $(X^{\op},\Sigma\beta^{\co},(\Sigma\alp)^{\co})$ is braided commutative.
\end{prop}
\pr
First, note that
\begin{align*}
(\id_{H} \odo \Sigma\beta^{\co})\circ(\Sigma\alp)^{\co} & = \Sigma_{23}\circ(S^{-1} \odo \dS \odo {}^{\op})\circ(\id \odo \beta)\circ\Sigma\alp\circ{}^{\op} \\
& = \Sigma_{23}\circ(S^{-1} \odo \dS \odo {}^{\op})\circ\Sigma_{12}\circ\Sigma_{23}\circ(\beta \odo \id)\circ\alp\circ{}^{\op} \\
& = \Sigma_{23}\circ\Sigma_{12}\circ\Sigma_{23}\circ(\dS \odo {}^{\op} \odo S^{-1})\circ(\beta \odo \id)\circ\alp\circ{}^{\op} 
\end{align*}
and
\begin{align*}
((\Sigma\alp)^{\co} \odo \id_{\du{H}})\circ\Sigma\beta^{\co} & = (S^{-1} \odo {}^{\op} \odo \id)\circ\Sigma_{12}\circ(\alp \odo \id)\circ({}^{\op} \odo \id)\circ\Sigma\circ(\dS \odo {}^{\op})\circ\beta\circ{}^{\op} \\
& = (S^{-1} \odo {}^{\op} \odo \id)\circ\Sigma_{12}\circ\Sigma_{23}\circ\Sigma_{12}\circ(\dS \odo \id \odo \id)\circ(\id \odo \alp)\circ\beta\circ{}^{\op} \\
& = \Sigma_{12}\circ\Sigma_{23}\circ\Sigma_{12}\circ(\dS \odo {}^{\op} \odo S^{-1})\circ(\id \odo \alp)\circ\beta\circ{}^{\op} \\
& = \Sigma_{23}\circ\Sigma_{12}\circ\Sigma_{23}\circ(\dS \odo {}^{\op} \odo S^{-1})\circ(\id \odo \alp)\circ\beta\circ{}^{\op}.
\end{align*}
Because, $(\dS \odo {}^{\op} \odo S^{-1})\ad(\Sigma(U)^{-1}_{13}) = \ad(\Sigma(U)_{13})(\dS \odo {}^{\op} \odo S^{-1})$, then the equality
\[
(\id_{\du{H}} \odo \alp)\beta = \ad(\Sigma(U)^{-1}_{13})(\beta \odo \id_{H})\alp.
\]
is equivalent to the equality $(\id_{H} \odo \Sigma\beta^{\co})(\Sigma\alp)^{\co} = \ad(U^{-1}_{13})((\Sigma\alp)^{\co} \odo \id_{\du{H}})\Sigma\beta^{\co}$. This implies, $(X,\alp,\beta) \in \aYD^{\,rl}(H,\du{H})$ if and only if $(X^{\op},(\Sigma\beta)^{\co},\Sigma\alp^{\co})  \in \aYD^{\,rl}(\du{H},H)$. Follows from Definition~\ref{df:rl-bc-mod} and Proposition~\ref{prop:dual-rl-bc-mod} the equivalence about the braided commutativity conditions.
\fin

\subsection{(left-left) Yetter--Drinfeld algebras}

Let $\alp : X \to H^{\op} \odo X$ a left coaction of $H^{\op}$ on $X$ and $\beta : X \to \du{H} \odo X$ left coaction of $\du{H}$ on $X$.

\begin{defi}\label{def:new-ll-YD}
We say that the tuple $(X,\alp,\beta)$ is a {\em left-left $(H,\du{H})$-Yetter--Drinfeld algebra} if
\begin{equation}
(\id_{\du{H}} \odo \alp)\beta = (\Sigma \odo \id_{X})(\ad({}^{\ops}U \odo \id_{X})(\id_{H^{\op}} \odo \beta)\alp.
\end{equation}
In other words, if the diagram
\[
\xymatrix@C=6pc@R=2pc{X\ar[r]^-{\beta}\ar[d]_-{\alp} & \du{H} \odo X\ar[r]^-{id_{\du{H}} \odo \alp} & \du{H} \odo H^{\op} \odo X \\ H^{\op} \odo X \ar[r]_-{id_{H^{\op}} \odo \beta} & H^{\op} \odo \du{H} \odo X\ar[r]_-{\ad({}^{\ops}U) \odo \id_{X}} & H^{\op} \odo \du{H} \odo X\ar[u]_-{\Sigma \odo \id_{X}}}
\]
is commutative.
\end{defi}

\begin{nota}
Sometimes we will denoted $(X,\alp,\beta) \in \aYD^{\,ll}(H,\du{H})$ for say that $(X,\alp,\beta)$ is a left-left $(H,\du{H})$-Yetter--Drinfeld algebra.
\end{nota}

\begin{prop}\label{prop:ll-YD-equiv}
$(X,\alp,\beta)$ is a left-left $(H,\du{H})$-Yetter--Drinfeld algebra if and only if $(X,\lhd_{\beta},\alp)$ is a (right-left) Yetter--Drinfeld algebra over $H$, where
\[
\begin{array}{lccc}
\lhd_{\beta} : & X \odo H & \to & X \\ & x \odo h & \mapsto & ({}_{h}\p \odo \id_{X})(\beta(x))
\end{array}
\]
is the right action of the Hopf algebra $(H,\com)$ on $X$ induced by the coaction $\beta$.
\end{prop}
\pr
Fix $h \in H$ and $x \in X$. It is obvious that $\alp(x \lhd_{\beta} h) = ({}_{h}\p \odo \id_{H^{\op}} \odo \id_{X})(\id_{\du{H}} \odo \alp)(\beta(x))$. Now, if we use the notation $\theta \odo y = \beta(x_{[0]})$ then
\begin{align*}
S^{-1}(h_{(3)})x_{[-1]}h_{(1)} \odo (x_{[0]} \lhd_{\beta} h_{(2)}) & = S^{-1}(h_{(3)})x_{[-1]}h_{(1)} \odo ({}_{h_{(2)}}\p \odo \id)(\theta \odo y) \\
& = S^{-1}(h_{(3)})x_{[-1]}h_{(1)}\p(h_{(2)},\theta) \odo y \\
& = S^{-1}(h_{(2)})x_{[-1]}(\theta \brhd h_{(1)}) \odo y \\
& = (\id_{H} \odo {}_{h}\p)(T_{U}(x_{[-1]} \odo \theta)) \odo y \\
& = ({}_{h}\p \odo \id_{H} \odo \id_{X})(\Sigma \odo \id_{X})(T_{U} \odo \id_{X})(x_{[-1]} \odo \theta \odo y) \\
& = ({}_{h}\p \odo \id_{H} \odo \id_{X})(\Sigma T_{U} ({}^{\op} \odo \id) \odo \id_{X})(\id_{H^{\op}} \odo \beta)(\alp(x)).
\end{align*}
This implies that the equality $\alp(x \lhd_{\beta} h) = (S^{-1}(h_{(3)})x_{[-1]}h_{(1)})^{\op} \odo (x_{[0]} \lhd_{\beta} h_{(2)})$ holds for all $h \in H$ and $x \in X$ if and only if
\[
(\id_{\du{H}} \odo \alp)\beta = (\Sigma({}^{\op} \odo \id_{\du{H}}) T_{U}({}^{\op} \odo \id_{\du{H}}) \odo \id_{X})(\id_{H^{\op}} \odo \beta)\alp.
\]
Finally, the proposition holds, because $({}^{\op} \odo \id_{\du{H}}) T_{U}({}^{\op} \odo \id_{\du{H}}) = \ad({}^{\ops}U)$.
\fin

\subsubsection*{Braided commutativity}\label{sec:ll-mod-bc}

Given a left-left $(H,\du{H})$-Yetter--Drinfeld algebra $(X,\alp,\beta)$, we known that $(X,\lhd_{\beta},\alp)$ is a (right-left) Yetter--Drinfeld algebra over $H$. Observe that the braided commutativity condition~\eqref{eq:rl-bc-1} for $(X,\lhd_{\beta},\alp)$, can be written as $\m_{X}\tau_{X} = \m_{X}$, where $\tau_{X} = (\p \odo \id \odo \id)\Sigma_{34}\Sigma_{23}(\dot\alp \odo \beta) : X \odo X \to X \odo X$. Similarly, the condition~\eqref{eq:rl-bc-2} can be written as the condition $\m_{X}\rho_{X} = \m_{X}$, where $\rho_{X} = (\p_{+} \odo \id \odo \id)\Sigma_{12}\Sigma_{23}\Sigma_{24}(\beta \odo \dot\alp) : X \odo X \to X \odo X$. Observe that $\tau_{X}\rho_{X} = \id_{X \odo X} = \rho_{X}\tau_{X}$.

\begin{rema}
We use the notation $\p: H \odo \du{H} \to \ku$ for the canonical pairing $\p(h \odo \theta)= \theta(h)$, and $\p_{+}: H \odo \du{H} \to \ku$ for the linear map $\p_{+}(h \odo \theta) = \theta(S(h)) = \du{S}(\theta)(h)$.
\end{rema}

Consider the algebra $\He(H) \odo X$ with multiplication given by
\[
\m_{\He(H) \odo X} = (\m_{\He(H)} \odo \m_{X})(\id \odo \Sigma_{X,\He(H)} \odo \id),
\]
here $\He(H) := H \# \du{H}$ is the Heisenberg algebra of $H$. Recall that in $\He(H)$, we have
\[
\m_{\He(H)}((h \# \theta) \odo (h' \# \theta')) = h(\theta_{(1)} \brhd h') \# \theta_{(2)}\theta'.
\]
Denote by $\iota_{1} : H \odo X^{\op} \to \He(H) \odo X^{\op}$, $h \odo x^{\op} \mapsto h \# 1 \odo x^{\op}$ and by $\iota_{2}: \du{H} \odo X^{\op} \to \He(H) \odo X^{\op}$, $\theta \odo x^{\op} \mapsto 1 \# \theta \odo x^{\op}$, the canonical inclusions induced by the inclusions of $H$ and $\du{H}$ in the Heisenberg algebra $\He(H)$. Note that the canonical pairing can be see as a linear map on $\He(H)$, $\p : H \# \du{H} \to \ku$, $h \# \theta \mapsto \theta(h)$.

\begin{defi}\label{df:ll-bc-mod}
We will say that a left-left $(H,\du{H})$-Yetter--Drinfeld algebra $(X,\alp,\beta)$ is {\em braided commutative} if
\small
\begin{equation}\label{eq:ll-bc-mod}
\m_{\He(H) \odo X^{^{\op}}}(\iota_{1,X^{\op}} \odo \iota_{2,X^{\op}})((\id \odo {}^{\op})\dot\alp(x) \odo (\dS \odo {}^{\op})\beta(y)) = \m_{\He(H) \odo X^{\op}}(\iota_{2,X^{\op}} \odo \iota_{1,X^{\op}})((\du{S} \odo {}^{\op})\beta(y) \odo (\id \odo {}^{\op})\dot\alp(x))
\end{equation}
\normalsize
for every $x,y \in X$.
\end{defi}

\begin{nota}
We will denote
\[
\aYD^{\,ll}_{bc}(H,\du{H}) = \bigl\{ (X,\alp,\beta) \in \aYD^{\,ll}(H,\du{H}) \;:\; (X,\alp,\beta) \text{ is braided commutative} \bigr\}
\]
\end{nota}

The braided commutativity condition~\eqref{eq:ll-bc-mod} for Yetter--Drinfeld algebras is equivalent to the braided commutativity condition in the standard characterization of Yetter--Drinfeld algebras as we prove in the next proposition. 

\begin{prop}\label{prop:ll-bc-equiv}
The left-left $(H,\du{H})$-Yetter--Drinfeld algebra $(X,\alp,\beta)$ is braided commutative if and only if the (right-left) Yetter--Drinfeld algebra over $H$, $(X,\lhd_{\beta},\alp)$, is braided commutative.
\end{prop}
\pr
Given $x,y \in X$, using the notations $x^{\op}_{[-1]} \odo x_{[0]} := \alp(x)$ and $y^{[-1]} \odo y^{[0]} := \beta(y)$, we have
\footnotesize
\begin{align*}
\m_{\He(H) \odo X^{^{\op}}}(\iota_{1,X^{\op}} \odo \iota_{2,X^{\op}})((\id \odo {}^{\op})\dot\alp(x) \odo (\dS \odo {}^{\op})\beta(y)) & = \m_{\He(H) \odo X^{\op}}(x_{[-1]}\# 1 \odo x^{\op}_{[0]} \odo 1 \# \dS(y^{[-1]}) \odo y^{[0]^{\op}}) \\
& = (\id_{\He(H)} \odo {}^{\op}\m_{X})(x_{[-1]} \# \dS(y^{[-1]}) \odo  y^{[0]} \odo x_{[0]})
\intertext{and}
\m_{\He(H) \odo X^{\op}}(\iota_{2,X^{\op}} \odo \iota_{1,X^{\op}})((\dS \odo {}^{\op})\beta(y) \odo (\id \odo {}^{\op})\dot\alp(x)) & = \m_{\He(H) \odo X^{\op}}(1 \# \dS(y^{[-1]}) \odo y^{[0]^{\op}} \odo x_{[-1]} \# 1 \odo x^{\op}_{[0]}) \\
& = (1 \# \dS(y^{[-1]}))(x_{[-1]} \# 1) \odo (x_{[0]}y^{[0]})^{\op} \\
& = (\dS(y^{[-1]}_{(2)}) \brhd x_{[-1]}) \# \dS(y^{[-1]}_{(1)}) \odo (x_{[0]}y^{[0]})^{\op} \\
& = x_{[-1]_{(1)}} \# \dS(y^{[-1]}_{(1)}) \odo \p_{+}(x_{[-1]_{(2)}} \odo y^{[-1]}_{(2)}) (x_{[0]}y^{[0]})^{\op} \\
& = x_{[-1]} \# \dS(y^{[-1]}) \odo \p_{+}(x_{[0]_{[-1]}} \odo y^{[0]^{[-1]}}) (x_{[0]_{[0]}}y^{[0]^{[0]}})^{\op} \\
& = x_{[-1]} \# \dS(y^{[-1]}) \odo (\p_{+} \odo {}^{\op}\m_{X})(x_{[0]_{[-1]}} \odo y^{[0]^{[-1]}} \odo x_{[0]_{[0]}} \odo y^{[0]^{[0]}}) \\
& = x_{[-1]} \# \dS(y^{[-1]}) \odo (\p_{+} \odo {}^{\op}\m_{X})\Sigma_{12}\Sigma_{23}\Sigma_{24}(\beta \odo \dot\alp)(y^{[0]} \odo x_{[0]}) \\
& = x_{[-1]} \# \dS(y^{[-1]}) \odo {}^{\op}\m_{X}\rho_{X}(y^{[0]} \odo x_{[0]}) \\
& = (\id_{\He(H)} \odo {}^{\op}\m_{X}\rho_{X})(x_{[-1]} \# \dS(y^{[-1]}) \odo y^{[0]} \odo x_{[0]}).
\end{align*}
\normalsize

\begin{itemize}
\item[($\vuelta$)] If $(X,\lhd_{\beta},\alp)$ is braided commutative, we have $\m_{X} = \m_{X}\rho_{X}$, follows immediately from the above equalities that the condition~\eqref{eq:ll-bc-mod} is true for every $x,y \in X$, i.e. the tuple $(X,\alp,\beta)$ is braided commutative.

\item[($\ida$)] Let $x',y' \in X$. Take $y \odo x = \tau_{X}\Sigma(y' \odo x')$, thus we have $y' \odo x' = \Sigma\rho_{X}(y \odo x)$. If $(X,\alp,\beta)$ satisfies the condition~\eqref{eq:ll-bc-mod}, then the equalities
\footnotesize
\begin{align*}
(\p \odo \id_{X^{\op}})\m_{\He(H) \odo X^{\op}}(\iota_{1,X^{\op}} \odo \iota_{2,X^{\op}})(\id \odo {}^{\op})\dot\alp(x) \odo (\dS \odo {}^{\op})\beta(y)) 
& = {}^{\op}\m_{X}(\p(x_{[-1]} \# \dS(y^{[-1]}))y^{[0]} \odo x_{[0]}) \\
& = {}^{\op}\m_{X}(\p(S(x_{[-1]}) \# y^{[-1]})y^{[0]} \odo x_{[0]}) \\
& = {}^{\op}\m_{X}((y \lhd_{\beta} S(x_{[-1]})) \odo x_{[0]}) \\
& = {}^{\op}\m_{X}\Sigma\rho_{X}(y \odo x) \\
& = {}^{\op}\m_{X}(y' \odo x'), \\
(\p \odo \id_{X^{\op}})\m_{\He(H) \odo X^{\op}}(\iota_{2,X^{\op}} \odo \iota_{1,X^{\op}})((\dS \odo {}^{\op})\beta(y) \odo (\id \odo {}^{\op})\dot\alp(x)) 
& = {}^{\op}\m_{X}\rho_{X}\Sigma\rho_{X}(y \odo x) \\
& = {}^{\op}\m_{X}\rho_{X}(y' \odo x'),
\end{align*}
\normalsize
implies $\m_{X}(y' \odo x') = \m_{X}\rho_{X}(y' \odo x')$. Thus $(X,\lhd_{\beta},\alp)$ is braided commutative.
\end{itemize}
\fin

We have an equivalent condition for the braided commutativity.

\begin{prop}\label{prop:dual-ll-bc-mod}
A left-left $(H,\du{H})$-Yetter--Drinfeld algebra $(X,\alp,\beta)$ is braided commutative if and only if
\[
\m_{\He(\du{H}) \odo X}(\du{\iota}_{2,X} \odo \du{\iota}_{1,X})((S \odo \id)\dot\alp(x) \odo \beta(y)) = \m_{\He(\du{H}) \odo X}(\du{\iota}_{1,X} \odo \du{\iota}_{2,X})(\beta(y) \odo (S \odo \id)\dot\alp(x))
\]
for every $x,y \in X$.
\end{prop}
\pr
Consider the anti-isomorphism $\mathcal{L}_{1}: \He(\du{H}) \to \He(H)$, see Proposition~\ref{prop:anti_heisen}. Observe that
\[
(\mathcal{L}_{1} \odo {}^{\op})\hat\iota_{1,X}\beta = \iota_{2,X^{\op}}(\dS \odo {}^{\op})\beta \quad \text{ and } \quad (\mathcal{L}_{1} \odo {}^{\op})\hat\iota_{2,X}(S \odo \id)\dot\alp = \iota_{1,X^{\op}}(\id \odo {}^{\op})\dot\alp.
\]
Because $(\mathcal{L}_{1} \odo {}^{\op})$ is an anti-isomorphism between $\He(\du{H}) \odo X$ and $\He(H) \odo X^{\op}$, it hold
\footnotesize
\begin{align*}
(\mathcal{L}_{1} \odo {}^{\op})\m_{\He(\du{H})^\odo X}(\du{\iota}_{1,X} \odo \du{\iota}_{2,X})(\beta(y) \odo (S \odo \id)\dot\alp(x)) = \m_{\He(H) \odo X^{\op}}(\iota_{1,X^{\op}} \odo \iota_{2,X^{\op}})((\id \odo {}^{\op})\dot\alp(x) \odo (\du{S} \odo {}^{\op})\beta(y))
\end{align*}
\normalsize
and
\footnotesize
\begin{align*}
(\mathcal{L}_{1} \odo {}^{\op})\m_{\He(\du{H}) \odo X}(\du{\iota}_{1,X} \odo \du{\iota}_{2,X})(\beta(y) \odo (S \odo \id)\dot\alp(x)) = \m_{\He(H) \odo X^{\op}}(\iota_{2,X^{\op}} \odo \iota_{1,X^{\op}})((\du{S} \odo {}^{\op})\beta(y) \odo (\id \odo {}^{\op})\dot\alp(x))
\end{align*}
\normalsize
for every $x, y \in X$. The proposition follows from this equalities and the fact that $(\mathcal{L}_{1} \odo {}^{\op})$ is an anti-isomorphism.
\fin

\begin{rema}
It is obvious that
\begin{enumerate}[label=\textup{(\roman*)}]
\item $\alp: X \to H^{\op} \odo X$ is a left coaction of $H^{\op}$ on $X$ if and only if $\alp^{\op}=(S \odo \id)\alp^{\co}: X^{\op} \to H \odo X^{\op}$ is a left coaction of $H$ on $X^{\op}$.
\item $\beta: X \to \du{H} \odo X$ is a left coaction of $\du{H}$ on $X$ if and only if $\beta^{\op}=(\dS^{-1} \odo \id)\beta^{\co}: X^{\op} \to \du{H}^{\op} \odo X^{\op}$ is a left coaction of $\du{H}^{\op}$ on $X^{\op}$.
\end{enumerate}
\end{rema}

With the last remark, we have the following proposition.

\begin{prop}\label{prop:chau}
We have
\[
(X,\alp,\beta) \in \aYD^{\,ll}(H,\du{H}) \qquad\text{ if and only if }\qquad (X^{\op},(\dS^{-1} \odo \id)\beta^{\co},(S \odo \id)\alp^{\co}) \in \aYD^{\,ll}(\du{H},H).
\]
Moreover, $(X,\alp,\beta)$ is braided commutative if and only if $(X^{\op},(\dS^{-1} \odo \id)\beta^{\co},(S \odo \id)\alp^{\co})$ is braided commutative.
\end{prop}
\pr
First, note that
\begin{align*}
(\id_{H} \odo (\dS^{-1} \odo \id)\circ\beta^{\co})\circ(S \odo \id)\circ\alp^{\co} & = ({}^{\op} \odo {}^{\op} \odo {}^{\op})\circ(\id \odo \beta)\circ\alp\circ{}^{\op} 
\end{align*}
and
\begin{align*}
(\id_{\du{H}^{\op}} \odo (S \odo \id)\circ\alp^{\co})\circ(\dS^{-1} \odo \id)\circ\beta^{\co} & = ({}^{\op} \odo {}^{\op} \odo {}^{\op})\circ(\id \odo \alp)\circ\beta\circ{}^{\op}.
\end{align*}
Because, $({}^{\op}\odo {}^{\op} \odo {}^{\op})(\ad({}^{\ops}U)^{-1} \odo \id)(\Sigma \odo \id) = (\Sigma \odo \id)(\ad({}^{\ops}\Sigma(U)) \odo \id)({}^{\op}\odo {}^{\op} \odo {}^{\op})$, then the equality
\[
(\id_{\du{H}} \odo \alp)\beta = (\Sigma \odo \id)(\ad({}^{\ops}U) \odo \id)(\id_{H^{\ops}} \odo \beta)\alp
\]
is equivalent to the equality
\[
(\id_{H} \odo (\dS^{-1} \odo \id)\beta^{\co})(S \odo \id)\alp^{\co} = (\Sigma \odo \id)(\ad({}^{\ops}\Sigma(U)) \odo \id)(\id_{\du{H}^{\op}} \odo (S \odo \id)\alp^{\co})(\dS^{-1} \odo \id)\beta^{\co}.
\]
This implies, $(X,\alp,\beta) \in \aYD^{\,ll}(H,\du{H})$ if and only if $(X^{\op},(\dS^{-1} \odo \id)\beta^{\co},(S \odo \id)\alp^{\co})\in \aYD^{\,ll}(\du{H},H)$. Follows from Definition~\ref{df:ll-bc-mod} and Proposition~\ref{prop:dual-ll-bc-mod} the equivalence about the braided commutativity conditions.
\fin

\section{Yetter--Drinfeld algebras, Drinfeld doubles and Drinfeld codoubles}\label{sec:yd_dd}

\subsection{Majid's version}

Consider the Drinfeld double characterization as defined by Majid in \cite{M95}, $\D_{M}(H) = \du{H}^{\op} \bicroi H$. The algebra and coalgebra structure on $\D_{M}(H)$ are given by
\[
(\ome^{\op} \bicroi h)(\theta^{\op} \bicroi g) = ((h_{(2)} \brhd \theta_{(2)})\ome)^{\op} \bicroi (h_{(1)} \blhd \dS(\theta_{(1)}))g
\]
and
\[
\com_{\D_{M}(H)}(\ome^{\op} \bicroi h) = \ome^{\op}_{(1)} \bicroi h_{(1)} \odo \ome^{\op}_{(2)} \bicroi h_{(2)},
\]
respectively. Its dual Hopf algebra is given by $\T_{M}(H) = H^{\co} \odo \du{H}$ with algebra and coalgebra structure given by
\[
(h \odo \ome)(g \odo \theta) = hg \odo \ome\theta
\]
and
\[
\com_{\T_{M}(H)}(h \odo \ome) = (\id_{H} \odo \Sigma \odo \id_{\du{H}})(\id_{H} \odo \ad(U^{-1}) \odo \id_{\du{H}})(h_{(2)} \odo h_{(1)} \odo \ome_{(1)} \odo \ome_{(2)}),
\]
respectively. It is known in the standard characterization that give a left-left (resp. right-right) Yetter--Drinfeld algebra over a Hopf algebra $H$ is equivalent to give an algebra endowed with an left (resp. right) action of $\D_{M}(H)$.

\begin{theo}[\cite{M95}]\label{th:majid_duality}
Let $H$ be a finite-dimensional Hopf algebra. There is an equivalence between the category of (left-left) Yetter--Drinfeld algebras over $H$ and the category of left $\D_{M}(H)$-module algebras. More explicitly, given $(X,\rhd,\alp) \in {}_{H}^{H}\aYD$, we have a structure of left $\D_{M}(H)$-module algebra on $X$ given by
\[
(\ome^{\op} \bicroi h) \cdot x = (\p(\cdot,\ome) \odo \id)(\alp(h \rhd x)) = (h \rhd x) \lhd_{\alp} \ome,
\]
for all $h \in H$, $\ome \in \du{H}$ and $x \in X$. Reciprocally, given $(X,\cdot) \in {}_{\D_{M}(H)}\alg$, we have a left action and a left coaction of $H$ on $X$ given by
\[
h \rhd x := (1 \bicroi h) \cdot x \quad \text{ and } \quad \alp(x) = U_{1} \odo (U^{\op}_{2} \bicroi 1) \cdot x,
\]
for every $h \in H$ and $x \in X$, respectively (Here $U_{1} \odo U_{2} := U \in H \odo \du{H}$ is the canonical element associate to the Hopf algebra $H$). We have that $(X,\rhd,\alp)$ is a (left-left) Yetter--Drinfeld algebra over $H$. The equivalence of categories is given by the functor
\[
\begin{array}{rccc}
\mathfrak{E}_{M} : & {}_{H}^{H}\aYD & \leadsto & {}_{\D_{M}(H)}\alg \\
& (X,\rhd,\alp) & \mapsto & (X,\cdot) \\
& f: (X,\rhd,\alp) \to (X',\rhd',\alp') & \mapsto & f: (X,\cdotp) \to (X,\cdot')
\end{array}
\]
\end{theo}

\begin{theo}[\cite{M95}]\label{th:majid_duality_right}
Let $H$ be a finite-dimensional Hopf algebra. There is an equivalence between the category of (right-right) Yetter--Drinfeld algebras over $H$ and the category of right $\D_{M}(H)$-module algebras. More explicitly, given $(X,\lhd,\alp) \in \aYD^{H}_{H}$, we have a structure of right $\D_{M}(H)$-module algebra on $X$ given by
\[
x \cdot (\ome^{\op} \bicroi h) = ((\id \odo \p(\cdot,\ome))(\alp(x))) \lhd h  = (\ome \rhd_{\alp} x) \lhd h,
\]
for all $h \in H$, $\ome \in \du{H}$ and $x \in X$. Reciprocally, given $(X,\cdot) \in \alg_{\D_{M}(H)}$, we have a right action and a right coaction of $H$ on $X$ given by
\[
x \lhd h := x \cdot (1 \bicroi h) \quad \text{ and } \quad \alp(x) = x \cdot (U^{\op}_{2} \bicroi 1) \odo U_{1},
\]
for every $h \in H$ and $x \in X$, respectively (Here $U_{1} \odo U_{2} := U \in H \odo \du{H}$ is the canonical element associate to the Hopf algebra $H$). We have that $(X,\lhd,\alp)$ is a (right-right) Yetter--Drinfeld algebra over $H$. The equivalence of categories is given by the functor
\[
\begin{array}{rccc}
\mathfrak{E}_{M} : & \aYD^{H}_{H} & \leadsto & \alg_{\D_{M}(H)} \\
& (X,\lhd,\alp) & \mapsto & (X,\cdot) \\
& f: (X,\lhd,\alp) \to (X',\lhd',\alp') & \mapsto & f: (X,\cdotp) \to (X,\cdot')
\end{array}
\]
\end{theo}

The next propositions show that using the ``only coaction'' characterization of Yetter--Drinfeld algebras, the same kind of equivalence also holds but instead of use actions of the Drinfeld double we use coactions of its dual, the Drinfeld codouble.


\begin{prop}\label{prop:right_equivalence_majid}
The functor
\[
\begin{array}{rccc}
\mathfrak{E}'_{M} : & \aYD^{\,lr}(H,\du{H}) & \leadsto & \alg^{\T_{M}(H)} \\
& (X,\alp,\beta) & \mapsto & (X,\gamma = (\Sigma\alp \odo \id)\beta) \\
& f: (X,\alp,\beta) \to (X',\alp',\beta') & \mapsto & f: (X,\gamma) \to (X,\gamma')
\end{array}
\]
is an equivalence of categories. Moreover, we have the following equivalence of categories:
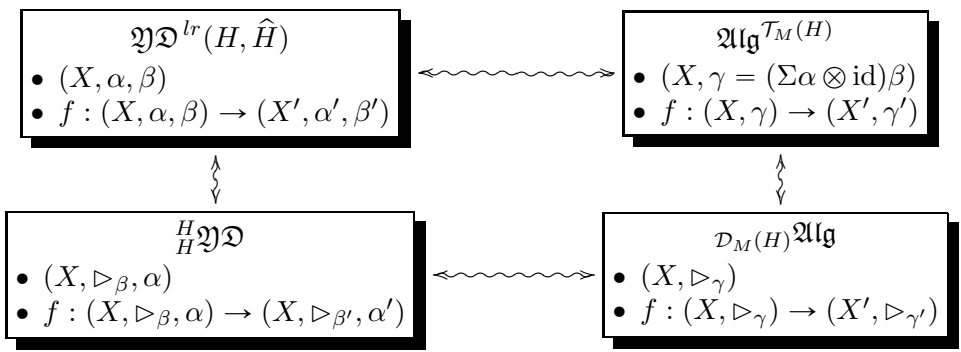
\begin{equation*}\label{eq:right_equivalence_majid}
\xymatrix@C=5pc@R=1.5pc{%
\shadowbox*{%
\begin{Bcenter}
$\;\aYD^{\,lr}(H,\du{H})\;$ \\[0.1cm]
\begin{Bitemize}
\item $(X,\alp,\beta)$
\item $f: (X,\alp,\beta) \to (X',\alp',\beta')$
\end{Bitemize}
\end{Bcenter}}\ar@{<~>}[r]\ar@{<~>}[d]%
&%
\shadowbox*{%
\begin{Bcenter}
$\;\alg^{\T_{M}(H)}\;$ \\[0.1cm]
\begin{Bitemize}
\item $(X,\gamma=(\Sigma\alp \odo \id)\beta)$
\item $f: (X,\gamma) \to (X',\gamma')$
\end{Bitemize}
\end{Bcenter}}\ar@{<~>}[d]%
\\
{\shadowbox*{%
\begin{Bcenter}
$\;{}_{H}^{H}\aYD\;$ \\[0.1cm]
\begin{Bitemize}
\item $(X,\rhd_{\beta},\alp)$
\item $f: (X,\rhd_{\beta},\alp) \to (X,\rhd_{\beta'},\alp')$
\end{Bitemize}
\end{Bcenter}}}\ar@{<~>}[r]%
&%
\shadowbox*{%
\begin{Bcenter}
$\;{}_{\D_{M}(H)}\alg\;$ \\[0.1cm]
\begin{Bitemize}
\item $(X,\rhd_{\gamma})$
\item $f: (X,\rhd_{\gamma}) \to (X',\rhd_{\gamma'})$
\end{Bitemize}
\end{Bcenter}}}
\end{equation*}
\end{prop}
\pr
Fix an element $(X,\alp,\beta) \in \aYD^{lr}(H,\du{H})$, and consider $\gamma = (\Sigma\alp \odo \id_{\du{H}})\beta : X \to X \odo H^{\co} \odo \du{H}$. We have
\begin{align*}
(\gamma \odo \id_{\T_M(H)})\gamma & = ((\Sigma\alp \odo \id_{\du{H}})\beta \odo \id_{\T_{M}(H)})(\Sigma\alp \odo \id_{\du{H}})\beta \\
& = (\Sigma\alp \odo \id_{\du{H}} \odo \id_{\T_{M}(H)})(\beta \odo \id_{\T_{M}(H)})(\Sigma\alp \odo \id_{\du{H}})\beta \\
& = (\Sigma\alp \odo \id_{\du{H}} \odo \id_{\T_{M}(H)})(\beta \odo \id_{\T_{M}(H)})(\Sigma \odo \id_{\du{H}})(\alp \odo \id_{\du{H}})\beta \\
& = (\Sigma\alp \odo \id_{\du{H}} \odo \id_{\T_{M}(H)})\Sigma_{23}\Sigma_{12}(\id \odo \beta \odo \id)(\alp \odo \id_{\du{H}})\beta \\
& = (\Sigma\alp \odo \id_{\du{H}} \odo \id_{\T_{M}(H)})\Sigma_{23}\Sigma_{12}\ad(U^{-1}_{13})(\alp \odo \id \odo \id)(\beta \odo \id_{\du{H}})\beta \\
& = (\Sigma\alp \odo \id_{\du{H}} \odo \id_{\T_{M}(H)})\ad((\Sigma(U^{-1}))_{23})\Sigma_{23}\Sigma_{12}(\alp \odo \id \odo \id)(\id \odo \dcom)\beta \\
& = \ad((\Sigma(U^{-1}))_{34})\Sigma_{34}(\Sigma\alp \odo \id_{H} \odo \id_{\T_{M}(H)})(\Sigma\alp \odo \id \odo \id)(\id \odo \dcom)\beta \\
& = \ad((\Sigma(U^{-1}))_{34})\Sigma_{34}(\id \odo \com^{\co} \odo \id_{\T_{M}(H)})(\Sigma\alp \odo \id \odo \id)(\id \odo \dcom)\beta \\
& = \Sigma_{34}\ad(U^{-1}_{34})(\id \odo \com^{\co} \odo \dcom)(\Sigma\alp \odo \id_{\du{H}})\beta\\
& = (\id_{X} \odo \Sigma_{23}\ad(U^{-1}_{23})(\com^{\co} \odo \dcom))\gamma \\
& = (\id_{X} \odo (\id_{H} \odo \Sigma \odo \id_{\du{H}})(\id_{H} \odo \ad(U^{-1}) \odo \id_{\du{H}})(\com^{\co} \odo \dcom))\gamma \\
& = (\id_{X}\odo \com_{\T_{M}(H)})\gamma \\
\end{align*}
and follows from the injectivity of $\alp$ and $\beta$ that $\gamma$ is injective. Then $(X,\gamma)$ is a right $\T_{M}(H)$-comodule algebra. Moreover, if we compute the left action of $\D_{M}(H)$ on $X$ induced by $\gamma$, we have
\begin{align*}
(\ome^{\op} \bicroi h) \rhd_{\gamma} x & = (\id \odo (\tilde{\bar{\p}} \odo \p)(\ome^{\op} \bicroi h, \cdot))(\gamma(x)) \\
& = (\id \odo (\tilde{\bar{\p}} \odo \p)(\ome^{\op} \bicroi h, \cdot))(\Sigma\alp \odo \id)(\beta(x)) \\
& = ((\id \odo \bar{\p}(\ome,\cdot))\Sigma\alp \odo \p(h,\cdot))(\beta(x)) \\
& = (\id \odo \bar{\p}(\ome,\cdot))\Sigma\alp((\id \odo \p(h,\cdot))(\beta(x))) \\
& = (\bar{\p}(\ome,\cdot) \odo \id)\alp(h \rhd_{\beta} x) = (h \rhd_{\beta} x ) \lhd_{\alp} \ome
\end{align*}
for every $\ome \in \du{H}$, $h \in H$ and $x \in X$, i.e. we have that $(X,\rhd_{\gamma})=\mathfrak{E}_{M}(X,\rhd_{\beta},\alp)$ by Theorem~\ref{th:majid_duality}. Now, fix $(X,\gamma) \in \alg^{\T_{M}(H)}$. Then because $(X,\rhd_{\gamma}) \in {}_{\D_{M}(H)}\alg$, there is a left action $\rhd: H \odo X \to X$ and a left coaction $\alp: X \to H \odo X$ such that $\mathfrak{E}_{M}(X,\rhd,\alp)=(X,\rhd_{\gamma})$. Then for every $x \in X$, we have
\begin{align*}
(\Sigma\alp \odo \id)(\beta_{\rhd}(x)) & = (\Sigma\alp \odo \id)(U_{1} \rhd x \odo U_{2}) = \Sigma\alp((1 \bicroi U_{1}) \rhd_{\gamma} x) \odo U_{2} \\
& = (U^{\op}_{2} \bicroi 1) \rhd_{\gamma} ((1 \bicroi U_{1}) \rhd_{\gamma} x) \odo U_{1} \odo U_{2} \\
& = ((U^{\op}_{2} \bicroi U_{1}) \rhd_{\gamma} x) \odo U_{1} \odo U_{2} \\
& = \gamma(x)
\end{align*}
because $U_{\D_{M}(H)} = U^{\op}_{2} \bicroi U_{1} \odo U_{1} \odo U_{2}$ is the canonical element for the canonical pairing $\p_{\D} = \tilde{\bar{\p}} \odo \p: \D_{M}(H) \times \T_{M}(H) \to \ku$, where $U = U_{1} \odo U_{2}$ is the canonical element for the canonical pairing $\p$ of $H$. Finally, the commutative diagram between categories follows from \ref{prop:rr-YD-equiv} and \ref{th:majid_duality}.
\fin

A similar statement exist in the case of right-left Yetter--Drinfeld algebras.

\begin{prop}\label{prop:left_equivalence_majid}
The functor
\[
\begin{array}{rccc}
\mathfrak{E}'_{M} : & \aYD^{\,rl}(H,\du{H}) & \leadsto & {}^{\T_{M}(H)}\alg \\
& (X,\alp,\beta) & \mapsto & (X,\gamma = (\id \odo \beta)\Sigma\alp) \\
& f: (X,\alp,\beta) \to (X',\alp',\beta') & \mapsto & f: (X,\gamma) \to (X,\gamma')
\end{array}
\]
is an equivalence of categories. Moreover, we have the following equivalence of categories:
\begin{equation*}\label{eq:left_equivalence_majid}
\xymatrix@C=5pc@R=1.5pc{%
\shadowbox*{%
\begin{Bcenter}
$\;\aYD^{\,rl}(H,\du{H})\;$ \\[0.1cm]
\begin{Bitemize}
\item $(X,\alp,\beta)$
\item $f: (X,\alp,\beta) \to (X',\alp',\beta')$
\end{Bitemize}
\end{Bcenter}}\ar@{<~>}[r]\ar@{<~>}[d]%
&%
\shadowbox*{%
\begin{Bcenter}
$\;{}^{\T_{M}(H)}\alg\;$ \\[0.1cm]
\begin{Bitemize}
\item $(X,\gamma=(\id \odo \beta)\Sigma\alp$
\item $f: (X,\gamma) \to (X',\gamma')$
\end{Bitemize}
\end{Bcenter}}\ar@{<~>}[d]%
\\
{\shadowbox*{%
\begin{Bcenter}
$\;\aYD{}_{H}^{H}\;$ \\[0.1cm]
\begin{Bitemize}
\item $(X,\lhd_{\beta},\alp)$
\item $f: (X,\lhd_{\beta},\alp) \to (X,\lhd_{\beta'},\alp')$
\end{Bitemize}
\end{Bcenter}}}\ar@{<~>}[r]%
&%
\shadowbox*{%
\begin{Bcenter}
$\;\alg_{\D_{M}(H)}\;$ \\[0.1cm]
\begin{Bitemize}
\item $(X,\lhd_{\gamma})$
\item $f: (X,\lhd_{\gamma}) \to (X',\lhd_{\gamma'})$
\end{Bitemize}
\end{Bcenter}}}
\end{equation*}
\end{prop}

\subsection{Radford's version}

We recall the original construction of the Drinfeld double as defined by Drinfeld in \cite{D87}. We use the characterization as was given by Radford in \cite{R93}. Let $\D_{R}(H) = \du{H}^{\co} \bicroi H$ the double crossproduct with product and coproduct given by
\[
(\ome \bicroi h)(\theta \bicroi g) = \ome(h_{(1)} \brhd \theta_{(3)} \blhd S^{-1}(h_{(2)})) \bicroi (S^{-1}(\theta_{(1)}) \brhd h_{(3)} \blhd \theta_{(2)})g
\]
and
\[
\com_{\D_{R}(H)}(\ome \bicroi h) = \ome_{(2)} \bicroi h_{(1)} \odo \ome_{(1)} \bicroi h_{(2)},
\]
respectively. Its dual Hopf algebra is given by $\T_{R}(H) = H^{\op} \odo \du{H}$ with algebra and coalgebra structure given by
\[
(h^{\op} \odo \ome)(g^{\op} \odo \theta) = (gh)^{\op} \odo \ome\theta
\]
and
\[
\com_{\T_{R}(H)}(h^{\op} \odo \ome) = (\id_{H^{\op}} \odo \Sigma\ad(U^{-1}) \odo \id_{\du{H}})(h^{\op}_{(1)} \odo h^{\op}_{(2)} \odo \ome_{(1)} \odo \ome_{(2)}),
\]
respectively. It is known in the standard setting that give a left-right (resp. right-left) Yetter--Drinfeld algebra over a Hopf algebra $H$ is equivalent to give an algebra endowed with an left (resp. right) action of $\D_{R}(H)$.

\begin{theo}[\cite{R93}]\label{th:radford_duality}
Let $H$ be a finite-dimensional Hopf algebra. There is an equivalence between the category of (left-right) Yetter--Drinfeld algebras over $H$ and the category of left $\D_{R}(H)$-module algebras. More explicitly, given $(X,\rhd,\alp) \in {}_{H}\aYD^{H}$, we have a structure of left $\D_{R}(H)$-module algebra on $X$ given by
\[
(\ome \bicroi h) \cdot x = (\id \odo \p(\cdot,\ome))(\alp(h \rhd x)) = \ome \rhd_{\alp}(h \rhd x).
\]
for all $h \in H$, $\ome \in \du{H}$ and $x \in X$. Reciprocally, given $(X,\cdot) \in {}_{\D_{R}(H)}\alg$, we have a left action of $H$ and a right coaction of $H^{\op}$ on $X$ given by
\[
h \rhd x := (1 \bicroi h) \cdot x \quad \text{ and } \quad \alp(x) = (U_{2} \bicroi 1) \cdot x \odo U^{\op}_{1},
\]
for every $h \in H$ and $x \in X$, respectively (Here $U_{1} \odo U_{2} := U \in H \odo \du{H}$ is the canonical element associate to the Hopf algebra $H$). We have that $(X,\rhd,\alp)$ is a (left-right) Yetter--Drinfeld algebra over $H$. The equivalence of categories is given by the functor
\[
\begin{array}{rccc}
\mathfrak{E}_{R} : & {}_{H}\aYD^{H} & \leadsto & {}_{\D_{R}(H)}\alg \\
& (X,\rhd,\alp) & \mapsto & (X,\cdot) \\
& f: (X,\rhd,\alp) \to (X',\rhd',\alp') & \mapsto & f: (X,\cdot) \to (X,\cdot')
\end{array}
\]
\end{theo} 

The ``right side'' version of the above theorem is the following:

\begin{theo}[\cite{R93}]\label{th:radford_duality_right}
Let $H$ be a finite-dimensional Hopf algebra. There is an equivalence between the category of (right-left) Yetter--Drinfeld algebras over $H$ and the category of right $\D_{R}(H)$-module algebras. More explicitly, given $(X,\lhd,\alp) \in {}^{H}\aYD_{H}$, we have a structure of right $\D_{R}(H)$-module algebra on $X$ given by
\[
x \cdot (\ome \bicroi h) = (\p(\cdot,\ome) \odo \id)(\alp(x)) \lhd h = (x \lhd_{\alp} \ome) \lhd h.
\]
for all $h \in H$, $\ome \in \du{H}$ and $x \in X$. Reciprocally, given $(X,\cdot) \in \alg_{\D_{R}(H)}$, we have a right action of $H$ and a left coaction of $H^{\op}$ on $X$ given by
\[
x \lhd h := x \cdot (1 \bicroi h) \quad \text{ and } \quad \alp(x) = U^{\op}_{1} \odo x \cdot (U_{2} \bicroi 1),
\]
for every $h \in H$ and $x \in X$, respectively (Here $U_{1} \odo U_{2} := U \in H \odo \du{H}$ is the canonical element associate to the Hopf algebra $H$). We have that $(X,\lhd,\alp)$ is a right-left Yetter--Drinfeld algebra over $H$. The equivalence of categories is given by the functor
\[
\begin{array}{rccc}
\mathfrak{E}_{R} : & {}^{H}\aYD_{H} & \leadsto & \alg_{\D_{R}(H)} \\
& (X,\lhd,\alp) & \mapsto & (X,\cdot) \\
& f: (X,\lhd,\alp) \to (X',\lhd',\alp') & \mapsto & f: (X,\cdot) \to (X,\cdot')
\end{array}
\]
\end{theo}

Now, we prove similar statements in the ``only coaction" characterization.


\begin{prop}\label{prop:right_equivalence_radford}
The functor
\[
\begin{array}{rccc}
\mathfrak{E}'_{R} : & \aYD^{\,rr}(H,\du{H}) & \leadsto & \alg^{\T_{R}(H)} \\
& (X,\alp,\beta) & \mapsto & (X,\gamma = (\alp \odo \id)\beta) \\
& f: (X,\alp,\beta) \to (X',\alp',\beta') & \mapsto & f: (X,\gamma) \to (X,\gamma')
\end{array}
\]
is an equivalence of categories. Moreover, we have the following equivalence of categories:
\begin{equation*}\label{eq:right_equivalence_radford}
\xymatrix@C=5pc@R=1.5pc{%
\shadowbox*{%
\begin{Bcenter}
$\;\aYD^{\,rr}(H,\du{H})\;$ \\[0.1cm]
\begin{Bitemize}
\item $(X,\alp,\beta)$
\item $f: (X,\alp,\beta) \to (X',\alp',\beta')$
\end{Bitemize}
\end{Bcenter}}\ar@{<~>}[r]\ar@{<~>}[d]%
&%
\shadowbox*{%
\begin{Bcenter}
$\;\alg^{\T_{R}(H)}\;$ \\[0.1cm]
\begin{Bitemize}
\item $(X,\gamma=(\alp \odo \id)\beta)$
\item $f: (X,\gamma) \to (X',\gamma')$
\end{Bitemize}
\end{Bcenter}}\ar@{<~>}[d]%
\\
{\shadowbox*{%
\begin{Bcenter}
$\;{}_{H}\aYD^{H}\;$ \\[0.1cm]
\begin{Bitemize}
\item $(X,\rhd_{\beta},\alp)$
\item $f: (X,\rhd_{\beta},\alp) \to (X,\rhd_{\beta'},\alp')$
\end{Bitemize}
\end{Bcenter}}}\ar@{<~>}[r]%
&%
\shadowbox*{%
\begin{Bcenter}
$\;{}_{\D_{R}(H)}\alg\;$ \\[0.1cm]
\begin{Bitemize}
\item $(X,\rhd_{\gamma})$
\item $f: (X,\rhd_{\gamma}) \to (X',\rhd_{\gamma'})$
\end{Bitemize}
\end{Bcenter}}}
\end{equation*}
\end{prop}
\pr
Fix an element $(X,\alp,\beta) \in \aYD^{rr}(H,\du{H})$, and consider the linear map $\gamma = (\alp \odo \id_{\du{H}})\beta : X \to X \odo H^{\op} \odo \du{H}$. We have
\begin{align*}
(\gamma \odo \id_{\T_{R}(H)})\gamma & = ((\alp \odo \id_{\du{H}})\beta \odo \id_{\T_{R}(H)})(\alp \odo \id_{\du{H}})\beta \\
& = (\alp \odo \id \odo \id_{\T_{R}(H)})(\beta \odo \id_{\T_{R}(H)})(\alp \odo \id_{\du{H}})\beta \\
& = (\alp \odo \id \odo \id_{\T_{R}(H)})((\beta \odo \id_{H^{\op}})\alp \odo \id_{\du{H}})\beta \\
& = (\alp \odo \id \odo \id_{\T_{R}(H)})((\id_{X} \odo \Sigma)(\id_{X} \odo \ad({}^{\ops}U))(\alp \odo \id)\beta \odo \id)\beta \\
& = (\alp \odo \id \odo \id_{\T_{R}(H)})(\Sigma \ad({}^{\ops}U))_{23}(\alp \odo \id \odo \id)(\beta \odo \id)\beta \\
& = (\id \odo (\Sigma \ad({}^{\ops}U))_{23})(\alp \odo \id \odo \id_{\T_{R}(H)})(\alp \odo \id \odo \id)(\id \odo \dcom)\beta \\
& = (\id \odo (\Sigma \ad({}^{\ops}U))_{23})(\id \odo \com^{\op} \odo \id_{\T_{R}(H)})(\alp \odo \id_{\T_{R}(H)})(\id \odo \dcom)\beta \\
& = (\id \odo (\Sigma \ad({}^{\ops}U))_{23})(\id \odo \com^{\op} \odo \id_{\T_{R}(H)})(\id \odo \id \odo \dcom)(\alp \odo \id)\beta \\
& = (\id \odo (\Sigma \ad({}^{\ops}U))_{23})(\id \odo \com^{\op} \odo \dcom)(\alp \odo \id)\beta \\
& = (\id \odo (\id \odo \Sigma \odo \id)(\id \odo \ad({}^{\ops}U) \odo \id)(\com^{\op} \odo \dcom))\gamma \\
& = (\id \odo \com_{\T_{R}(H)})\gamma
\end{align*}
and follows from the injectivity of $\alp$ and $\beta$ that $\gamma$ is injective. Then $(X,\gamma)$ is a right $\T_{R}(H)$-comodule algebra. Moreover, if we compute the left action of $\D_{R}(H)$ on $X$ induced by $\gamma$, we have
\begin{align*}
(\ome \bicroi h) \rhd_{\gamma} x & = (\id \odo (\tilde{\p} \odo \p)(\ome \bicroi h, \cdot))(\gamma(x)) \\
& = (\id \odo (\tilde{\p} \odo \p)(\ome \bicroi h, \cdot))(\alp \odo \id)(\beta(x)) \\
& = ((\id \odo \tilde{\p}(\ome,\cdot))\alp \odo \p(h,\cdot))(\beta(x)) \\
& = (\id \odo \tilde{\p}(\ome,\cdot))\alp((\id \odo \p(h,\cdot))(\beta(x))) \\
& = (\id \odo \tilde{\p}(\ome,\cdot))\alp(h \rhd_{\beta} x) = \ome \rhd_{\alp} (h \rhd_{\beta} x)
\end{align*}
for every $\ome \in \du{H}$, $h \in H$ and $x \in X$, i.e. we have that $(X,\rhd_{\gamma})=\mathfrak{E}_{R}(X,\rhd_{\beta},\alp)$ by Theorem~\ref{th:radford_duality}. Now, fix $(X,\gamma) \in \alg^{\T_{R}(H)}$. Then because $(X,\rhd_{\gamma}) \in {}_{\D_{R}(H)}\alg$, there is a left action $\rhd: H \odo X \to X$ and a right coaction $\alp: X \to X \odo H^{\op}$ such that $\mathfrak{E}_{R}(X,\rhd,\alp) = (X,\rhd_{\gamma})$. Then for every $x \in X$, we have
\begin{align*}
(\alp \odo \id)(\beta_{\rhd}(x)) & = (\alp \odo \id)(U_{1} \rhd x \odo U_{2}) = \alp((1 \bicroi U_{1}) \rhd_{\gamma} x) \odo U_{2}\\
& = (U_{2} \bicroi 1) \rhd_{\gamma} ((1 \bicroi U_{1}) \rhd_{\gamma} x) \odo U^{\op}_{1} \odo U_{2} \\
& = ((U_{2} \bicroi U_{1}) \rhd_{\gamma} x) \odo U^{\op}_{1} \odo U_{2} \\
& = \gamma(x),
\end{align*}
because $U_{\D_{R}(H)} = U_{2} \bicroi U_{1} \odo U^{\op}_{1} \odo U_{2}$ is the canonical element for the canonical pairing $\p_{\D} = \tilde{\p} \odo \p: \D_{R}(H) \times \T_{R}(H) \to \ku$, where $U = U_{1} \odo U_{2}$ is the canonical element for the canonical pairing $\p$ of $H$. Finally, the commutative diagram between categories follows from \ref{prop:rr-YD-equiv} and \ref{th:radford_duality}.
\fin

\begin{prop}\label{prop:left_equivalence_radford}
The functor
\[
\begin{array}{rccc}
\mathfrak{E}'_{R} : & \aYD^{\,ll}(H,\du{H}) & \leadsto & {}^{\T_{R}(H)}\alg \\
& (X,\alp,\beta) & \mapsto & (X,\gamma = (\id \odo \beta)\alp) \\
& f: (X,\alp,\beta) \to (X',\alp',\beta') & \mapsto & f: (X,\gamma) \to (X,\gamma')
\end{array}
\]
is an equivalence of categories. Moreover, we have the following equivalence of categories:
\begin{equation*}\label{eq:left_equivalence_radford}
\xymatrix@C=5pc@R=1.5pc{%
\shadowbox*{%
\begin{Bcenter}
$\;\aYD^{\,ll}(H,\du{H})\;$ \\[0.1cm]
\begin{Bitemize}
\item $(X,\alp,\beta)$
\item $f: (X,\alp,\beta) \to (X',\alp',\beta')$
\end{Bitemize}
\end{Bcenter}}\ar@{<~>}[r]\ar@{<~>}[d]%
&%
\shadowbox*{%
\begin{Bcenter}
$\;{}^{\T_{R}(H)}\alg\;$ \\[0.1cm]
\begin{Bitemize}
\item $(X,\gamma=(\id \odo \beta)\alp)$
\item $f: (X,\gamma) \to (X',\gamma')$
\end{Bitemize}
\end{Bcenter}}\ar@{<~>}[d]%
\\
{\shadowbox*{%
\begin{Bcenter}
$\;{}^{H}\aYD_{H}\;$ \\[0.1cm]
\begin{Bitemize}
\item $(X,\lhd_{\beta},\alp)$
\item $f: (X,\lhd_{\beta},\alp) \to (X,\lhd_{\beta'},\alp')$
\end{Bitemize}
\end{Bcenter}}}\ar@{<~>}[r]%
&%
\shadowbox*{%
\begin{Bcenter}
$\;\alg_{\D_{R}(H)}\;$ \\[0.1cm]
\begin{Bitemize}
\item $(X,\lhd_{\gamma})$
\item $f: (X,\lhd_{\gamma}) \to (X',\lhd_{\gamma'})$
\end{Bitemize}
\end{Bcenter}}}
\end{equation*}
\end{prop}
\pr
Fix an element $(X,\alp,\beta) \in \aYD^{ll}(H,\du{H})$, and consider the linear  map $\gamma = (\id_{H^{\op}} \odo \beta)\alp: X \to H^{\op} \odo \du{H} \odo X$. We have
\begin{align*}
(\id_{\T_{R}(H)} \odo \gamma)\gamma & = (\id_{\T_{R}(H)} \odo (\id_{H^{\op}} \odo \beta)\alp)(\id_{H^{\op}} \odo \beta)\alp \\
& = (\id_{\T_{R}(H)} \odo \id \odo \beta)(\id_{\T_{R}(H)} \odo \alp)(\id \odo \beta)\alp \\
& = (\id_{\T_{R}(H)} \odo \id \odo \beta)(\id_{H^{\op}} \odo (\id \odo \alp)\beta)\alp \\
& = (\id_{\T_{R}(H)} \odo \id \odo \beta)(\id_{H^{\op}} \odo (\Sigma \odo \id)(\ad({}^{\ops}U) \odo \id)(\id \odo \beta)\alp)\alp \\
& = ((\Sigma \ad({}^{\ops}U))_{23} \odo \beta)(\id \odo \id \odo \beta)(\id \odo \alp)\alp \\
& = ((\Sigma \ad({}^{\ops}U))_{23} \odo \beta)(\id \odo \id \odo \beta)(\com^{\op} \odo \id)\alp \\
& = ((\Sigma \ad({}^{\ops}U))_{23} \odo \id)(\com^{\op} \odo \id \odo \beta)(\id \odo \beta)\alp \\
& = ((\Sigma \ad({}^{\ops}U))_{23} \odo \id)(\com^{\op} \odo \dcom \odo \id)(\id \odo \beta)\alp \\
& = ((\id \odo \Sigma \odo \id)(\id \odo \ad({}^{\ops}U) \odo \id)(\com^{\op} \odo \dcom) \odo \id)(\id \odo \beta)\alp \\
& = (\com_{\T_{R}(H)} \odo \id)\gamma
\end{align*}
and follows from the injectivity of $\alp$ and $\beta$ that $\gamma$ is injective. Then $(X,\gamma)$ is a left $\T_{R}(H)$-comodule algebra. Moreover, if we compute the right action of $\D_{R}(H)$ on $X$ induced by $\gamma$, we have
\begin{align*}
x \lhd_{\gamma} (\ome \bicroi h) & = ((\tilde{\p} \odo \p)(\ome \bicroi h, \cdot) \odo \id )(\gamma(x)) \\
& = ((\tilde{\p} \odo \p)(\ome \bicroi h, \cdot) \odo \id)(\id \odo \beta)(\alp(x)) \\
& = (\tilde{\p}(\ome,\cdot) \odo (\p(h,\cdot) \odo \id)\beta)(\alp(x)) \\
& = (\p(h,\cdot) \odo \id)\beta)(x \lhd_{\alp} \ome) \\
& = (x \lhd_{\alp} \ome) \lhd_{\beta} h
\end{align*}
for every $\ome \in \du{H}$, $h \in H$ and $x \in X$, i.e. we have that $(X,\lhd_{\gamma})=\mathfrak{E}_{R}(X,\lhd_{\beta},\alp)$ by Theorem~\ref{th:radford_duality_right}. Now, fix $(X,\gamma) \in {}^{\T_{R}(H)}\alg$. Then because $(X,\lhd_{\gamma}) \in \alg_{\D_{R}(H)}$, there is a right action $\lhd: X \odo H \to X$ and a left coaction $\alp: X \to H^{\op} \odo X$ such that $\mathfrak{E}_{R}(X,\lhd,\alp) = (X,\lhd_{\gamma})$. Then for every $x \in X$, we have
\begin{align*}
(\id \odo \beta_{\lhd})(\alp(x)) & = (\id \odo \beta_{\lhd})(U^{\op}_{1} \odo x \lhd_{\gamma} (U_{2} \bicroi 1)) = U^{\op}_{1} \odo \beta_{\lhd}(x \lhd_{\gamma} (U_{2} \bicroi 1))\\
& = U^{\op}_{1} \odo U_{2} \odo (x \lhd_{\gamma} (U_{2} \bicroi 1)) \lhd_{\gamma} (1 \bicroi U_{1}) \\
& = U^{\op}_{1} \odo U_{2} \odo (x \lhd_{\gamma} (U_{2} \bicroi U_{1})) \\
& = \gamma(x),
\end{align*}
because $U_{\D_{R}(H)} = U_{2} \bicroi U_{1} \odo U^{\op}_{1} \odo U_{2}$ is the canonical element for the canonical pairing $\p_{\D} = \tilde{\p} \odo \p: \D_{R}(H) \times \T_{R}(H) \to \ku$, where $U = U_{1} \odo U_{2}$ is the canonical element for the canonical pairing $\p$ of $H$. Finally, the commutative diagram between categories follows from \ref{prop:ll-YD-equiv} and \ref{th:radford_duality_right}.
\fin

\section{Equivalence of Yetter--Drinfeld categories}\label{sec:equi_yd_cat}

\subsection{Standard characterization}

An adapted formulation of the equivalence obtained for Yetter--Drinfeld modules in \cite{RT93}, is given for the standard characterization of Yetter--Drinfeld algebras over finite-dimensional Hopf algebras.

Given a left action $\rhd: H \odo X \to X$, we define the map $\lhd_{S^{-1}}: X \odo H^{\co} \to X$, $x \lhd_{S^{-1}} h := S^{-1}(h) \rhd x$. Similarly, given a right action $\lhd: X \odo H \to X$, we define the map ${}_{S^{-1}}\rhd: H^{\co}\odo X \to X$, $h \rhd_{S^{-1}} x := x \lhd S^{-1}(h)$. Given a left coaction $\delta: X \to H^{\op} \odo X$, $x \mapsto x^{\op}_{[-1]} \odo x_{[0]}$, we define the map ${}_{S}\delta: X \to X \odo H$, ${}_{S}\delta(x) = x_{[0]} \odo S(x_{[-1]}) = \Sigma(S \odo \id)({}^{\op} \odo \id)\delta(x) = \Sigma(S \odo \id)\dot{\delta}(x)$. Similarly, given a right coaction $\delta: X \to X \odo H^{\op}$, $x \mapsto x_{[0]} \odo x^{\op}_{[1]}$, we define the map $\delta_{S}: X \to H \odo X$, $\delta_{S}(x) = S(x_{[1]}) \odo x_{[0]} = \Sigma(\id \odo S)(\id \odo {}^{\op})\delta(x) = \Sigma(\id \odo S)\dot{\delta}(x)$.

\begin{rema}
Given a right coaction $\delta: X \to X \odo H^{\op}$, $x \mapsto x_{[0]} \odo x^{\op}_{[1]}$, the map $\delta_{S}: X \to H \odo X$, $\delta_{S} = \Sigma(\id \odo S)\dot{\delta}$, is a left coaction of $H$ on $X$. Indeed, we have
\begin{align*}
(\id \odo \delta_{S})\delta_{S} & = (\id \odo \Sigma(\id \odo S)\dot{\delta})\Sigma(\id \odo S)\dot{\delta} = \Sigma_{23}(\id \odo \id \odo S)\Sigma_{12}\Sigma_{23}(\dot{\delta} \odo \id)(\id \odo S)\dot{\delta} \\
& = \Sigma_{13}(\id \odo S \odo S)(\dot{\delta} \odo \id)\dot{\delta} = \Sigma_{23}\Sigma_{12}\Sigma_{23}(\id \odo S \odo S)(\id \odo \com_{H})\dot{\delta} \\
& = \Sigma_{23}\Sigma_{12}(\id \odo \com_{H})(\id \odo S)\dot{\delta} = (\com_{H} \odo \id)\Sigma(\id \odo S)\dot{\delta} = (\com_{H} \odo \id)\delta_{S}, \\\\
(\cou_{H} \odo \id)\delta_{S} & = (\cou_{H} \odo \id)\Sigma(\id \odo S)\dot{\delta} = (\id \odo \cou_{H}S)\dot{\delta} = (\id \odo \cou_{H})\dot{\delta} = \id
\end{align*}
and
\begin{align*}
\delta_{S}(xy) & = \Sigma(\id \odo S)\dot{\delta}(xy) = \Sigma(\id \odo S)(\id \odo {}^{\op})(\delta(x)\delta(y)) = \Sigma(\id \odo S)(x_{[0]}y_{[0]} \odo y_{[1]}x_{[1]}) \\
& = (S(x_{[1]}) \odo x_{[0]})(S(y_{[1]}) \odo y_{[0]}) = \delta_{S}(x)\delta_{S}(y)
\end{align*}
for any $x,y \in X$.
\end{rema}


\begin{rema}
We have the following equivalences between actions and coactions.
\begin{equation*}\label{eq:actions_coactions_1}
\xymatrix@C=5pc@R=1.5pc{%
\shadowbox*{%
\begin{Bcenter}
$\;\alg^{H}\;$ \\[0.1cm]
$\begin{array}{lccc}
\delta: & X & \to & X \odo H \\
& x & \mapsto & x_{[0]} \odo x_{[1]}
\end{array}$
\end{Bcenter}}\ar@{<~>}[r]\ar@{<~>}[d]%
&%
\shadowbox*{%
\begin{Bcenter}
$\;{}^{H^{\co}}\alg\;$ \\[0.1cm]
$\begin{array}{lccc}
\Sigma\delta: & X & \to & H^{\co} \odo X \\
& x & \mapsto & x_{[1]} \odo x_{[0]}
\end{array}$
\end{Bcenter}}\ar@{<~>}[d]%
\\
{\shadowbox*{%
\begin{Bcenter}
$\;{}_{\du{H}}\alg\;$ \\[0.1cm]
$\begin{array}{lccc}
\rhd_{\delta}: & \du{H} \odo X  & \to & X \\
& \ome \odo x & \mapsto & \p(x_{[1]},\ome)x_{[0]}
\end{array}$
\end{Bcenter}}}\ar@{<~>}[r]%
&%
\shadowbox*{%
\begin{Bcenter}
$\;\alg_{\du{H}^{\op}}\;$ \\[0.1cm]
$\begin{array}{lccc}
\dot{\lhd}_{\delta} : & X \odo \du{H}^{\op} & \to & X \\
& x \odo \ome^{\op} & \mapsto & \ome \rhd_{\delta} x
\end{array}$
\end{Bcenter}}}
\end{equation*}

\begin{equation*}\label{eq:actions_coactions_2}
\xymatrix@C=5pc@R=1.5pc{%
\shadowbox*{%
\begin{Bcenter}
$\;\alg{}^{H}\;$ \\[0.1cm]
$\begin{array}{lccc}
\delta: & X & \to & X \odo H \\
& x & \mapsto & x_{[0]} \odo x_{[1]}
\end{array}$
\end{Bcenter}}\ar@{<~>}[r]\ar@{<~>}[d]%
&%
\shadowbox*{%
\begin{Bcenter}
$\;{}^{H^{\op}}\alg\;$ \\[0.1cm]
$\begin{array}{lccc}
\delta_{S^{\ops}}: & X & \to & H^{\op} \odo X \\
& x & \mapsto & S^{-1}(x_{[-1]})^{\op} \odo x_{[0]}
\end{array}$
\end{Bcenter}}\ar@{<~>}[d]%
\\
{\shadowbox*{%
\begin{Bcenter}
$\;{}_{\du{H}}\alg\;$ \\[0.1cm]
$\begin{array}{lccc}
\rhd_{\delta}: & \du{H} \odo X & \to & X \\
& \ome \odo x
& \mapsto & \p(x_{[1]},\ome)x_{[0]}
\end{array}$
\end{Bcenter}}}\ar@{<~>}[r]%
&%
\shadowbox*{%
\begin{Bcenter}
$\;\alg_{\du{H}^{\co}}\;$ \\[0.1cm]
$\begin{array}{lccc}
{}_{\dS^{-1}}(\lhd_{\delta}) : & X \odo \du{H}^{\co} & \to & X \\
& x \odo \ome & \mapsto & \dS^{-1}(\ome) \rhd_{\delta} x
\end{array}$
\end{Bcenter}}}
\end{equation*}
\end{rema}

\begin{rema}
Given a left action $\rhd: H \odo X \to X$. Consider $\dot{\lhd}: X \odo H^{\op} \to X$, $x \dot{\lhd} h^{\op} := h \rhd x$. Take the map $S: H^{\op} \to H^{\co}$, $h^{\op} \mapsto S(h)$ which is a Hopf algebra isomorphism, we have $\lhd_{S^{-1}}(\id \odo S) = \dot{\lhd}$, i.e. we have the following commutative diagram
\[
\xymatrix{X \odo H^{\op}\ar[rr]^-{\id \odo S}\ar[rd]_-{\dot{\lhd}} & & X \odo H^{\co}\ar[ld]^-{\lhd_{S^{-1}}} \\ & X &}
\]

Moreover, it can be prove that the map $S \# \id: H^{\op}\# X \to H^{\co}\# X$, $h^{\op}\# x \mapsto S(h)\# x$ is a canonical isomorphism between the smash products.
\end{rema}

\begin{rema}
For the right coaction $\delta: X \to X \odo H$, $x \mapsto x_{[0]} \odo x_{[1]}$, the map $\delta_{S^{\ops}}: X \to H^{\op} \odo X$ is given by $\delta_{S^{\ops}}(x) = \Sigma(\id \odo S^{\ops})\dot{\delta}(x) = S^{-1}(x_{[1]})^{\op} \odo x_{[0]}$. This follows because the antipode on the Hopf algebra $H^{\op}$ is given by $S^{\ops}(h^{\op}) := S^{-1}(h)^{\op}$ for all $h\in H$.
\end{rema}

\begin{prop}\label{prop:duality_equivalence}
Let $H$ be a finite-dimensional Hopf algebra. The functors
\[
\begin{array}{rccc}
\mathfrak{F}_{H} : & {}_{H}\aYD^{H} & \leadsto & {}^{H^{\op}}\aYD_{H^{\op}} \\
& (X,\rhd,\delta) & \mapsto & (X,\dot{\lhd},\delta_{S}) \\
& f: (X,\rhd,\delta) \to (X',\rhd',\delta') & \mapsto & f: (X,\dot{\lhd},\delta_{S}) \to (X,\dot{\lhd}',\delta'_{S})
\end{array}
\]
\[
\begin{array}{rccc}
\mathfrak{F}'_{H} : & {}^{H}\aYD_{H} & \leadsto & {}_{H^{\op}}\aYD^{H^{\op}} \\
& (X,\lhd,\delta) & \mapsto & (X,\dot{\rhd},{}_{S}\delta) \\
& f: (X,\lhd,\alp) \to (X',\lhd',\alp') & \mapsto & f: (X,\dot{\rhd},{}_{S}\alp) \to (X,\dot{\rhd}',{}_{S}\alp')
\end{array}
\]
are equivalence of categories. Moreover, $\mathfrak{F}'_{H^{\op}}\mathfrak{F}_{H} = \mathfrak{Id}_{{}_{H}\aYD^{H}}$ and $\mathfrak{F}_{H^{\op}}\mathfrak{F}'_{H} = \mathfrak{Id}_{{}^{H}\aYD_{H}}$.
\end{prop}
\pr
Let $(X,\rhd,\delta) \in {}_{H}\aYD^{H}$, then we have $\delta(h \rhd x) = (h_{(2)} \rhd x_{[0]}) \odo (h_{(3)}x_{[1]}S^{-1}(h_{(1)}))^{\op}$ for all $h \in H$, $x \in X$. Then, it is straighforward, see that
\[
\delta_{S}(x \dot{\lhd} h^{\op}) = h_{(1)}S(x_{[1]})S(h_{(3)}) \odo (h_{(2)} \rhd x_{[0]}) = ((S^{\ops})^{-1}(h^{\op}_{(3)})S(x_{[1]})^{\op}h^{\op}_{(1)})^{\op} \odo (x_{[0]} \dot{\lhd} h^{\op}_{(2)})
\]
for all $h \in H$, $x \in X$. This implies $(X,\dot{\lhd},\delta_{S}) \in {}^{H^{\op}}\aYD_{H^{\op}}$. Similarly, if $(X,\lhd,\delta) \in {}^{H}\aYD_{H}$, then we have $\delta(x \lhd h) = (S^{-1}(h_{(3)})x_{[-1]}h_{(1)})^{\op} \odo (x_{[0]} \lhd h_{(2)})$ for all $h \in H$, $x \in X$. Then, it is straighforward, see that
\[
{}_{S}\delta(h^{\op} \dot{\rhd} x) = (x_{[0]} \lhd h_{(2)}) \odo S(h_{(1)})S(x_{[-1]})h_{(3)} = (h^{\op}_{(2)} \dot{\rhd} x_{[0]}) \odo (h^{\op}_{(3)}S(x_{[1]})^{\op}(S^{\ops})^{-1}(h^{\op}_{(1)}))^{\op}
\]
for all $h \in H$, $x \in X$. This implies $(X,\dot{\lhd},\delta_{S}) \in {}^{H^{\op}}\aYD_{H^{\op}}$. To check the last statement of the proposition, note that
\[
{}_{S^{\ops}}(\delta_{S}) = \Sigma( S^{\ops} \odo \id)\dot{\delta_{S}} = \Sigma({}^{\op} \odo id)(S^{-1} \odo \id)\Sigma(\id \odo S)\dot{\delta} = (\id \odo {}^{\op})\dot{\delta} = \delta
\]
for a right coaction $\delta: X \to X \odo H$.
\fin

\begin{prop}\label{prop:duality_equivalence_cop}
Let $H$ be a finite-dimensional Hopf algebra. The functors
\[
\begin{array}{rccc}
\mathfrak{G}_{H} : & {}_{H}\aYD^{H} & \leadsto & {}^{H^{\co}}\aYD_{H^{\co}} \\
& (X,\rhd,\delta) & \mapsto & (X,\lhd_{S},\Sigma\delta) \\
& f: (X,\rhd,\delta) \to (X',\rhd',\delta') & \mapsto & f: (X,\lhd_{S},\Sigma\delta) \to (X,\lhd'_{S},\Sigma\delta')
\end{array}
\]
\[
\begin{array}{rccc}
\mathfrak{G}'_{H} : & {}^{H}\aYD_{H} & \leadsto & {}_{H^{\co}}\aYD^{H^{\co}} \\
& (X,\lhd,\delta) & \mapsto & (X,{}_{S}\rhd,\Sigma\delta) \\
& f: (X,\lhd,\alp) \to (X',\lhd',\alp') & \mapsto & f: (X,{}_{S}\rhd,\Sigma\delta) \to (X,{}_{S}\rhd',\Sigma\delta')
\end{array}
\]
are equivalence of categories. Moreover, $\mathfrak{G}'_{H^{\co}}\mathfrak{G}_{H} = \mathfrak{Id}_{{}_{H}\aYD^{H}}$ and $\mathfrak{G}_{H^{\co}}\mathfrak{G}'_{H} = \mathfrak{Id}_{{}^{H}\aYD_{H}}$.
\end{prop}


The other variants of the above proposition are the following

\begin{prop}
Let $H$ be a finite-dimensional Hopf algebra. The functors
\[
\begin{array}{rccc}
\mathfrak{F}_{H} : & {}_{H}\aYD^{H} & \leadsto & {}^{H}_{H}\aYD \\
& (X,\rhd,\delta) & \mapsto & (X,\rhd,\delta_{S}) \\
& f: (X,\rhd,\delta) \to (X',\rhd',\delta') & \mapsto & f: (X,\rhd,\delta_{S}) \to (X,\rhd',\delta'_{S})
\\\\
\mathfrak{F}'_{H} : & {}^{H}_{H}\aYD & \leadsto & {}_{H}\aYD^{H} \\
& (X,\rhd,\delta) & \mapsto & (X,\rhd,{}_{S^{-1}}\delta) \\
& f: (X,\rhd,\delta) \to (X',\rhd',\delta') & \mapsto & f: (X,\rhd,{}_{S^{-1}}\delta) \to (X,\rhd',{}_{S^{-1}}\delta')
\\\\
\mathfrak{G}_{H} : & {}_{H}\aYD^{H} & \leadsto & \aYD^{H^{\op}}_{H^{\op}} \\
& (X,\rhd,\delta) & \mapsto & (X,\dot{\lhd},\delta) \\
& f: (X,\rhd,\delta) \to (X',\rhd',\delta') & \mapsto & f: (X,\dot{\lhd},\delta) \to (X,\dot{\lhd}',\delta')
\\\\
\mathfrak{G}'_{H} : & \aYD^{H}_{H} & \leadsto & {}_{H^{\op}}\aYD^{H^{\op}} \\
& (X,\lhd,\delta) & \mapsto & (X,\dot{\rhd},\delta) \\
& f: (X,\rhd,\delta) \to (X',\rhd',\delta') & \mapsto & f: (X,\dot{\rhd},\delta) \to (X,\dot{\rhd}',\delta')
\end{array}
\]
are equivalence of categories. Moreover, $\mathfrak{F}'_{H}\mathfrak{F}_{H} = \mathfrak{Id}_{{}_{H}\aYD^{H}}$ and $\mathfrak{F}_{H}\mathfrak{F}'_{H} = \mathfrak{Id}_{{}^{H}_{H}\aYD}$.
\end{prop}
\pr
It is straightforward.
\fin

\begin{prop}
Let $H$ be a finite-dimensional Hopf algebra. The functors
\[
\begin{array}{rccc}
\mathfrak{F}_{H} : & \aYD^{H}_{H} & \leadsto & {}^{H^{\op}}_{H^{\op}}\aYD \\
& (X,\lhd,\delta) & \mapsto & (X,\dot{\rhd},\delta_{S}) \\
& f: (X,\lhd,\delta) \to (X',\lhd',\delta') & \mapsto & f: (X,\dot{\rhd},\delta_{S}) \to (X,\dot{\rhd}',\delta'_{S})
\\\\
\mathfrak{G}_{H} : & \aYD^{H}_{H} & \leadsto & {}^{H^{\co}}_{H^{\co}}\aYD \\
& (X,\lhd,\delta) & \mapsto & (X,{}_{S^{-1}}\rhd,\Sigma\delta) \\
& f: (X,\lhd,\delta) \to (X',\lhd',\delta') & \mapsto & f: (X,{}_{S^{-1}}\rhd,\Sigma\delta) \to (X,{}_{S^{-1}}\rhd',\Sigma\delta')
\\\\
\mathfrak{F}'_{H} : & {}^{H}_{H}\aYD & \leadsto & \aYD^{H^{\op}}_{H^{\op}} \\
& (X,\rhd,\delta) & \mapsto & (X,\dot{\lhd},{}_{S}\delta) \\
& f: (X,\rhd,\delta) \to (X',\rhd',\delta') & \mapsto & f: (X,\dot{\lhd},{}_{S}\delta) \to (X,\dot{\lhd}',{}_{S}\delta')
\\\\
\mathfrak{G}'_{H} : & {}^{H}_{H}\aYD & \leadsto & \aYD^{H^{\co}}_{H^{\co}} \\
& (X,\rhd,\delta) & \mapsto & (X,\lhd_{S^{-1}},\Sigma\delta) \\
& f: (X,\rhd,\delta) \to (X',\rhd',\delta') & \mapsto & f: (X,{}_{S^{-1}}\rhd,\Sigma\delta) \to (X,{}_{S^{-1}}\rhd',\Sigma\delta')
\end{array}
\]
are equivalence of categories. Moreover,
\[
\mathfrak{F}'_{H^{\op}}\mathfrak{F}_{H} = \mathfrak{Id}_{\aYD^{H}_{H}} = \mathfrak{G}'_{H^{\co}}\mathfrak{G}_{H}, \quad \text{ and } \quad \mathfrak{F}_{H^{\op}}\mathfrak{F}'_{H} = \mathfrak{Id}_{{}^{H}_{H}\aYD} = \mathfrak{G}_{H^{\co}}\mathfrak{G}'_{H}.
\]
\end{prop}
\pr
It is straightforward.
\fin

We end this section, given an equivalence between the four versions of a Yetter--Drinfeld algebra in the standard characterization.

\begin{prop}\label{eq:left_equivalence_all_yd}
Let $H$ be a finite-dimensional Hopf algebra. We have the following equivalence of categories:
\begin{equation*}\label{eq:left_equivalence_all_yd}
\xymatrix@C=5pc@R=1.5pc{%
\shadowbox*{%
\begin{Bcenter}
$\;\aYD^{H}_{H}\;$ \\[0.1cm]
$(X,\lhd_{\beta},\alp)$
\end{Bcenter}}\ar@{<~>}[r]\ar@{<~>}[d]%
&%
\shadowbox*{%
\begin{Bcenter}
$\;{}^{\du{H}}_{\du{H}}\aYD\;$ \\[0.1cm]
$(X,\rhd_{\alp},\beta)$
\end{Bcenter}}\ar@{<~>}[d]%
\\
{\shadowbox*{%
\begin{Bcenter}
$\;{}^{H}\aYD_{H}\;$ \\[0.1cm]
$(X,\lhd_{\beta},(S^{-1} \odo \id)\Sigma\alp)$
\end{Bcenter}}}\ar@{<~>}[r]%
&%
\shadowbox*{%
\begin{Bcenter}
$\;{}_{\du{H}}\aYD^{\du{H}}\;$ \\[0.1cm]
$(X,\rhd_{\alp},(\id \odo \du{S}^{-1})\Sigma\beta)$
\end{Bcenter}}}
\end{equation*}
\end{prop}
Moreover, the above equivalence respect the braided commutativity conditions.
\pr
It is straightforward.
\fin

\subsection{``Only coaction'' characterization}

Before to give the statement about the equivalence of categories in the only coaction characterization, we will introduced two definitions in order to give a simplified version of the equivalence. Assume, we have two Hopf algebras $H$, $H'$ and a pairing $\p: H \times H' \to \ku$ which is endowed with a canonical element $U\in H\odo H'$ such that $\p^{2}(U,h' \odo h) = \p(h,h')$ for all $h\in H$, $h'\in H'$.

\begin{defi}\label{def:new-ll-YD-gen}
Let $\alp : X \to H^{\op} \odo X$ a left coaction of $H^{\op}$ on $X$ and $\beta : X \to H' \odo X$ left coaction of $H'$ on $X$. We say that the tuple $(X,\alp,\beta)$ is a {\em left-left $(H,H')$-Yetter--Drinfeld algebra} if
\begin{equation}
(\id_{H'} \odo \alp)\beta = (\Sigma \odo \id_{X})(\ad({}^{\ops}U \odo \id_{X})(\id_{H^{\op}} \odo \beta)\alp.
\end{equation}
In other words, if the diagram
\[
\xymatrix@C=6pc@R=2pc{X\ar[r]^-{\beta}\ar[d]_-{\alp} & H' \odo X\ar[r]^-{id_{H'} \odo \alp} & H' \odo H^{\op} \odo X \\ H^{\op} \odo X \ar[r]_-{id_{H^{\op}} \odo \beta} & H^{\op} \odo H' \odo X\ar[r]_-{\ad({}^{\ops}U) \odo \id_{X}} & H^{\op} \odo H' \odo X\ar[u]_-{\Sigma \odo \id_{X}}}
\]
is commutative. We will denoted $(X,\alp,\beta) \in \aYD^{\,ll}(H,H')$ for say that $(X,\alp,\beta)$ is a left-left $(H,H')$-Yetter--Drinfeld algebra.
\end{defi}

Similarly,

\begin{defi}\label{def:new-pairing-rr-YD}
Let $\alp : X \to X \odo H^{\op}$ be a right coaction of $H^{\op}$ on $X$ and $\beta : X \to X \odo H'$ be right coaction of $H'$ on $X$. We say that the tuple $(X,\alp,\beta)$ is a {\em right-right $(H,H')$-Yetter--Drinfeld algebra} if
\begin{equation}
(\alp \odo \id_{H'})\beta = (\id_{X} \odo \ad({}^{\ops}U)^{-1})(\id_{X} \odo \Sigma)(\beta \odo \id_{H^{\op}})\alp.
\end{equation}
In other words, if the diagram
\[
\xymatrix@C=6pc@R=2pc{X\ar[r]^-{\beta}\ar[d]_-{\alp} & X \odo H'\ar[r]^-{\alp \odo \id_{H'}} & X \odo H^{\op} \odo H' \\ X \odo H^{\op}\ar[r]_-{\beta \odo \id_{H}} & X \odo H' \odo H^{\op}\ar[r]_-{id_{X} \odo \Sigma} & X \odo H^{\op} \odo H'\ar[u]_-{id_{X} \odo \ad({}^{\ops}U)^{-1}}}
\]
is commutative. We will denoted $(X,\alp,\beta) \in \aYD^{\,rr}(H,H')$ for say that $(X,\alp,\beta)$ is a right-right $(H,H')$-Yetter--Drinfeld algebra.
\end{defi}

\begin{rema}
Note that if we take the canonical pairing of a Hopf algebra $H$, $\p: H \times \du{H} \to \ku$, the above definitions are generalizations of the Definitions~\ref{def:new-ll-YD} and \ref{def:new-rr-YD}, 
\end{rema}

For a Hopf algebra $H$, we consider the pairing $\tilde{\p}: H^{\op} \times \du{H}^{\op} \to \ku$ defined by $\tilde{\p}(h^{\op},\ome^{\op}) = \p(h,\du{S}(\ome)) = \p(S(h),\ome)$. The canonical element for $\tilde{\p}$ is given by $\tilde{U} = ({}^{\op} \odo {}^{\op})(\id \odo \du{S}^{-1})(U)$, where $U$ is the canonical element of the canonical pairing $\p: H \times \du{H} \to \ku$.

\begin{prop}
Let $H$ be a finite-dimensional Hopf algebra. The functors
\[
\begin{array}{rccc}
\mathfrak{F}_{H} : & \aYD^{ll}(H,\du{H}) & \leadsto & \aYD^{rr}(\du{H},H) \\
& (X,\alp,\beta) & \mapsto & (X,{}_{\du{S}^{\ops}}\beta,{}_{S}\alp) \\
& f: (X,\alp,\beta) \to (X',\alp',\beta') & \mapsto & f: (X,{}_{\du{S}^{\ops}}\beta,{}_{S}\alp) \to (X,{}_{\du{S}^{\ops}}\beta',{}_{S}\alp')
\end{array}
\]
\[
\begin{array}{rccc}
\mathfrak{F}'_{H} : & \aYD^{rr}(H,\du{H}) & \leadsto & \aYD^{ll}(\du{H},H) \\
& (X,\alp,\beta) & \mapsto & (X,\beta_{\du{S}^{\ops}},\alp_{S}) \\
& f: (X,\alp,\beta) \to (X',\alp',\beta') & \mapsto & f: (X,\beta_{\du{S}^{\ops}},\alp_{S}) \to (X,\beta'_{\du{S}^{\ops}},\alp'_{S})
\end{array}
\]
are equivalence of categories. Moreover, $\mathfrak{F}'_{\du{H}}\mathfrak{F}_{H} = \mathfrak{Id}_{\aYD^{ll}(H,\du{H})}$ and $\mathfrak{F}_{\du{H}}\mathfrak{F}'_{H} = \mathfrak{Id}_{\aYD^{rr}(H,\du{H})}$.
\end{prop}
\pr
It is straightforward.
\fin

\begin{prop}
Let $H$ be a finite-dimensional Hopf algebra. The functors
\[
\begin{array}{rccc}
\mathfrak{F}_{H} : & \aYD^{ll}(H,\du{H}) & \leadsto & \aYD^{rr}(H^{\op},\du{H}^{\op}) \\
& (X,\alp,\beta) & \mapsto & (X,{}_{S}\alp,{}_{\du{S}^{\ops}}\beta) \\
& f: (X,\alp,\beta) \to (X',\alp',\beta') & \mapsto & f: (X,{}_{S}\alp',{}_{\du{S}^{\ops}}\beta') \to (X,{}_{S}\alp',{}_{\du{S}^{\ops}}\beta')
\end{array}
\]
\[
\begin{array}{rccc}
\mathfrak{F}'_{H} : & \aYD^{rr}(H,\du{H}) & \leadsto & \aYD^{ll}(H^{\op},\du{H}^{\op}) \\
& (X,\alp,\beta) & \mapsto & (X,\alp_{S},\beta_{\du{S}^{\ops}}) \\
& f: (X,\alp,\beta) \to (X',\alp',\beta') & \mapsto & f: (X,\alp_{S},\beta_{\du{S}^{\ops}}) \to (X,\alp'_{S},\beta'_{\du{S}^{\ops}})
\end{array}
\]
are equivalence of categories. Moreover, $\mathfrak{F}'_{H^{\op}}\mathfrak{F}_{H} = \mathfrak{Id}_{{}_{H}\aYD^{H}}$ and $\mathfrak{F}_{H^{\op}}\mathfrak{F}'_{H} = \mathfrak{Id}_{{}^{H}\aYD_{H}}$.
\end{prop}
\pr
It is straightforward.
\fin

We end this section, given an equivalence between the four versions of a Yetter--Drinfeld algebra in ``only coaction'' characterization.

\begin{prop}
Let $H$ be a finite-dimensional Hopf algebra. We have the following equivalence of categories
\begin{equation*}\label{eq:left_equivalence_all_yd_m}
\xymatrix@C=5pc@R=1.5pc{%
\shadowbox*{%
\begin{Bcenter}
$\;\aYD^{rl}(H,\du{H})\;$ \\[0.1cm]
$(X,\alp,\beta)$
\end{Bcenter}}\ar@{<~>}[r]\ar@{<~>}[d]%
&%
\shadowbox*{%
\begin{Bcenter}
$\;\aYD^{lr}(\du{H},H)\;$ \\[0.1cm]
$(X,\beta,\alp)$
\end{Bcenter}}\ar@{<~>}[d]%
\\
{\shadowbox*{%
\begin{Bcenter}
$\;\aYD^{ll}(H,\du{H})\;$ \\[0.1cm]
$(X,(S^{-1} \odo \id)\Sigma\alp,\beta)$
\end{Bcenter}}}\ar@{<~>}[r]%
&%
\shadowbox*{%
\begin{Bcenter}
$\;\aYD^{rr}(\du{H},H)\;$ \\[0.1cm]
$(X,(\id \odo \du{S}^{-1})\Sigma\beta,\alp)$
\end{Bcenter}}}
\end{equation*}
which is compatible with the equivalence given in Proposition~\ref{eq:left_equivalence_all_yd}. Moreover, the above equivalence respect the braided commutativity conditions.
\end{prop}
\pr
It is straightforward.
\fin

\section{Examples of Yetter--Drinfeld algebras}\label{sec:examples_yd}

\subsection{Trivial Yetter--Drinfeld algebras}\label{exa:trivial}

\begin{enumerate}[label=(\roman*)]
\item Let $H$ be a Hopf algebra and $X$ be an algebra. Consider the trivial left action $\rhd_{tr}: H \odo X \to X$, $h \odo x \mapsto \cou(h)x$ and the trivial left coaction $\alp_{tr}: X \to H \odo X$, $x \mapsto 1_{H} \odo x$. Then the tuple $(X,\rhd_{tr},\alp_{tr})$ is a left-left Yetter--Drinfeld algebra. Moreover, $(X,\rhd_{tr},\alp_{tr})$ is braided commutative if and only if $X$ is commutative. In particular, we always have $(\ku,\rhd_{tr},\alp_{tr}) \in ({}^{H}_{H}\aYD)^{bc}$. Then, we have
\begin{enumerate}[label=(\alph*)]
\item $(X,\alp_{tr},\beta_{tr})$ is a left-right $(H,\du{H})$-Yetter--Drinfeld algebra.
\item $(X,\alp_{tr},\beta_{tr})$ is braided commutative if and only if $X$ is commutative.
\end{enumerate}

\item Let $H$ be a commutative Hopf algebra and $X$ be an algebra. Consider the trivial left action $\rhd_{tr}: H \odo X \to X$, $h \odo x \mapsto \cou(h)x$ and any left coaction $\alp: X \to H \odo X$, then the tuple $(X,\rhd_{tr},\alp)$ is a left-left Yetter--Drinfeld algebra over $H$. Moreover, $(X,\rhd_{tr},\alp)$ is braided commutative if and only if $X$ is commutative. Then, for any coaction $\alp: X \to H \odo X$, we have
\begin{enumerate}[label=(\alph*)]
\item $(X,\alp,\beta_{tr})$ is a left-right $(H,\du{H})$-Yetter--Drinfeld algebra.
\item $(X,\alp,\beta_{tr})$ is braided commutative if and only if $X$ is commutative.
\end{enumerate}

\item Let $H$ be a cocommutative Hopf algebra and $X$ be an algebra. Consider any left action $\rhd: H \odo X \to X$ and the trivial left coaction $\alp_{tr}: X \to H \odo X$, $x \mapsto 1_{H} \odo x$, then the tuple $(X,\rhd,\alp_{tr})$ is a left-left Yetter--Drinfeld algebra over $H$. Moreover, $(X,\rhd,\alp_{tr})$ is braided commutative if and only if $X$ is commutative. Then, for any coaction $\beta: X \to X \odo \du{H}$, we have:
\begin{enumerate}[label=(\alph*)]
\item $(X,\alp_{tr},\beta)$ is a left-right $(H,\du{H})$-Yetter--Drinfeld algebra.
\item $(X,\alp_{tr},\beta)$ is braided commutative if and only if $X$ is commutative.
\end{enumerate}

\end{enumerate}

\subsection{Hopf algebras}\label{exa:hopf_algebras}

Let $H$ be a Hopf algebra. Consider the left and right adjoint actions of $H$ on itself, $\brhd_{ad}, \blhd_{ad}: H \odo H \to H$, and the comultiplication $\com_{H}: H \to H \odo H$ see as left or right coaction of $H$ on $H$. Denote by $\com_{H,S^{\ops}}: H \to H^{\op} \odo H$ and ${}_{S^{\ops}}\com_{H}: H \to H \odo H^{\op}$ the coactions $\com_{H,S^{\ops}}(h) = S^{-1}(h_{(2)})^{\op} \odo h_{(1)}$ and ${}_{S^{\ops}}\com_{H}(h) = h_{(2)} \odo S^{-1}(h_{(1)})^{\op}$, respectively. It holds
\begin{enumerate}[label=(\roman*)]
\item $(H,\brhd_{ad},\com_{H})$ is a left-left braided commutative Yetter--Drinfeld algebra over $H$.
\item $(H,\brhd_{ad},{}_{S^{\ops}}\com_{H})$ is a left-right braided commutative Yetter--Drinfeld algebra over $H$.
\item $(H,\blhd_{ad},\com_{H,S^{\ops}})$ is a right-left braided commutative Yetter--Drinfeld algebra over $H$.
\item $(H,\blhd_{ad},\com_{H})$ is a right-right braided commutative Yetter--Drinfeld algebra over $H$.
\end{enumerate}
On the other hand, because $\ad^{\du{H}}_{r} = \alp_{\brhd_{ad}} = \alp_{\Sigma(U^{-1})}$ and $\ad^{\du{H}}_{l} = \alp_{\blhd_{ad}} = \alp_{U}$, it holds
\begin{enumerate}[label=(\roman*)]
\item $(H,\com_{H},\ad^{\du{H}}_{r})$ is a left-right braided commutative $(H,\du{H})$-Yetter--Drinfeld algebra.
\item $(H,{}_{S^{\ops}}\com_{H},\ad^{\du{H}}_{r})$ is a right-right braided commutative $(H,\du{H})$-Yetter--Drinfeld algebra.
\item $(H,\com_{H,S^{\ops}},\ad^{\du{H}}_{l})$ is a left-left braided commutative $(H,\du{H})$-Yetter--Drinfeld algebra.
\item $(H,\com_{H},\ad^{\du{H}}_{l})$ is a right-left braided commutative $(H,\du{H})$-Yetter--Drinfeld algebra.
\end{enumerate}

\begin{rema}
Given a Hopf algebra $H$, the left and right adjoint coactions of $\du{H}$ on $H$ are given by
\[
\begin{array}{rccc}
\ad^{\du{H}}_{l} : & H & \to & \du{H} \odo H \\
& h & \mapsto & \Sigma(U^{-1})(1 \odo h)\Sigma(U)
\end{array},
\qquad
\begin{array}{rccc}
\ad^{\du{H}}_{r} : & H & \to & H \odo \du{H} \\
& h & \mapsto & U(h \odo 1)U^{-1}
\end{array},
\]
respectively. Here $U$ denotes the canonical element for the canonical pairing $\p$. The actions of $H$ associated to this coactions of $\du{H}$ are given by $\lhd_{\ad^{\du{H}}_{l}} = \blhd_{ad}$ and $\rhd_{\ad^{\du{H}}_{r}} = \brhd_{ad}$ respectively. Note that
\begin{equation}\label{eq:6}
(\com_{H} \odo \id_{\du{H}})\ad^{\du{H}}_{r} = U_{13}U_{23}\com_{H}(\cdot)_{12}U^{-1}_{23}U^{-1}_{13} = \ad(U_{13})(\id_{H} \odo \ad^{\du{H}}_{r})\com_{H}
\end{equation}
and
\begin{equation}\label{eq:7}
(\id_{\du{H}} \odo \com_{H})\ad^{\du{H}}_{l} = \Sigma(U)^{-1}_{13}\Sigma(U)^{-1}_{12}\com_{H}(\cdot)_{23}\Sigma(U)_{12}\Sigma(U)_{13} = \ad(\Sigma(U)^{-1}_{13})(\ad^{\du{H}}_{l} \odo \id_{H})\com_{H}.
\end{equation}
\end{rema}

\subsection{Coideal sub-algebras}\label{exa:coideals}

Let $H$, $H'$ be two Hopf algebras and $\pi: H \to H'$ be a surjective Hopf algebra homomorphism. The linear maps $\lambda_{l} = (\pi \odo \id_{H})\com_{H}$ and $\lambda_{r} = (\id_{H} \odo \pi)\com_{H}$ are a left coaction and right coaction of the Hopf algebra $H'$ on $H$ respectively. Moreover, we have the equalities
\begin{equation}\label{eq:3}
(\lambda_{l} \odo \id_{H})\com_{H} = (\id_{H'} \odo \com_{H})\lambda_{l}, \quad\quad (\id_{H} \odo \lambda_{r})\com_{H} = (\com_{H} \odo \id_{H'})\lambda_{r},
\end{equation}
\begin{equation}\label{eq:4}
(\id_{\du{H}} \odo \lambda_{l})\ad^{\du{H}}_{l} = \Sigma(U^{-1})_{13}(\id \odo \pi)(\Sigma(U^{-1}))_{12}\lambda_{l}(\cdot)_{23}(\id \odo \pi)(\Sigma(U))_{12}\Sigma(U)_{13}
\end{equation}
and
\begin{equation}\label{eq:5}
(\lambda_{r} \odo \id_{\du{H}})\ad^{\du{H}}_{r} = U_{13}(\pi \odo \id)(U)_{23}\lambda_{r}(\cdot)_{12}(\pi \odo \id)(U^{-1})_{23}U^{-1}_{13}.
\end{equation}
Consider now the sub-algebras of fixed points relate to the coaction $\lambda_{l}$ and $\lambda_{r}$ respectively:
\[
I_{l}(\pi) := \{ h \in H \;:\; \lambda_{l}(h) = 1_{H'} \odo h \}, \quad\quad I_{r}(\pi) := \{ h \in H \;:\; \lambda_{r}(h) = h \odo 1_{H'} \}.
\]
Using the equalities~\eqref{eq:3},~\eqref{eq:4} and~\eqref{eq:5}, we have
\[
\com_{H}(I_{r}(\pi)) \subseteq H \odo I_{r}(\pi), \qquad \ad^{\du{H}}_{r}(I_{r}(\pi)) \subseteq I_{r}(\pi) \odo \du{H},
\]
and
\[
\com_{H}(I_{l}(\pi)) \subseteq I_{l}(\pi) \odo H, \qquad \ad^{\du{H}}_{l}(I_{l}(\pi)) \subseteq \du{H} \odo I_{l}(\pi).
\]

\begin{prop}
With the above notations, $\alp_{I_{r}} := \com_{H}|_{I_{r}(\pi)}: I_{r}(\pi) \to H \odo I_{r}(\pi)$ is a left coaction of $H$ on $I_{r}(\pi)$ and $\beta_{I_{r}} := \ad^{\du{H}}_{r}|_{I_{r}(\pi)}: I_{r}(\pi) \to I_{r}(\pi) \odo \du{H}$ is a right coaction of $\du{H}$ on $I_{r}(\pi)$ such that $(I_{r}(\pi),\alp_{I_{r}},\beta_{I_{r}})$ is a braided commutative left-right $(H,\du{H})$-Yetter--Drinfeld algebra. Similarly, $\alp'_{I_{l}} := \com_{H}|_{I_{l}(\pi)}: I_{l}(\pi) \to I_{l}(\pi) \odo H$ is a right coaction of $H$ on $I_{l}(\pi)$ and $\beta'_{I_{l}} := \ad^{\du{H}}_{l}|_{I_{l}(\pi)}: I_{l}(\pi) \to \du{H} \odo I_{l}(\pi)$ is a left coaction of $\du{H}$ on $I_{l}(\pi)$ such that $(I_{l}(\pi),\alp'_{I_{l}},\beta'_{I_{l}})$ is a braided commutative right-left $(H,\du{H})$-Yetter--Drinfeld algebra.
\end{prop}
\pr
Take the equality~\eqref{eq:6}, then $(I_{r}(\pi),\alp_{I_{r}},\beta_{I_{r}})$ is left-right $(H,\du{H})$-Yetter--Drinfeld algebra. On the other hand, for any $x,y \in H$, we have
\[
(x_{(1)} \brhd_{ad} y)x_{(2)} = x_{(1)}yS(x_{(2)})x_{(3)} = \cou(x_{(2)})x_{(1)}y = xy,
\]
this says that $(H,\brhd_{ad},\com_{H})$ is a left-left braided commutative Yetter--Drinfeld algebra, then by Proposition~\ref{prop:ll-bc-equiv}, we have that
\begin{align*}
\m_{\He(H) \odo H}(\iota_{1,H} \odo \iota_{2,H})(\com_{H}(x) \odo \Sigma\ad^{\du{H}}_{r}(y)) = \m_{\He(H) \odo H}(\iota_{2,H} \odo \iota_{1,H})(\Sigma\ad^{\du{H}}_{r}(y) \odo \com_{H}(x))
\end{align*}
for any $x,y \in H$. In particular, because $I_{r}(\pi) \subseteq H$, we have
\begin{align*}
\m_{\He(H) \odo I_{r}(\pi)}(\iota_{1,I_{r}(\pi)} \odo \iota_{2,I_{r}(\pi)})(\alp_{I_{r}}(x) \odo \Sigma\beta_{I_{r}}(y)) = \m_{\He(H) \odo I_{r}(\pi)}(\iota_{2,I_{r}(\pi)} \odo \iota_{1,I_{r}(\pi)})(\Sigma\beta_{I_{r}}(y) \odo \alp_{I_{r}}(x))
\end{align*}
for all $x, y \in I_{r}(\pi)$. This implies that $(I_{r}(\pi),\alp_{I_{r}},\beta_{I_{r}})$ is braided commutative. Similarly, the right-left version can be proved.
\fin

\begin{rema}
If we consider the surjective Hopf algebra homomorphism $\cou: H \to \ku$. We have that $(I_{r}(\cou),\alp_{I_{r}},\beta_{I_{r}}) = (H,\com_{H},\ad^{\du{H}}_{r})$ as elements in $\aYD^{lr}_{bc}(H,\du{H})$ and $(I_{l}(\cou),\alp'_{I_{l}},\beta'_{I_{l}}) = (H,\com_{H},\ad^{\du{H}}_{l})$ as elements in $\aYD^{rl}_{bc}(H,\du{H})$.
\end{rema}

\subsection{Heisenberg algebras}\label{sec:heis-algebra-over}

Let $H$ be a Hopf algebra. Take the Heisenberg algebra $\He(\du{H}) := \du{H} \# H$ and the Drinfeld double $\D(H) := \D_{R}(H) = \du{H}^{\co} \bicroi H$ with Drinfeld codouble $\T(H) := \T_{R}(H) = H^{\op} \odo \du{H}$. Following~\cite{Se11}, if we consider the left action
\[
\begin{array}{rccc}
\rhd : & \D(H) \odo \He(\du{H}) & \to & \He(\du{H}) \\
& (\ome \bicroi h) \odo (\theta \# g) & \mapsto & \ome_{(3)}(h_{(1)} \brhd \theta)\du{S}^{-1}(\ome_{(2)}) \# (h_{(2)}gS(h_{(3)}) \blhd \du{S}^{-1}(\ome_{(1)}))
\end{array}
\]
and the left coaction
\[
\begin{array}{rccc}
\delta : & \He(\du{H}) & \to & \D(H) \odo \He(\du{H}) \\
& \ome \# h & \mapsto & \ome_{(2)} \bicroi h_{(1)} \odo \ome_{(1)} \# h_{(2)} 
\end{array}
\]
then $(\He(\du{H}),\rhd,\delta)$ is a left-left braided commutative Yetter--Drinfeld algebra over the Hopf algebra $\D(H)$.

In what follows, we will give a new way to construct this same Yetter--Drinfeld algebra structure using our ``only coaction'' characterization.

\subsubsection*{The left coaction of $\D(H)$ on $\He(\du{H})$} Consider the twisting
\[
\begin{array}{rccc}
\eta : & \D(H) \odo \D(H) & \to & \ku \\
& \theta \odo h \odo \theta' \odo h' & \mapsto & \bar{\p}(\theta,1)\p(h,\theta')\bar{\p}(1,h')
\end{array}
\]
induced by the canonical pairing $\p$ between $H$ and $\du{H}$. Following Theorem~\ref{theo:Double_Heisenberg}, consider the isomorphism induced by the twisting $\eta$, $i_{\eta} : \D(H) \to \D(H)_{\eta}$ and the Lu's isomorphism $L_{2} : \D(H)_{\eta} \to \He(\du{H})$. Using the notation $\bar{L}_{2} := L_{2} i_{\eta}: \D(H) \to \He(\du{H})$, the algebra homomorphism
\[
\begin{array}{lccc}
\Gamma' = (\id_{\D(H)} \odo \bar{L}_{2})\com_{\D(H)} \bar{L}^{-1}_{2} : & \He(\du{H}) & \to & \D(H) \odo \He(\du{H}) \\
& \theta \# h & \mapsto & \theta_{(2)} \bicroi h_{(1)} \odo \theta_{(1)} \# h_{(2)}
\end{array}
\]
is a left coaction of the Drinfeld algebra $\D(H)$ on $\He(\du{H})$. Indeed,
\begin{align*}
(\id_{\D(H)} \odo \Gamma')\Gamma' & = (\id_{\D(H)} \odo (\id_{\D(H)} \odo \bar{L}_{2})\com_{\D(H)} \bar{L}^{-1}_{2})(\id_{\D(H)} \odo \bar{L}_{2})\com_{\D(H)} \bar{L}^{-1}_{2} \\
& = (\id_{\D(H)} \odo \id_{\D(H)} \odo \bar{L}_{2})(\id_{\D(H)} \odo \com_{\D(H)})\com_{\D(H)} \bar{L}^{-1}_{2} \\
& = (\com_{\D(H)} \odo \id_{\He(\du{H})})(\id_{\D(H)} \odo \bar{L}_{2})\com_{\D(H)} \bar{L}^{-1}_{2} \\
& = (\com_{\D(H)} \odo \id_{\He(\du{H})})\Gamma',
\end{align*}
and
\begin{align*}
(\cou_{\D(H)} \odo \id_{\He(\du{H})})\Gamma' & = (\cou_{\D(H)} \odo \id_{\He(\du{H})})(\id_{\D(H)} \odo \bar{L}_{2})\com_{\D(H)} \bar{L}^{-1}_{2} \\
& = \bar{L}_{2}(\cou_{\D(H)} \odo \id_{\He(\du{H})})\com_{\D(H)} \bar{L}^{-1}_{2} = \id_{\He(\du{H})}.
\end{align*}

\begin{rema}
Recall $\m_{\D(H)_{\eta}}(i_{\eta} \odo i_{\eta}) = i_{\eta} \m_{\D(H)}$, $\m_{\He(\du{H})}(L_{2} \odo L_{2}) = L_{2} \m_{\D(H)_{\eta}}$, where $\m_{\D(H)_{\eta}}$ denotes the multiplication map of the twisting algebra $\D(H)_{\eta}$, $\m_{\D(H)}$ denotes the multiplication map of the algebra $\D(H)$ and $\m_{\He(\du{H})}$ denotes the multiplication map of the algebra $\He(\du{H})$.
\end{rema}

\subsubsection*{The right coaction of $\T(H)$ on $\He(\du{H})$}

Here, we will consider a canonical left coaction of $H^{\op}$ and a canonical coaction of $\du{H}$ which allow to construct a right coaction of $\T(H)$ on $\He(\du{H})$.

\begin{itemize}
\item Consider the linear map
\[
\begin{array}{lccc}
\alp : & \He(\du{H}) & \to & \He(\du{H}) \odo H^{\op} \\
& \theta \# h & \mapsto & (\du{\iota}_{1} \odo \id)(\Sigma({}^{\ops}U))(\theta \# h \odo 1^{\op})(\du{\iota}_{1} \odo \id)(\Sigma(({}^{\ops}U)^{-1}))
\end{array}
\]

\begin{itemize}
\item[($*$)] Take $\V := (\du{\iota}_{1} \odo \id)(\Sigma({}^{\ops}U)) \in \He(\du{H}) \odo H^{\op}$. Because
\begin{align*}
(\id \odo \com_{H})(\V) & = (\id \odo \com_{H})(\du{\iota}_{1} \odo \id)(\Sigma({}^{\ops}U)) \\
& = (\du{\iota}_{1} \odo \id \odo \id)(\id \odo \com_{H})(\Sigma({}^{\ops}U)) \\
& = (\du{\iota}_{1} \odo \id \odo \id)\Sigma_{12}\Sigma_{23}(\com_{H} \odo \id)({}^{\ops}U) \\
& = (\du{\iota}_{1} \odo \id \odo \id)\Sigma_{12}\Sigma_{23}({}^{\ops}U_{13}{}^{\ops}U_{23}) \\
& = (\du{\iota}_{1} \odo \id \odo \id)(\Sigma({}^{\ops}U)_{12}\Sigma({}^{\ops}U)_{13}) = \V_{12}\V_{13},
\end{align*}
then $\alp$ is a right coaction of $H^{\op}$  on the Heisenberg algebra $\He(\du{H})$.

\item[($*$)] For any $\ome, \theta \in \du{H}$ and $h \in H$, follows from the equality
\begin{align*}
\ome \rhd_{\alp} (\theta \# h) & = (\id \odo \tilde{\p}(\cdot,\ome))(\alp(\theta \# h)) \\
& = (\id \odo \tilde{\p}(\cdot,\ome))(\V(\theta \# h \odo 1^{\op})\V^{-1}) \\
& = (\id \odo \tilde{\p}(\cdot,\ome_{(3)}))(\V)(\id \odo \tilde{\p}(\cdot,\ome_{(2)}))(\theta \# h \odo 1^{\op})(\id \odo \tilde{\p}(\cdot,\ome_{(1)}))(\V^{-1}) \\
& = \du{\iota}_{1}((\tilde{\p}(\cdot,\ome_{(2)}) \odo \id)({}^{\ops}U))(\theta \# h)\du{\iota}_{1}((\tilde{\p}(\cdot,\ome_{(1)}) \odo \id)(({}^{\ops}U)^{-1})) \\
& = \du{\iota}_{1}(\ome_{(2)})(\theta \# h)\du{\iota}_{1}(\du{S}^{-1}(\ome_{(1)})).
\end{align*}
that, the left action induced by $\alp$ is given by
\[
\begin{array}{lccc}
\rhd_{\alp} : & \du{H}^{\co} \odo \He(\du{H}) & \to & \He(\du{H})  \\
& \ome \odo \theta \# h & \mapsto & \du{\iota}_{1}(\ome_{(2)})(\theta \# h)\du{\iota}_{1}(\du{S}^{-1}(\ome_{(1)}))
\end{array}.
\]
Moreover, we can see that
\begin{equation}\label{eq:alpha}
\ome \rhd_{\alp} (\theta \# h) = (\ome_{(2)} \# 1)(\theta \# h)(\du{S}^{-1}(\ome_{(1)}) \# 1) = \ome_{(3)}\theta\du{S}^{-1}(\ome_{(2)}) \# (h \blhd \du{S}^{-1}(\ome_{(1)})).
\end{equation}
\end{itemize}

\item Consider the linear map
\[
\begin{array}{lccc}
\beta : & \He(\du{H}) & \to & \He(\du{H}) \odo \du{H} \\
& \theta \# h & \mapsto & (\du{\iota}_{2} \odo \id)(U)(\theta \# h \odo 1)(\du{\iota}_{2} \odo \id)(U^{-1})
\end{array}
\]

\begin{itemize}
\item[($*$)] Denote $\U = (\du{\iota}_{2} \odo \id)(U) \in \He(\du{H}) \odo \du{H}$. Because
\begin{align*}
(\id \odo \com_{\du{H}})(\U) & = (\id \odo \com_{\du{H}})(\du{\iota}_{2} \odo \id)(U) = (\du{\iota}_{2} \odo \id \odo \id)(\id \odo \com_{\du{H}})(U) \\
& = (\du{\iota}_{2} \odo \id \odo \id)(U_{12}U_{13}) = \U_{12}\U_{13},
\end{align*}
then $\beta$ is a right coaction of $\du{H}$ on the Heisenberg algebra $\He(\du{H})$.

\item[($*$)] Because
\begin{align*}
g \rhd_{\beta} (\theta \# h) & = (\id \odo \p(g,\cdot))(\beta(\theta \# h)) \\
& = (\id \odo \p(g,\cdot))(\U(\theta \# h \odo 1)\U^{-1}) \\
& = (\id \odo \p(g_{(1)},\cdot))(\U)(\id \odo \p(g_{(2)},\cdot))(\theta \# h \odo 1)(\id \odo \p(g_{(3)},\cdot))(\U^{-1}) \\
& = \du{\iota}_{2}((\id \odo \p(g_{(1)},\cdot))(U))(\theta \# h)\du{\iota}_{2}((\id \odo \p(g_{(2)},\cdot))(U^{-1})) \\
& = \du{\iota}_{2}(g_{(1)})(\theta \# h)\du{\iota}_{2}(\du{S}(g_{(2)}))
\end{align*}
for any $\theta \in \du{H}$ and $g, h \in H$. Then, the induced left action by $\beta$ is given by
\[
\begin{array}{lccc}
\rhd_{\beta} : & H \odo \He(\du{H}) & \to & \He(\du{H})  \\
& g \odo \theta \# h & \mapsto & \du{\iota}_{2}(g_{(1)})(\theta \# h)\du{\iota}_{2}(S(g_{(2)}))
\end{array}.
\]
Moreover, we can see that
\begin{equation}\label{eq:beta}
g \rhd_{\beta} (\theta \# h) = (1 \# g_{(1)})(\theta \# h)(1 \# S(g_{(2)})) = (g_{(1)} \brhd \theta) \# g_{(2)}hS(g_{(3)}).
\end{equation}
\end{itemize}
\end{itemize}

\begin{lemm}\label{lemm:hola}
Let $U_{1} \odo U_{2} = U \in H \odo \du{H}$. Consider
\[
\U = (\du{\iota}_{2} \odo \id_{\du{H}})(U) \in \He(\du{H}) \odo \du{H}, \qquad \V = (\du{\iota}_{1} \odo \id_{H^{\op}})(\Sigma({}^{\ops}U)) \in \He(\du{H}) \odo H^{\op}.
\]
Then
\[
\V_{12}\U_{13} = (\bar{L}_{2} \odo \id_{\T(H)})(\W) = U_{2} \# U_{1} \odo U^{\op}_{1} \odo U_{2}
\]
where $\W := U_{2} \bicroi U_{1} \odo U^{\op}_{1} \odo U_{2} \in \D(H) \odo \T(H)$ is the canonical element associated to the canonical pairing $\p_{\D}: \D(H) \times \T(H) \to \ku$.
\end{lemm}
\pr
The canonical element associated to the pairing $\p_{\D} = \tilde{\p} \odo \p$ is given by the element
\begin{align*}
\W & := (\iota_{\D(H)} \odo \id_{H^{\op}})(\Sigma({}^{\ops}U))_{12}(\iota_{\D(H)} \odo \id_{\du{H}})(U)_{13} \\
& = (\iota_{\D(H)} \odo \id_{H^{\op}})(U_{2} \odo U^{\op}_{1})_{12}(\iota_{\D(H)} \odo \id_{\du{H}})(U_{1} \odo U_{2})_{13} \\
& = U_{2} \bicroi U_{1} \odo U^{\op}_{1} \odo U_{2}.
\end{align*}
Thus, we have
\begin{align*}
\V_{12}\U_{13} & = (\du{\iota}_{1}\du{\iota}_{2} \odo \id_{H^{\op}} \odo \id_{\du{H}})(\Sigma({}^{\circ}U)_{13}U_{24}) = (\bar{L}_{2}\iota_{\D(H)} \odo \id_{H^{\op}} \odo \id_{\du{H}})(\Sigma({}^{\circ}U)_{13}U_{24}) \\
& = (\bar{L}_{2} \odo \id_{H^{\op}} \odo \id_{\du{H}})(\iota_{\D(H)} \odo \id_{H^{\op}} \odo \id_{\du{H}})(\Sigma({}^{\circ}U)_{13}U_{24}) \\
& = (\bar{L}_{2} \odo \id_{H^{\op}} \odo \id_{\du{H}})((\iota_{\D(H)} \odo \id_{H^{\op}})(\Sigma({}^{\circ}U))_{12}(\iota_{\D(H)} \odo \id_{\du{H}})(U)_{13}) \\
& = (\bar{L}_{2} \odo \id_{\T(H)})(\W) = (\bar{L}_{2} \odo \id_{\T(H)})(U_{2} \bicroi U_{1} \odo U^{\op}_{1} \odo U_{2}) \\
& = U_{2} \# U_{1} \odo U^{\op}_{1} \odo U_{2}.
\end{align*}
\fin

\begin{prop}\label{prop:right_YD_hola}
Using the coactions $\alp$ and $\beta$ as defined above. The tuple $(\He(\du{H}),\alp,\beta)$ is a right-right $(H,\du{H})$-Yetter--Drinfeld algebra.
\end{prop}
\pr
Using $U \in H \odo \du{H}$, consider the invertible elements $\U := (\du{\iota}_{2} \odo \id)(U) \in \He(\du{H}) \odo \du{H}$ and $\V := (\du{\iota}_{1} \odo \id)(\Sigma({}^{\ops}U)) \in \He(\du{H}) \odo H^{\op}$. Follows from Lemma~\ref{lemm:hola} that
\begin{align*}
(\id_{\He(\du{H})} \odo \com_{\T(H)})(\V_{12}\U_{13}) & = (\id_{\He(\du{H})} \odo \com_{\T(H)})(\bar{L}_{2} \odo \id_{\T(H)})(\W) \\
& = (\bar{L}_{2} \odo \id_{\T(H)} \odo \id_{\T(H)})(\id_{\D(H)} \odo \com_{\T(H)})(\W) \\
& = (\bar{L}_{2} \odo \id_{\T(H)} \odo \id_{\T(H)})(\W_{12}\W_{13}) \\
& = (\bar{L}_{2} \odo \id_{\T(H)} \odo \id_{\T(H)})(\W)_{12}(\bar{L}_{2} \odo \id_{\T(H)} \odo \id_{\T(H)})(\W)_{13} \\
& = (\V_{12}\U_{13})_{12}(\V_{12}\U_{13})_{13},
\end{align*}
and because $\U$ and $\V$ are invertible elements, we have that $(\alp \odo \id_{\du{H}})\beta: \He(\du{H}) \to \He(\du{H}) \odo \T(H)$, $x \mapsto \ad(\V_{12}\U_{13})(x \odo 1^{\op}_{H} \odo 1_{\du{H}})$ is a right coaction $\T(H)$ on $\He(\du{H})$. This implies that $(\He(\du{H}),\alp,\beta)$ is a right-right $(H,\du{H})$-Yetter--Drinfeld algebra by Proposition~\ref{prop:right_equivalence_radford}.
\fin

\begin{coro}\label{coro:hola}
There is a right coaction of $\T(H)$ on $\He(\du{H})$ given by
\[
\begin{array}{lccc}
\Gamma : & \He(\du{H}) & \to & \He(\du{H}) \odo \T(H) \\
& \theta \# g & \mapsto & (\alp \odo \id)(\beta(\theta \# g))
\end{array}
\]
which induced a left action of $\D(H)$ on $\He(\du{H})$ given by
\[
\begin{array}{lccc}
\rhd_{\Gamma}  : & \D(H) \odo \He(\du{H}) & \to & \He(\du{H})  \\
& \ome \bicroi h \odo \theta \# g & \mapsto & \ome_{(3)}(h_{(1)} \brhd \theta)\du{S}^{-1}(\ome_{(2)}) \# (h_{(2)}gS(h_{(3)}) \blhd \du{S}^{-1}(\ome_{(1)}))
\end{array}.
\]
\end{coro}
\pr
Follows directly by Proposition~\ref{prop:right_equivalence_radford} and Proposition~\ref{prop:right_YD_hola}. For the induced action, using the equalities~\eqref{eq:alpha} and~\eqref{eq:beta}, we have
\begin{align*}
(\ome \bicroi h) \rhd_{\Gamma} (\theta \# g) & = \ome \rhd_{\alp} (h \rhd_{\beta} (\theta \# g)) = \ome \rhd_{\alp} ((h_{(1)} \brhd \theta) \# h_{(2)}gS(h_{(3)})) \\
& = \ome_{(3)}(h_{(1)} \brhd \theta)\du{S}^{-1}(\ome_{(2)}) \# (h_{(2)}gS(h_{(3)}) \blhd \du{S}^{-1}(\ome_{(1)})),
\end{align*}
for all $h,g \in H$ and $\ome, \theta \in \du{H}$.
\fin


We finish given a new prove to show that the Heisenberg algebra $\He(\du{H})$ is a braided commutative Yetter--Drinfeld algebra over $\D(H)$.

\begin{theo}
Consider the left coaction of $\D(H)$ on $\He(\du{H})$,
\[
\begin{array}{lccc}
\Gamma' = (\id_{\D(H)} \odo \bar{L}_{2})\com_{\D(H)} \bar{L}^{-1}_{2} : & \He(\du{H}) & \to & \D(H) \odo \He(\du{H}) \\
& \theta \# a & \mapsto & \theta_{(2)} \bicroi h_{(1)} \odo \theta_{(1)} \# h_{(2)}
\end{array}
\]
and the right coaction of $\T(H)$ on $\He(\du{H})$ construct in Corollary~\ref{coro:hola},
\[
\begin{array}{lccc}
\Gamma : & \He(\du{H}) & \to & \He(\du{H}) \odo \T(H) \\
& \theta \# g & \mapsto & \V_{12}\U_{13}((\theta \# g) \odo 1_{\T(H)})\U^{-1}_{13}\V^{-1}_{12}
\end{array}
\]
where $\U = (\du{\iota}_{2} \odo \id_{\du{H}})(U) \in \He(\du{H}) \odo \du{H}$, $\V = (\du{\iota}_{1} \odo \id_{H^{\op}})(\Sigma({}^{\ops}U)) \in \He(\du{H}) \odo H^{\op}$ and $U = U_{1} \odo U_{2} \in H \odo \du{H}$ is the canonical element of $H$. Then, we have
\begin{enumerate}[label=\textup{(\roman*)}]
\item The tuple $(\He(\du{H}),\Gamma',\Gamma)$ is an element in the category $\aYD^{lr}(\D(H),\T(H))$.
\item The tuple $(\He(\du{H}),\rhd_{\Gamma},\Gamma')$ is an element in the category ${}^{\D(H)}_{\D(H)}\aYD$.
\end{enumerate}
\end{theo}
\pr
Consider the map $\Theta = (\Sigma\Gamma' \odo \id)\Gamma: \He(\du{H}) \to \He(\du{H}) \odo \D(H)^{\co} \odo \T(H)$. Take now the Majid-Drinfeld codouble of $\D(H)$ associated to the skew-copairing $\Sigma(\W)^{-1} \in \T(H) \odo \D(H)^{\co}$, this Hopf algebra will be denoted by $\T_{M}(\D(H))$. Its comultiplication is given by the map
\[
\com_{\W} = (\id_{\D(H)} \odo \Sigma\ad(\W)^{-1} \odo \id_{\T(H)})(\com^{\co}_{\D(H)} \odo \com_{\T(H)}).
\]
By Proposition~\ref{prop:right_equivalence_majid}, we have: $\Theta$ is a right coaction of $\T_{M}(\D(H))$ on $\He(\du{H})$ if and only if $(\He(\du{H}),\Gamma',\Gamma)$ is a left-right $(\D(H),\T(H))$-Yetter--Drinfeld algebra if and only if $(\He(\du{H}),\rhd_{\Gamma},\Gamma')$ is a left-left Yetter--Drinfeld algebra over $\D(H)$.

Then it is enough to prove that $\Theta$ is a right coaction of $\T_{M}(\D(H))$ on $\He(\du{H})$. For that, take $X = (\bar{L}_{2} \odo \id)(\W) = \V_{12}\U_{13} \in \He(\du{H}) \odo \T(H)$. We have $\Gamma(x) = \ad(X)(x \odo 1)$ for all $x \in \He(\du{H})$, and the equalities
\small
\begin{align*}
(\Sigma\Gamma' \odo \id)(X) & = (\Sigma \odo \id)((\id \odo \bar{L}_{2})\com_{\D(H)}\bar{L}^{-1}_{2} \odo \id)(\V_{12}\U_{13}) \\
& = (\Sigma \odo \id)((\id \odo \bar{L}_{2})\com_{\D(H)}\bar{L}^{-1}_{2} \odo \id)(\bar{L}_{2} \odo \id)(\W) \\
& = (\Sigma \odo \id)(\id \odo \bar{L}_{2} \odo \id)(\com_{\D(H)} \odo \id)(\W) \\
& = (\Sigma \odo \id)(\id \odo \bar{L}_{2} \odo \id)(\W_{13}\W_{23}) = \W_{23}X_{13}, \\
& \\
(\id \odo \com_{\W})(\W_{23}) & = 1 \odo \com_{\W}(\W) = 1 \odo (\id \odo \Sigma\ad(\W^{-1}) \odo \id)(\com^{\co}_{\D(H)} \odo \com)(\W) \\
& = 1 \odo (\id \odo \Sigma\ad(\W^{-1}) \odo \id)(\com^{\co}_{\D(H)} \odo \id \odo \id)(\W_{12}\W_{13}) \\
& = 1 \odo (\id \odo \Sigma\ad(\W^{-1}) \odo \id)((\com^{\co}_{\D(H)} \odo \id)(\W)_{123}(\com^{\co}_{\D(H)} \odo \id)(\W)_{124}) \\
& = 1 \odo (\id \odo \Sigma\ad(\W^{-1}) \odo \id)((\W_{23}\W_{13})_{123}(\W_{23}\W_{13})_{124}) \\
& = 1 \odo \Sigma_{23}(\W^{-1}_{23}\W_{23}\W_{13}\W_{24}\W_{14}\W_{23}) = 1 \odo \W_{12}\W_{34}\W_{14}\W_{32} \\
& = \W_{23}\W_{45}\W_{25}\W_{43}
\end{align*}
and
\begin{align*}
(\id \odo \com_{\W})(X_{13}) & = (\id \odo \id \odo \Sigma\ad(\W^{-1}) \odo \id)(\id \odo \com^{\co}_{\D(H)} \odo \com)(X_{13}) \\
& = (\id \odo \id \odo \Sigma\ad(\W^{-1}) \odo \id)((\id \odo \com)(X)_{145}) \\ 
& = (\id \odo \id \odo \Sigma\ad(\W^{-1}) \odo \id)((\id \odo \com)((\bar{L}_{2} \odo \id)(\W))_{145}) \\
& = (\id \odo \id \odo \Sigma\ad(\W^{-1}) \odo \id)(((\bar{L}_{2} \odo \id)(\W)_{12}(\bar{L}_{2} \odo \id)(\W)_{13})_{145}) \\
& = (\id \odo \id \odo \Sigma\ad(\W^{-1}) \odo \id)((\bar{L}_{2} \odo \id)(\W)_{14}(\bar{L}_{2} \odo \id)(\W)_{15}) \\
& = (\bar{L}_{2} \odo \id \odo \id \odo \id \odo \id)\Sigma_{34}(\W^{-1}_{34}\W_{14}\W_{15}\W_{34}) \\
& = (\bar{L}_{2} \odo \id \odo \id \odo \id \odo \id)(\W^{-1}_{43}\W_{13}\W_{15}\W_{43}) = \W^{-1}_{43}X_{13}X_{15}\W_{43}.
\end{align*}

Fix $x \in \He(\p)$. We have
\small
\begin{align*}
(\id \odo \com_{\W})\Theta(x) & = (\id \odo \com_{\W})(\Sigma\Gamma' \odo \id)\Gamma(x) \\
& = (\id \odo \com_{\W})(\Sigma \odo \id)((\id\odo \bar{L}_{2})\com_{\D(H)}\bar{L}^{-1}_{2} \odo \id)\ad(\V_{12}\U_{13})(x \odo 1_{\T(\bar{\p})}) \\
& = (\id \odo \com_{\W})(\W_{23}X_{13}(\Sigma\Gamma'(x) \odo 1_{\T(\bar{\p})})(\W_{23}X_{13})^{-1}) \\
& = (\id \odo \com_{\W})(\W_{23}X_{13})(\id \odo \com_{\W})(\Sigma\Gamma'(x) \odo 1)(\id \odo \com_{\W})(\W_{23}X_{13})^{-1} \\ 
& = \ad(\W_{23}\W_{45}\W_{25}X_{13}X_{15}\W_{43})((\id \odo \com_{\W})(\Sigma\Gamma'(x) \odo 1))
\end{align*}
\normalsize
and
\small
\begin{align*}
(\Theta \odo \id)\Theta(x) & = (\Theta \odo \id)(\W_{23}X_{13}(\Sigma\Gamma'(x) \odo 1)(\W_{23}X_{13})^{-1}) \\
& = (\W_{23}X_{13})_{123}((\Sigma\Gamma' \odo \id \odo \id)(\W_{23}X_{13}(\Sigma\Gamma'(x) \odo 1)(\W_{23}X_{13})^{-1}))_{1245}(\W_{23}X_{13})^{-1}_{123} \\
& = \W_{23}X_{13}\W_{45}((\Sigma\Gamma' \odo \id \odo \id)(X_{13}(\Sigma\Gamma'(x) \odo 1)X^{-1}_{13}))_{1245}\W^{-1}_{45}X^{-1}_{13}\W^{-1}_{23} \\
& = \W_{23}X_{13}\W_{45}(\Sigma\Gamma' \odo \id)(X)_{125}((\Sigma\Gamma' \odo \id)\Sigma\Gamma'(x) \odo 1)_{1245}(\Sigma\Gamma' \odo \id)(X)^{-1}_{125}\W^{-1}_{45}X^{-1}_{13}\W^{-1}_{23} \\
& = \W_{23}X_{13}\W_{45}(\Sigma\Gamma' \odo \id)(X)_{125}((\Sigma\Gamma' \odo \id)\Sigma\Gamma'(x))_{124}(\Sigma\Gamma' \odo \id)(X)^{-1}_{125}\W^{-1}_{45}X^{-1}_{13}\W^{-1}_{23} \\
& = \W_{23}X_{13}\W_{45}(\Sigma\Gamma' \odo \id)(X)_{125}((\id \odo \com^{\co}_{\D(H)})\Sigma\Gamma'(x))_{124}(\Sigma\Gamma' \odo \id)(X)^{-1}_{125}\W^{-1}_{45}X^{-1}_{13}\W^{-1}_{23} \\
& = \W_{23}X_{13}\W_{45}\W_{25}X_{15}((\id \odo \com^{\co}_{\D(H)})\Sigma\Gamma'(x))_{124}X^{-1}_{15}\W^{-1}_{25}\W^{-1}_{45}X^{-1}_{13}\W^{-1}_{23} \\
& = \ad(\W_{23}X_{13}\W_{45}\W_{25}X_{15})(((\id \odo \com^{\co}_{\D(H)})\Sigma\Gamma'(x))_{124}) \\
& = \ad(\W_{23}X_{13}\W_{45}\W_{25}X_{15})(((\id \odo \com^{\co}_{\D(H)})\Sigma\Gamma'(x))_{124}).
\end{align*}
\normalsize

Because,
\begin{align*}
(\id \odo \com^{\co}_{\D(H)})\Sigma\Gamma' & = (\id \odo \com^{\co}_{\D(H)})\Sigma(\id \odo \bar{L}_{2})\com_{\D(H)}\bar{L}^{-1}_{2} \\
& = (\id \odo \com^{\co}_{\D(H)})(\bar{L}_{2} \odo \id)\com^{\co}_{\D(H)}\bar{L}^{-1}_{2} \\
& = (\bar{L}_{2} \odo \id \odo \id)(\id \odo \com^{\co}_{\D(H)})\com^{\co}_{\D(H)}\bar{L}^{-1}_{2} \\
& = (\bar{L}_{2} \odo \id \odo \id)(\com^{\co}_{\D(H)} \odo \id)\com^{\co}_{\D(H)}\bar{L}^{-1}_{2},
\end{align*}
then
\small
\begin{align*} (\id \odo \com_{\W})(\Sigma\Gamma'(x) \odo 1) & = (\id \odo \id \odo \Sigma\ad(\W^{-1}) \odo \id)(\id \odo \com^{\co}_{\D(H)} \odo \com)(\Sigma\Gamma'(x) \odo 1) \\
& = (\id \odo \id \odo \Sigma\ad(\W^{-1}) \odo \id)(\id \odo \com^{\co}_{\D(H)} \odo \com)((\bar{L}_{2} \odo \id)\com^{\co}_{\D(H)}\bar{L}^{-1}_{2} \odo \id)(x \odo 1) \\
& = (\bar{L}_{2} \odo \id \odo \Sigma\ad(\W^{-1}) \odo \id)(\com^{\co}_{\D(H)} \odo \id \odo \com)(\com^{\co}_{\D(H)}\bar{L}^{-1}_{2} \odo \id)(x \odo 1) \\
& = (\bar{L}_{2} \odo \id \odo \id \odo \id^{\odo 2})(\com^{\co}_{\D(H)} \odo \id \odo \id^{\odo 2})\Sigma_{23}\ad(\W^{-1})_{23}(\com^{\co}_{\D(H)} \odo \id^{\odo 2})(\bar{L}^{-1}_{2}(x) \odo 1^{\odo 2}) \\
& = \Sigma_{34}(\bar{L}_{2} \odo \id \odo \id \odo \id^{\odo 2})(\com^{\co}_{\D(H)} \odo \id \odo \id^{\odo 2})(\com^{\co}_{\D(H)} \odo \id^{\odo 2})(\bar{L}^{-1}_{2}(x) \odo 1^{\odo 2}) \\
& = \Sigma_{34}((\id \odo \com^{\co}_{\D(H)})\Sigma\Gamma'(x) \odo 1 \odo 1) = ((\id \odo \com^{\co}_{\D(H)})\Sigma\Gamma'(x))_{124},
\end{align*}
\normalsize
thus $(\Theta \odo \id)\Theta(x) = (\Theta \odo \id)\Theta(x)$. The last equality and the injectivity of $\Theta$ imply that $\Theta$ is a right coaction of $\T_{M}(\D(H))$ on $\He(\du{H})$.

\td{Check why $(1 \odo \W^{-1})(\com^{\co}_{\D(H)} \odo \id)(1 \odo \W) = (\com^{\co}_{\D(H)} \odo \id)$}

\td{The braided commutativity condition comes from the standard characterization}

\fin


\bibliographystyle{abbrv}
\bibliography{biblio}

\addresseshere

\newpage

\appendix

\section{The two characterizations of Yetter--Drinfeld objects}\label{sec:appendix_equivalences}

The equivalence between the standard and the only coaction characterization of Yetter--Drinfeld modules can be summarized in the following table:

\small
\begin{center}
\begin{tabular}{|c|c|c||c|c|c|}
\hline
\multicolumn{3}{|c||}{\textsc{Standard}} & \multicolumn{3}{|c|}{\textsc{Only coaction}} \\
\hline
\hline
category name & objects & equivalence & category name & objects & equivalence \\
\hline
&\multirow{3}{8em}{\centering $X$ \\ $\rhd : H \odo X \to X$ $\alp: X \to H \odo X$}&&&\multirow{3}{8em}{\centering $X$ \\ $\alp : X \to H \odo X$ $\beta: X \to X \odo \du{H}$}& \\ ${}^{H}_{H}\YD$&&${}_{\D_{M}(H)}\text{Mod}$&$\YD^{\,lr}(H,\du{H})$&&$\text{Mod}^{\T_{M}(H)}$ \\ &&&&& \\
\hline
&\multirow{3}{8em}{\centering $X$ \\ $\lhd : X \odo H \to X$ $\alp: X \to X \odo H$}&&&\multirow{3}{8em}{\centering $X$ \\ $\alp : X \to X \odo H$ $\beta: X \to \du{H} \odo X$}& \\ $\YD{}^{H}_{H}$&&$\text{Mod}_{\D_{M}(H)}$&$\YD^{\,rl}(H,\du{H})$&&${}^{\T_{M}(H)}\text{Mod}$ \\ &&&&& \\
\hline
&\multirow{3}{8em}{\centering $X$ \\ $\lhd : X \odo H \to X$ $\alp: X \to H \odo X$}&&&\multirow{3}{8em}{\centering $X$ \\ $\alp : X \to H \odo X$ $\beta: X \to \du{H} \odo X$}& \\ ${}^{H}\YD_{H}$&&$\text{Mod}_{\D_{R}(H)}$&$\YD^{\,ll}(H,\du{H})$&&${}^{\T_{R}(H)}\text{Mod}$ \\ &&&&& \\
\hline
&\multirow{3}{8em}{\centering $X$ \\ $\rhd : H \odo X \to X$ $\alp: X \to X \odo H$}&&&\multirow{3}{8em}{\centering $X$ \\ $\alp : X \to X \odo H$ $\beta: X \to X \odo \du{H}$}& \\ ${}_{H}\YD^{H}$&&${}_{\D_{R}(H)}\text{Mod}$&$\YD^{\,rr}(H,\du{H})$&&$\text{Mod}^{\T_{R}(H)}$ \\ &&&&& \\
\hline
\end{tabular}
\end{center}
\normalsize

\medskip

In the special case of algebraic objects in a Yetter--Drinfeld category or the so-called Yetter--Drinfeld algebras, we have:

\small
\begin{center}
\begin{tabular}{|c|c|c||c|c|c|}
\hline
\multicolumn{3}{|c||}{\textsc{Standard}} & \multicolumn{3}{|c|}{\textsc{Only coaction}} \\
\hline
\hline
category name & algebra objects & equivalence & category name & algebra objects & equivalence \\
\hline
&\multirow{3}{8em}{\centering $X$ \\ $\rhd : H \odo X \to X$ $\alp: X \to H \odo X$}&&&\multirow{3}{8em}{\centering $X$ \\ $\alp : X \to H \odo X$ $\beta: X \to X \odo \du{H}$}& \\ ${}^{H}_{H}\aYD$&&${}_{\D_{M}(H)}\alg$&$\aYD^{\,lr}(H,\du{H})$&&$\alg^{\T_{M}(H)}$ \\ &&&&& \\
\hline
&\multirow{3}{8em}{\centering $X$ \\ $\lhd : X \odo H \to X$ $\alp: X \to X \odo H$}&&&\multirow{3}{8em}{\centering $X$ \\ $\alp : X \to X \odo H$ $\beta: X \to \du{H} \odo X$}& \\ $\aYD{}^{H}_{H}$&&$\alg_{\D_{M}(H)}$&$\aYD^{\,rl}(H,\du{H})$&&${}^{\T_{M}(H)}\alg$ \\ &&&&& \\
\hline
&\multirow{3}{8em}{\centering $X$ \\ $\lhd : X \odo H \to X$ $\alp: X \to H^{\op} \odo X$}&&&\multirow{3}{8em}{\centering $X$ \\ $\alp : X \to H^{\op} \odo X$ $\beta: X \to \du{H} \odo X$}& \\ ${}^{H}\aYD_{H}$&&$\alg_{\D_{R}(H)}$&$\aYD^{\,ll}(H,\du{H})$&&${}^{\T_{R}(H)}\alg$ \\ &&&&& \\
\hline
&\multirow{3}{8em}{\centering $X$ \\ $\rhd : H \odo X \to X$ $\alp: X \to X \odo H^{\op}$}&&&\multirow{3}{8em}{\centering $X$ \\ $\alp : X \to X \odo H^{\op}$ $\beta: X \to X \odo \du{H}$}& \\ ${}_{H}\aYD^{H}$&&${}_{\D_{R}(H)}\alg$&$\aYD^{\,rr}(H,\du{H})$&&$\alg^{\T_{R}(H)}$ \\ &&&&& \\
\hline
\end{tabular}
\end{center}
\normalsize

\subsection*{Braided commutativity of Yetter--Drinfeld algebras}

In this section, we summarize the braided commutativity conditions for the standard characterization of Yetter--Drinfeld algebras.

\small
\begin{center}
\begin{tabular}{|c|c||c|c|c|}
\hline
\multicolumn{2}{|c||}{\textsc{Yetter--Drinfeld categories}} & \multicolumn{2}{c|}{\textsc{braided commutativity condition}} \\
\hline
\hline
category & objects & using the braiding map $\tau$ & using the braiding map $\rho$ \\
\hline
& \multirow{3}{10em}{\centering $(X,\rhd,\alpha)$ \\ $\alp: x \mapsto x_{[-1]} \odo x_{[0]}$} && \\ ${}^{H}_{H}\aYD$ && $(x_{[-1]} \rhd y)x_{[0]} = xy$ & $y_{[0]}(S^{-1}(y_{[-1]}) \rhd x) = xy$ \\ &&& \\
\hline
& \multirow{3}{10em}{\centering $(X,\lhd,\alpha)$ \\ $\alp: x \mapsto x_{[0]} \odo x_{[1]}$} && \\ $\aYD^{H}_{H}$ && $(y \lhd S^{-1}(x_{[1]}))x_{[0]} = xy$ & $y_{[0]}(x \lhd y_{[1]}) = xy$ \\ &&& \\
\hline
& \multirow{3}{10em}{\centering $(X,\lhd,\alpha)$ \\ $\alp: x \mapsto x^{\op}_{[-1]} \odo x_{[0]}$} && \\ ${}^{H}\aYD_{H}$ && $(y \lhd x_{[-1]})x_{[0]} = xy$ & $y_{[0]}(x \lhd S(y_{[-1]})) = xy$ \\ &&& \\
\hline
& \multirow{3}{10em}{\centering $(X,\rhd,\alpha)$ \\ $\alp: x \mapsto x_{[0]} \odo x^{\op}_{[1]}$} && \\ ${}_{H}\aYD^{H}$ && $(S(x_{[1]}) \rhd y)x_{[0]} = xy$ & $y_{[0]}(y_{[1]} \rhd x) = xy$ \\ &&& \\
\hline
\end{tabular}
\end{center}
\normalsize

\section{Comparison of the Drinfeld double constructions}\label{sec:appendix_comparison}

In the literature, there are many equivalent constructions for the Drinfeld double. In this section, we summarize the important characterizations used in this work and we show how they are related.

\begin{defi}
Let $H$ be a Hopf algebra. The Drinfeld double can be define in the following ways:
\begin{itemize}
\item[*] (Drinfeld's original version) $\D(H) = H \bicroi \du{H}^{\co}$
\item[*] (Majid's characterization) $\D_{M}(H) = \du{H}^{\op} \bicroi H$
\item[*] (Radford's characterization) $\D_{R}(H) = \du{H}^{\co} \bicroi H$
\item[*] (Taipe's characterization) $\D_{T}(H) = \du{H} \bicroi H^{\op}$
\item[*] (Delvaux-Van Daele's characterization) $\D_{DVD}(H) = \du{H} \bicroi H^{\co}$
\end{itemize}
\end{defi}

\begin{rema}
The Drinfeld double construction used in Lu's works is the Radford's characterization.
\end{rema}

\begin{prop}
Let $H$ be a Hopf algebra. We have the following relations:
\begin{enumerate}[label=\textup{(\roman*)}]
\item $\D_{T}(H) = \D_{R}(H^{\op}) \iso \D_{M}(H^{\op})$
\item $\D_{T}(H) \iso \D_{DVD}(H) = \D(\du{H}) = \D_{M}(H^{\co})$
\end{enumerate}
\end{prop}
\pr
It is straightforward.
\fin

\section{Back to the braided commutativity condition}\label{sec:appendix_bc}

In this section, we prove that the braided commutativity condition in the standard characterization of Yetter--Drinfeld algebras over a Hopf algebra $H$ depend on a canonical action of the Drinfeld double $\D_{R}(H)$ arising from the Yetter--Drinfeld structure and the canonical element $U$ of $H$. Recall the canonical isomorphism
\begin{equation}\label{eq:canonical_DMR}\tag{CI}
\begin{array}{lccc}
\natural :& \D_{M}(H) & \iso & \D_{R}(H) \\
& \ome^{\op} \bicroi h & \mapsto & \dS(\ome) \bicroi h
\end{array}
\end{equation}
and the canonical embeddings $\mathfrak{i}_{1}: \ome \in \du{H}^{\co} \mapsto \ome \bicroi 1 \in \D_{R}(H)$, $\mathfrak{i}_{2}: h \in H \mapsto \cou \bicroi h \in \D_{R}(H)$.

\subsection*{Left-left case}\label{sec:left-left-case}

Let $(X,\rhd,\alp)$ be a (left-left) Yetter--Drinfeld algebra over $H$. By Theorem~\ref{th:majid_duality}, $X$ is a left $\D_{M}(H)$-module algebra with action map
\[
(\ome^{\op} \bicroi h) \cdot x = (\p(\cdot,\ome)) \odo \id)(\alp(h \rhd x))
\]
for all $h \in H$, $\ome \in \du{H}$ and $x \in X$. Using the canonical isomorphism~\eqref{eq:canonical_DMR} between $\D_{M}(H)$ and $\D_{R}(H)$, hence $X$ yields a left $\D_{R}(H)$-module with module map
\[
(\ome \bicroi h) \rhd x = (\dS^{-1}(\ome)^{\op} \bicroi h) \cdot x = (\p(\cdot,\dS^{-1}(\ome)) \odo \id)(\alp(h \rhd x))
\]
for all $h \in H$, $\ome \in \du{H}$ and $x \in X$.

\begin{prop}
The last statement shows the equivalence
\[
(X,\rhd,\alp) \in {}^{H}_{H}\aYD \quad \sii \quad (X,\rhd,(\id \odo S^{-1})(\Sigma\alp)) \in {}_{H}\aYD^{H}
\]
\end{prop}

Consider the left action
\[
\begin{array}{lccccccc}
\rhd^{2} \; : & \D_{R}(H) \odo \D_{R}(H) & \odo & X \odo X & \to & X & \odo & X \\
& M \odo N & \odo & x \odo y & \mapsto & (M \rhd x) & \odo & (N \rhd y),
\end{array}
\]
then, for any element $P = \sum_{i} h_{i} \odo \theta_{i} \in H \odo \du{H}$, we denote by
\[
\begin{array}{lccc}
\lambda_{P} : & X \odo X & \to & X \odo X \\
& x \odo y & \mapsto & \sum_{i} (\mathfrak{i}_{2}(h_{i}) \rhd x) \odo (\mathfrak{i}_{1}(\theta_{i}) \rhd y)
\end{array}
\]
the linear map induced by the left action $\rhd^{2}$. Evidently $\lambda_{P} = \sum_{i} \lambda_{h_{i} \odo \theta_{i}}$ as linear maps. Moreover, if $P$ is invertible in $H \odo \du{H}$, thus $\lambda_{P}$ is an isomorphism with inverse given by $\lambda^{-1}_{P} = \lambda_{P^{-1}}$.

\begin{prop}
Let $U \in H \odo \du{H}$ be the canonical element for the canonical pairing of $H$. A (left-left) Yetter--Drinfeld algebra $(X,\rhd,\alp)$ is  braided commutative if and only if we have $\m_{X} = \m_{X}\Sigma_{X,X}\lambda_{U}$, or equivalently $\m_{X} = \m_{X}\lambda_{U^{-1}}\Sigma_{X,X}$.
\end{prop}
\pr
Because $U = U_{1} \odo U_{2} \in H \odo \du{H}$ is the canonical element for the canonical pairing of $H$, the equality $\m_{X} = \m_{X}\Sigma_{X,X}\lambda_{U}$ is equivalent to the condition
\begin{align*}
xy & = \m_{X}(x \odo y) = \m_{X}\Sigma_{X,X}\lambda_{U}(x \odo y) \\
& = \m_{X}\Sigma_{X,X}(((\cou \bicroi U_{1}) \rhd x) \odo ((U_{2} \bicroi 1) \rhd y)) = ((U_{2} \bicroi 1) \rhd y)((\cou \bicroi U_{1}) \rhd x) \\
& = (\p(y_{[-1]},\dS^{-1}(U_{2}))y_{[0]})(U_{1} \rhd x) = y_{[0]}(S^{-1}(y_{[-1]}) \rhd x)
\end{align*}
for every $x,y \in X$. Similarly, $\m_{X} = \m_{X}\lambda_{U^{-1}}\Sigma_{X,X}$ is equivalent to the condition $xy = (x_{[-1]} \rhd y)x_{[0]}$ for every $x,y \in X$. Here, we are using the Sweedler type leg notation $\alp(x) = x_{[-1]} \odo x_{[0]}$.
\fin

\subsection*{Left-right case}

Let $(X,\rhd,\alp)$ be a (left-right) Yetter--Drinfeld algebra over $H$. By Theorem~\ref{th:radford_duality}, $X$ is a left $\D_{R}(H)$-module algebra with action map
\[
(\ome \bicroi h) \rhd x = (\id \odo \p(\cdot,\ome))(\alp(h \rhd x))
\]
for all $h \in H$, $\ome \in \du{H}$ and $x \in X$.

\begin{prop}
Let $U \in H \odo \du{H}$ be the canonical element for the canonical pairing of $H$. A (left-right) Yetter--Drinfeld algebra $(X,\rhd,\alp)$ is braided commutative if and only if we have $\m_{X} = \m_{X}\Sigma_{X,X}\lambda_{U}$, or equivalently $\m_{X} = \m_{X}\lambda_{U^{-1}}\Sigma_{X,X}$.
\end{prop}
\pr
Because $U = U_{1} \odo U_{2} \in H \odo \du{H}$ is the canonical element for the canonical pairing of $H$, the equality $\m_{X} = \m_{X}\Sigma_{X,X}\lambda_{U}$ is equivalent to the condition
\begin{align*}
xy & = \m_{X}(x \odo y) = \m_{X}\Sigma_{X,X}\lambda_{U}(x \odo y) \\
& = \m_{X}\Sigma_{X,X}(((\cou \bicroi U_{1}) \rhd x) \odo ((U_{2} \bicroi 1) \rhd y)) = ((U_{2} \bicroi 1) \rhd y)((\cou \bicroi U_{1}) \rhd x) \\
& = (\p(y_{[1]},U_{2})y_{[0]})(U_{1} \rhd x) = y_{[0]}(y_{[1]} \rhd x)
\end{align*}
for every $x,y \in X$. Similarly, $\m_{X} = \m_{X}\lambda_{U^{-1}}\Sigma_{X,X}$ is equivalent to the condition $xy = (S(x_{[1]}) \rhd y)x_{[0]}$ for every $x,y \in X$. Here, we are using the Sweedler type leg notation $\alp(x) = x_{[0]} \odo x^{\op}_{[1]}$.
\fin

\subsection*{Right-right case}

Let $(X,\lhd,\alp)$ be a (right-right) Yetter--Drinfeld algebra over $H$. By Theorem~\ref{th:majid_duality_right}, $X$ is a right $\D_{M}(H)$-module algebra with action map
\[
x \cdot (\ome^{\op} \bicroi h) = (\id \odo \p(\cdot,\ome))(\alp(x)) \lhd h
\]
for all $h \in H$, $\ome \in \du{H}$ and $x \in X$. Using the canonical isomorphism~\eqref{eq:canonical_DMR} between $\D_{M}(H)$ and $\D_{R}(H)$, hence $X$ yields a right $\D_{R}(H)$-module with module map
\[
x \lhd (\ome \bicroi h) = x \cdot (\dS^{-1}(\ome)^{\op} \bicroi h) = (\id \odo \p(\cdot,\dS^{-1}(\ome)))(\alp(x)) \lhd h
\]
for all $h \in H$, $\ome \in \du{H}$ and $x \in X$.

\begin{prop}
The last statement shows the equivalence
\[
(X,\lhd,\alp) \in \aYD^{H}_{H} \quad \sii \quad (X,\lhd,(S^{-1} \odo \id)(\Sigma\alp)) \in {}^{H}\aYD_{H}
\]
\end{prop}

Consider the right action
\[
\begin{array}{lccccccc}
\lhd^{2} \; : & X \odo X & \odo & \D_{R}(H) \odo \D_{R}(H) & \to & X & \odo & X \\
& x \odo y & \odo & M \odo N & \mapsto & (x \lhd M) & \odo & (y \lhd N),
\end{array}
\]
then, for any element $P = \sum_{i} h_{i} \odo \theta_{i} \in H \odo \du{H}$, we denote by
\[
\begin{array}{lccc}
\rho_{P} : & X \odo X & \to & X \odo X \\
& x \odo y & \mapsto & \sum_{i} (x \lhd \mathfrak{i}_{2}(h_{i})) \odo (y \lhd \mathfrak{i}_{1}(\theta_{i}))
\end{array}
\]
the linear map induced by the left action $\lhd^{2}$. Evidently $\lambda_{P} = \sum_{i} \rho_{h_{i} \odo \theta_{i}}$ as linear maps. Moreover, if $P$ is invertible in $H \odo \du{H}$, thus $\rho_{P}$ is an isomorphism with inverse given by $\rho^{-1}_{P} = \rho_{P^{-1}}$.

\begin{prop}
Let $U \in H \odo \du{H}$ be the canonical element for the canonical pairing of $H$. A (right-right) Yetter--Drinfeld algebra $(X,\lhd,\alp)$ is braided commutative if and only if we have $\m_{X} = \m_{X}\Sigma_{X,X}\rho_{U^{-1}}$ or equivalently $\m_{X} = \m_{X}\rho_{U}\Sigma_{X,X}$.
\end{prop}
\pr
Because $U = U_{1} \odo U_{2} \in H \odo \du{H}$ is the canonical element for the canonical pairing of $H$ and $U^{-1}=(\id \odo \dS)(U)$, the equality $\m_{X} = \m_{X}\Sigma_{X,X}\rho_{U^{-1}}$ is equivalent to the condition
\begin{align*}
xy & = \m_{X}(x \odo y) = \m_{X}\Sigma_{X,X}\rho_{U^{-1}}(x \odo y) \\
& = \m_{X}\Sigma_{X,X}((x \odo y) \lhd^{2} (\cou \bicroi U_{1} \odo \dS(U_{2}) \bicroi 1)) = (y \lhd (\dS(U_{2}) \bicroi 1))(x \lhd (\cou \bicroi U_{1})) \\
& = (y_{[0]}\p(y_{[1]},U_{2}))(x \lhd U_{1}) = y_{[0]}(x \lhd y_{[1]})
\end{align*}
for every $x,y \in X$. Similarly, $\m_{X} = \m_{X}\rho_{U}\Sigma_{X,X}$ is equivalent to the condition $xy = (y \lhd S^{-1}(x_{[1]}))x_{[0]}$ for every $x,y \in X$. Here, we are using the Sweedler type leg notation $\alp(x) = x_{[0]} \odo x_{[1]}$.
\fin

\subsection*{Right-left case}

Let $(X,\lhd,\alp)$ be a (right-left) Yetter--Drinfeld algebra over $H$. By Theorem~\ref{th:radford_duality_right}, $X$ is a right $\D_{R}(H)$-module algebra with action map
\[
x \lhd (\ome \bicroi h) = (\p(\,\cdot\,,\ome) \odo \id)(\alp(x)) \lhd h 
\]
for all $h \in H$, $\ome \in \du{H}$ and $x \in X$.

\begin{prop}
Let $U \in H \odo \du{H}$ be the canonical element for the canonical pairing of $H$. A (right-left) Yetter--Drinfeld algebra $(X,\lhd,\alp)$ in braided commutative if and only if we have $\m_{X} = \m_{X}\Sigma_{X,X}\rho_{U^{-1}}$, or equivalently $\m_{X} = \m_{X}\rho_{U}\Sigma_{X,X}$.
\end{prop}
\pr
Because $U = U_{1} \odo U_{2} \in H \odo \du{H}$ is the canonical element for the canonical pairing of $H$ and $U^{-1}=(\id \odo \dS)(U)$, the equality $\m_{X} = \m_{X}\Sigma_{X,X}\rho_{U^{-1}}$ is equivalent to the condition
\begin{align*}
xy & = \m_{X}(x \odo y) = \m_{X}\Sigma_{X,X}\rho_{U^{-1}}(x \odo y) \\
& = \m_{X}\Sigma_{X,X}((x \odo y) \lhd^{2} (\cou \bicroi U_{1} \odo \dS(U_{2}) \bicroi 1)) = (y \lhd (\dS(U_{2}) \bicroi 1))(x \lhd (\cou \bicroi U_{1})) \\
& = (\p(y_{[-1]},\dS(U_{2}))y_{[0]})(x \lhd U_{1}) = y_{[0]}(x \lhd S(y_{[-1]}))
\end{align*}
for every $x,y \in X$. Similarly, $\m_{X} = \m_{X}\rho_{U}\Sigma_{X,X}$ is equivalent to the condition $xy = (y \lhd x_{[-1]})x_{[0]}$ for every $x,y \in X$. Here, we are using the Sweedler type leg notation $\alp(x) = x^{\op}_{[-1]} \odo x_{[0]}$.
\fin

\section{The categorical point of view}\label{sec:appendix_categorical}

\subsection*{The standard characterization}

Consider the category ${}^{H}\YD_{H}$ of (right-left) Yetter--Drinfeld modules over a Hopf algebra $H$. Recall the data of the category ${}^{H}\YD_{H}$
\begin{enumerate}[label=\textup{(\roman*)}]
\item The objects are the tuples $(X,\lhd:X \odo H \to X,\alp:X \to H \odo X)$ such that $\lhd : X \odo \du{H} \odo H \to X$ is an right action of $\D_{R}(H)$ on $X$.
\item The morphisms are given by $f: (X,\lhd,\alp) \to (X',\lhd',\alp')$ such that $f: X \to X'$ is a intertwiner for the right actions and for the left coactions, i.e. $f\lhd = \lhd'(f \odo \id)$ and $(\id \odo f)\alp = \alp' f$.
\end{enumerate}

\begin{defi}
Let $\mathcal{C}$ be a monoidal category. The center of $\mathcal{C}$ is given by the category $\mathcal{Z}(\mathcal{C})$: 
\begin{enumerate}[label=\textup{(\arabic*)}]
\item the {\em objects} are given by the pairs $(X,i)$, where $X$ is an object of $\mathcal{C}$, $i = \{i_{Y}\}_{Y \in \mathcal{C}}$ is a natural isomorphism (for each $Y \in \mathcal{C}$, $i_{Y}: X \oti Y \to Y \oti X$ is an isomorphism in $\mathcal{C}$) and for each $Y, Y' \in \mathcal{C}$, we have $i_{Y \oti Y'} = (\id_{Y} \oti i_{Y'})(i_{Y} \oti \id_{Y'})$.

\item the {\em morphisms} are given by
\[
\text{mor}_{\mathcal{Z}(\mathcal{C})}((X,i),(X',i')) = \{f \in \text{mor}_{\mathcal{C}}(X,X') \,:\, i'_{Y}(f \oti \id_{Y}) = (\id_{Y} \oti f) i_{Y} \text{ for each } Y \in \mathcal{C}\}.
\]
\item the {\em monoidal operation} is given by
\[
(X,i) \oti (X',i') := (X \oti X',\{(i_{Y} \oti \id_{X'})(\id_{X} \oti i'_{Y})\}_{Y \in \mathcal{C}}).
\]
\end{enumerate}
The category $\mathcal{Z}(\mathcal{C})$ yields a {\em braided monoidal category}.
\end{defi}

\begin{prop}
Let $H$ be a finite-dimensional Hopf algebra. We have
\[
{}^{H}\YD_{H} \iso \mathcal{Z}(\text{Mod}_{H}).
\]
\end{prop}

\subsection*{The ``only coaction'' characterization}

Let $H$ be a Hopf algebra. The category of left-left $(H,\du{H})$-Yetter--Drinfeld modules, denoted by $\YD^{ll}(H,\du{H})$, consists of the following data:
\begin{enumerate}[label=\textup{(\roman*)}]
\item The objects are the tuples $(X,\alp:X \to H \odo X,\beta:X \to \du{H} \odo X)$ such that $(\id \odo \beta)\alp : X \to H \odo \du{H} \odo X$ yields a module coaction of $\T_{R}(H)$ on $X$.
\item The morphisms $f: (X,\alp,\beta) \to (X',\alp',\beta')$ are given by maps $f: X \to X'$ such that $(\id \odo f)\alp = \alp' f$ and $(\id \odo f)\beta = \beta' f$ or equivalently such that $(\id \odo \id \odo f)(\id \odo \beta)\alp = (\id \odo \beta')\alp' f$.
\end{enumerate}

Given two objects $(X,\alp,\beta)$ and $(X',\alp',\beta')$ in $\YD^{ll}(H,\du{H})$, the {\em braided product of $X$ and $Y$}, denoted by $X \bt Y$, is defined as the vector space $X \odo Y$ endowed with the left $H$-comodule map $\alp \bt \alp' = (\m_{H} \odo \id \odo \id)\Sigma_{12}\Sigma_{23}(\alp \odo \alp')$ and the left $\du{H}$-comodule map $\beta \bt \beta' = (\m_{\du{H}} \odo \id \odo \id)\Sigma_{23}(\beta \odo \beta')$. The tuple $(X \bt Y,\alp^{\bt},\beta^{\bt})$ yields an object of $\YD^{ll}(H,\du{H})$. In general, we have the bifunctor
\[
\begin{array}{lccccc}
\bt: & \YD^{ll}(H,\du{H}) & \times & \YD^{ll}(H,\du{H}) & \to & \YD^{ll}(H,\du{H}) \\
& (X,\alp,\beta) & \times & (X,\alp',\beta') & \mapsto & (X \bt Y,\alp \bt \alp',\beta \bt \beta')
\end{array}.
\]
The tuple $(\ku,\alp_{tr},\beta_{tr})$ is called {\em the unit of $\YD^{ll}(H,\du{H})$} with relation to the braided product $\bt$, because
\[
(\ku,\alp_{tr},\beta_{tr}) \bt (X,\alp,\beta) = (X,\alp,\beta) \quad \text{ and } \quad (X,\alp,\beta) \bt (\ku,\alp_{tr},\beta_{tr}) = (X,\alp,\beta).
\]
for every object $(X,\alp,\beta) \in \YD^{ll}(H,\du{H})$.

\begin{defi}
An object $(X,\alp,\beta) \in \YD^{ll}(H,\du{H})$ is called {\em an algebra in $\YD^{ll}(H,\du{H})$}, if there is morphisms in $\YD^{ll}(H,\du{H})$, $\m_{X}: X \bt X \to X$ and $\mu_{X}: \ku \to X$, such that
\[
\m_{X}(\id_{X} \bt \mu_{X}) = \m_{X}(\mu_{X} \bt \id_{X}) = \id_{X} \quad \text{ and } \quad \m_{X}(\id_{X} \bt \m_{X}) = \m_{X}(\m_{X} \bt \id_{X}).
\]
\end{defi}

\begin{prop}
The subcategory formed by the algebra objects of $\YD^{ll}(H,\du{H})$ is $\aYD^{ll}(H,\du{H})$. In others words, $(X,\alp,\beta)$ is an algebra object in $\YD^{ll}(H,\du{H})$ if and only if $(X,\alp,\beta)$ is a left-left $(H,\du{H})$-Yetter--Drinfeld algebra.
\end{prop}
\pr
Straightforward.
\fin

\end{document}